\documentclass[a4paper, 10pt]{article}
\usepackage[pdfencoding=auto,psdextra]{hyperref}
\usepackage{tikz}
\usepackage{tikz-cd}

\usepackage{latexsym}
\usepackage{verbatim}
\usepackage{enumitem}
\usepackage{amsmath,amsthm,amssymb}
\usepackage[mathscr]{euscript}
\usepackage{yfonts}
\usepackage{stmaryrd}
\usepackage{multicol}
\usepackage{bm}
\usepackage[all]{xy}
\numberwithin{equation}{section}
\newtheorem{thm}{Theorem}[section]
\newtheorem{prop}[thm]{Proposition}
\newtheorem{lem}[thm]{Lemma}
\newtheorem{cor}[thm]{Corollary}
\newtheorem{claim}{Claim}{\bf}{\it}

\newtheorem{fthm}{Theorem}{\bf}{\it}
{\bf}{\it}
\newtheorem{fcor}[fthm]{Corollary}{\bf}{\it}
{\bf}{\it}
{\bf}{\it}

\theoremstyle{definition}
\newtheorem{defn}[thm]{Definition}

{\bf}{\rm}

\theoremstyle{remark}

\newtheorem{rem}[thm]{Remark}
{\bf}{\it}

\newtheorem{definition and corollary}[thm]{Definition and Corollary}

\newtheorem{fex}[fthm]{Example}{\it}{\rm}

\DeclareMathOperator{\Hom}{Hom}
\DeclareMathOperator{\shom}{hom}
\DeclareMathOperator{\Ext}{Ext}
\DeclareMathOperator{\ext}{ext}

\DeclareMathOperator*{\gch}{\mathrm{gch}}

\DeclareMathOperator*{\gdim}{\mathrm{gdim}}

\DeclareMathOperator{\ch}{ch}
\DeclareMathOperator{\gr}{gr}

\DeclareMathOperator{\hd}{\mathsf{hd}}
\DeclareMathOperator{\soc}{\mathsf{soc}}

\newcommand{\fB}{\mathfrak{B}}
\newcommand{\fC}{\mathfrak{C}}
\newcommand{\al}{\alpha}
\newcommand{\af}{\mathrm{af}}

\newcommand{\C}{{\mathbb C}}

\newcommand{\bq}{{\mathsf q}}

\newcommand{\bD}{\mathbb{D}}

\newcommand{\la}{\lambda}
\newcommand{\La}{\Lambda}

\newcommand{\bd}{\mathrm{bdd}}

\newcommand{\g}{\mathfrak{g}}
\newcommand{\tg}{\widetilde{\mathfrak{g}}}
\newcommand{\gb}{\mathfrak{b}}
\newcommand{\tb}{\widetilde{\mathfrak{b}}}

\newcommand{\h}{\mathfrak{h}}

\newcommand{\wth}{\widetilde{\mathfrak{h}}}

\newcommand{\tI}{\mathtt{I}}

\newcommand{\gn}{\mathfrak{n}}
\newcommand{\tn}{\widetilde{\mathfrak{n}}}

\newcommand{\tp}{\widetilde{\mathfrak{p}}}

\newcommand{\bv}{\mathbf{v}}
\newcommand{\wg}{\widehat{\mathfrak{g}}}

\newcommand{\Q}{\mathbb{Q}}
\newcommand{\R}{\mathbb{R}}
\newcommand{\bW}{\mathbb{W}}

\newcommand{\Z}{\mathbb{Z}}

\newcommand{\Gm}{\mathbb G_m}

\title{Higher level BGG reciprocity for current algebras}
\author{Syu \textsc{Kato}\footnote{Department of Mathematics, Kyoto University, Oiwake Kita-Shirakawa Sakyo Kyoto 606-8502 JAPAN}\footnote{\tt{E-mail:syuchan@math.kyoto-u.ac.jp, ORCID:https://orcid.org/0000-0002-4593-3205}}}

\begin{document}
\maketitle
\begin{abstract}
We establish a higher-level analogue of the Bernstein--Gelfand--Gelfand (BGG)
reciprocity for twisted current algebras at every positive integer level,
extending the level-one results of Bennett et al. We show that these reciprocity
relations form a hierarchy linking adjacent levels, thereby providing a
systematic framework to relate representations across arbitrary levels by
iteration.

Our formulation of the BGG reciprocity naturally includes infinite-dimensional
(thick) Demazure modules. Consequently, our main theorem provides the first
effective means to analyze the module-theoretic structure of thick Demazure
modules at any level. As an application, we realize level-restricted generalized
Kostka polynomials as branching polynomials in the language of symmetric
polynomials. The resulting level-dependent picture also naturally incorporates
theta functions and modular forms. Furthermore, we prove that every
finite-dimensional Demazure module admits a filtration by Demazure modules of
the underlying simple Lie algebra, resolving a long-standing speculation in
non-simply-laced types.
\end{abstract}
\setcounter{tocdepth}{2}
\tableofcontents

\section*{Introduction}
An affine Kac-Moody algebra $\tg$ over $\mathbb{C}$ has a non-negative part $\tg_{\ge 0}$, which contains a simple Lie algebra $\g$ whose Dynkin diagram is obtained by removing the zeroth vertex from the Dynkin diagram of $\tg$ (\cite{Kac}). 
Let $P$ denote the set of integral weights of $\g$, and let $P^+$ be the set of dominant integral weights of $\g$. 
Following the work of Chari-Pressley~\cite{CP01} on $\tg$, Chari-Loktev~\cite{CL06} introduced the local Weyl modules of $\tg_{\ge 0}$, and Chari-Fourier-Khandai~\cite{CFK} introduced the global Weyl modules of $\tg_{\ge 0}$, denoted by $W_{\la}$ and $\bW_{\la}$ $(\la \in P^+)$, respectively.
Let $V_{\la}$ ($\la \in P^+$) be the irreducible finite-dimensional $\g$-module of highest weight $\la$. 
Bennett, Berenstein, Chari, Ion, Khoroshkin, Loktev, and Manning~\cite{BCM,BBCKL,CI15} showed the following Ext-orthogonality:
\begin{equation}
\mathrm{Ext}^{i} ( \bW_{\la}, W_{\mu}^{*}) \cong 
\begin{cases} 
\mathbb{C} & \text{if } i=0 \text{ and } V_{\la} \cong V_{\mu}^* \\ 
0 & \text{otherwise}
\end{cases}.
\label{Intro-Extorth}
\end{equation}
This relation is tightly connected to the BGG reciprocity for $\tg_{\ge 0}$. 
The graded characters of $W_{\la}$ and $\bW_{\la}$ are proportional to each other, with the proportionality factor being an invertible element of $\mathbb{Z}[\![q]\!]$. They are realized as specializations of a Macdonald polynomial (\cite{Mac95B}) at $t=0$.
In addition, the numerical counterpart of~\eqref{Intro-Extorth} is precisely the orthogonality relation with respect to the Hall inner product.

Replacing the subalgebra $\tg_{\ge 0}$ with the Iwahori subalgebra of $\tg$ (i.e., the upper triangular part in the sense of Kac-Moody algebras) yields another class of modules, $D_{\la}$ and $\bD_{\la}$, indexed by $\la \in P$. The modules $D_\la$ are the level-one Demazure modules, and their graded characters are shown to satisfy
\[
    \gch D_{\la}= E_{\la}(q, 0)
\]
by Sanderson~\cite{San00} and Ion~\cite{Ion03}, where $E_{\la}(q,t)$ is a non-symmetric Macdonald polynomial (\cite{Che95}). The modules $\bD_{\la}$ were introduced in~\cite{FKM,CK18}, and their graded characters are given by
\[
    \gch \bD_{\la} = c_{\la} \cdot E_{-\la}^{\dag}(q^{-1}, \infty),
\]
where $E_{\la}^{\dag}(q, t)$ is the character conjugate of $E_{\la}(q,t)$, and $c_{\la}$ is an invertible element in $\mathbb{Z}[\![q]\!]$. These families of modules possess $\mathrm{Ext}$-orthogonality properties analogous to~\eqref{Intro-Extorth}.

From this perspective, the relation~\eqref{Intro-Extorth} may be viewed as an
affine analogue of the duality between Demazure modules, extending the
classical result of van der Kallen~\cite{vdK89} (see also
\cite{Mat89c,Pol89}). This suggests that one should seek a similar
generalization for arbitrary Demazure modules in integrable highest-weight
modules over Kac--Moody algebras. However, such an extension is not
straightforward in the infinite-dimensional setting. The main difficulty is
that, for a general Kac--Moody algebra, there are two natural versions of
Demazure modules, which essentially coincide only in finite type. Outside
finite type, a proper formulation of~\eqref{Intro-Extorth} necessarily
involves both versions: one finite-dimensional and the other
infinite-dimensional.

The finite-dimensional Demazure modules, which we refer to as \emph{thin Demazure modules}, are obtained by successive applications of Demazure functors starting from one-dimensional modules.
This fact is quite helpful in their analysis (see Remark~\ref{rem:thin}).
(These are the modules usually called ``Demazure modules'' in the literature; see, e.g.,~\cite{Kum02}.)
Infinite-dimensional Demazure modules, which we refer to as \emph{thick Demazure modules}, are much more intractable.
Unlike their thin counterparts, they are defined in terms of irreducible highest-weight integrable representations and cannot be constructed from finite-dimensional objects via successive applications of Demazure functors.
In fact, because Demazure functors lose information, there have been practically no effective means to analyze the precise module-theoretic structure of thick Demazure modules.

Integrable highest-weight representations of an affine Lie algebra come equipped with a level, and this level provides an additional parameter for organizing the corresponding $\g$-stable Demazure modules, viewed as a family indexed by $P^+$.
For simply-laced $\g$, the branching property of thin Demazure modules
is known by~\cite{Jos06,Nao12}. As the level varies, it implies that we have a family of symmetric polynomials interpolating between $\{\gch W_\la\}_{\la \in P^+}$, corresponding to Macdonald polynomials at $t=0$, and $\{\gch V_\la\}_{\la \in P^+}$, corresponding to Schur polynomials.
Furthermore, the notion of level also leads to the combinatorial study of level-restricted generalized Kostka polynomials, and some explicit descriptions are known (\cite{HKOTY,HKKOTY,SW99,SS01,KL17}).
Thus, at least since the appearance of the foundational works \cite{CF13,CI15,FKM,CK18}, the extension of the reciprocity in \eqref{Intro-Extorth} to the higher-level case has been a natural question.

In this paper, we build upon the framework developed in~\cite{CK18,Kat18b} to address the aforementioned difficulties. 
We begin by considering the category $\mathcal{C}$, which consists of $[\tg_{\ge 0},\tg_{\ge 0}]$-modules that are semisimple when viewed as $\g$-modules. 
In this setting, we disregard the central charge and the grading.

\begin{fthm}[$\doteq$ Theorems~\ref{thm:dual},~\ref{thm:critfilt}, and~\ref{thm:branch}]\label{f-filt}
For each $k \in \mathbb{Z}_{>0}$ and $\la \in P^+$, there exist two modules, $W_\la^{(k)}$ and $\bW_\la^{(k)}$, satisfying the following properties:
\begin{enumerate}
    \item The modules $W_\la^{(k)}$ and $\bW_\la^{(k)}$ are $\tg_{\ge 0}$-modules whose heads are isomorphic to $V_\la$.

    \item The module $W_\la^{(k)}$ is a thin Demazure module of an irreducible integrable highest-weight module of $\tg$ of level $k$.

    \item The following Ext-orthogonality holds in the category $\mathcal{C}$:
    \begin{equation}
    \mathrm{Ext}^{i}_{\mathcal C} ( \bW_{\la}^{(k)}, (W_{\mu}^{(k)})^{*}) \cong 
    \begin{cases} 
    \mathbb{C} & \text{if } i=0 \text{ and } V_{\la} \cong V_{\mu}^* \\ 
    0 & \text{otherwise}
    \end{cases}.
    \label{Intro-Extorth2}
    \end{equation}

    \item The module $W_\la^{(k)}$ has a filtration whose factors are of the form $W_{\mu}^{(k+1)}$ for $\mu \in P^+$, and the module $\bW_\la^{(k+1)}$ has a filtration whose factors are of the form $\bW_\mu^{(k)}$ for $\mu \in P^+$.

    \item For each $\mu \in P^+$, the following multiplicity formulas hold:
    \begin{equation*}
    (P_\la : \bW_\mu^{(k)}) = [ W_\mu^{(k)} : V_{\la} ] \quad \text{and} \quad (\bW_\la^{(k+1)} : \bW_\mu^{(k)}) = (W_\mu^{(k)} : W_{\la}^{(k+1)}),
    \end{equation*}
    where $(M:N)$ denotes the multiplicity of $N$ in a given filtration of $M$, and $P_\la$ is the projective cover of $V_\la$ in $\mathcal{C}$.
\end{enumerate}
\end{fthm}

Taken together, these results show that the level-one orthogonality in~\eqref{Intro-Extorth} extends not merely to an isolated higher-level statement, but to a hierarchy of reciprocity relations linking adjacent levels.

The $k=1$ case of~\eqref{Intro-Extorth2} is~\eqref{Intro-Extorth}, and its
non-symmetric counterpart is contained in~\cite{CK18}; we also present the
non-symmetric counterparts of Theorem~\ref{f-filt} in \S\ref{sec:statement}.
For $k \ge 2$, the only previously known case of Theorem~\ref{f-filt}, apart
from the first two items, is the existence of a filtration of $W_\la^{(k)}$ by
$\{W_\mu^{(k+1)}\}_{\mu \in P^+}$ when $\g$ is of type $\mathsf{ADE}$
(Joseph--Naoi~\cite{Jos06,Nao12}; see also~\cite{CSSW}). In particular,
Theorem~\ref{f-filt} is new even for $\g = \mathfrak{sl}(2)$ in its full
generality. Let us add that our proof here (and consequently our proof of
Corollary~\ref{fcor:cbr} below) is independent of that of Joseph--Naoi.

Theorem~\ref{f-filt} asserts that $\mathcal{C}$ possesses a structure similar to
that of highest-weight categories~\cite{CPS88,BS24}. In fact, if we also
incorporate the grading into the structure, then we deduce that $W_\la^{(k)}$
and $\bigl( \bW_\la^{(k)} \bigr)^{\vee}$ can be naturally understood as
standard/costandard modules in a certain graded category
(Corollary~\ref{cor:Whw} and Remark~\ref{rem:hw}). From this perspective, our
relation~\eqref{Intro-Extorth2} exactly matches the $\Ext$-vanishing expected
from the framework of highest-weight theory: Theorem~\ref{f-filt}(3) is the
higher-level analogue of~\eqref{Intro-Extorth}, while Theorem~\ref{f-filt}(5)
plays the role of a higher-level BGG reciprocity (\cite{CI15}).

We can compute the multiplicity $(W_{\mu}^{(k)}: W_{\la}^{(k+1)})$ by counting
highest-weight elements in suitable tensor products of crystals
(Joseph~\cite{Jos03}); there is also another combinatorial formula for
$\g = \mathfrak{sl}(2)$ by Biswal--Kus~\cite{BK21}. Such descriptions, however,
are confined to the thin side: the description in \cite[\S3.3]{Jos03} {\it
cannot} capture $(\bW_{\la}^{(k+1)}: \bW_{\mu}^{(k)})$, as the modules
$\bW^{(k)}_\la$ do not contain highest-weight vectors. In this sense,
Theorem~\ref{f-filt} yields the first effective means to analyze the structure
of the infinite-dimensional Demazure modules studied
in~\cite{KS09,CF13,Kat18b,CK18} for an arbitrary level.

This level-dependent picture has a particularly transparent asymptotic form: for
a fixed $\la \in P^+$ and for $k \gg 0$ we have $W_{\la}^{(k)} = V_\la$, and
$\bW_\la^{(k)}$ is isomorphic to a level-$(1{-}k)$ integrable lowest weight
module of $\tg$ (Lemma~\ref{lem:idDL}). Thus, Theorem~\ref{f-filt} enlarges the
scope of the theory of symmetric polynomials (\cite{Mac95}) to include theta
functions and modular forms (Remark~\ref{rem:theta}).

The reciprocity in Theorem~\ref{f-filt} also leads to concrete structural consequences for finite-dimensional Demazure modules.
For instance, the following result resolves a long-standing speculation, which generalizes a result from~\cite{Jos06} for simply-laced types:

\begin{fthm}[$\doteq$ Corollary~\ref{cor:kumar}]\label{fcor:cbr}
Every finite-dimensional Demazure module of $\tg$ admits a filtration whose successive quotients are isomorphic to Demazure modules of $\g$ by restriction.
\end{fthm}

As a further consequence of our analysis, we recover the level-$k$-restricted
generalized Kostka polynomials in the sense of~\cite{KL17} as a subfamily of the
branching polynomials arising from Theorem~\ref{f-filt}. These are the
$q$-polynomials attached to tensor products of Kirillov--Reshetikhin crystals via
the corresponding one-dimensional sums; see~\eqref{eqn:Xkdef}, Corollary~\ref{cor:kostka}, and
Remark~\ref{rem:kostka} for the precise identification. In particular, the graded
characters of our modules provide a natural realization of these polynomials:

\begin{fcor}[$\doteq$ Corollary~\ref{cor:KR} and Remark~\ref{rem:kostka}]
Let $k \in \mathbb{Z}_{>0}$, and let $B$ be a tensor product of
Kirillov--Reshetikhin crystals of level at most $k$ (see~\cite[\S3.3]{HKOTT}).
Suppose that its character is given by
\[
    \sum_{\la \in P^+} X_\la(q)\cdot \ch V_\la,
    \qquad X_\la(q)\in\mathbb{Z}[q].
\]
We expand this character in terms of the graded characters of
$W_\la^{(k+1)}$ as
\[
    \sum_{\la \in P^+} X_\la(q)\cdot \ch V_\la
    =
    \sum_{\la \in P^+} X^{(k)}_\la(q)\cdot \gch W^{(k+1)}_\la,
    \qquad X^{(k)}_\la(q)\in\mathbb{Z}[q].
\]
Then, for dominant weights $\la$ satisfying the level condition
$\langle \vartheta^\vee, \la \rangle \le k$, the coefficients
$X^{(k)}_\la(q)$ coincide with the level-$k$-restricted generalized Kostka
polynomials associated with $B$.
\end{fcor}

Here, the level shift $k \mapsto (k+1)$ appearing in $W^{(k+1)}_\la$ may be explained by the following branching formula:

\begin{fcor}[$\doteq$ Corollary~\ref{cor:num}]\label{fcor:num}
Let $k \in \mathbb{Z}_{>0}$. For each $\la, \mu \in P^+$ such that $\langle \vartheta^\vee, \la \rangle \le k$, we have
\begin{equation*}
(W_{\mu}^{(k)} : W_{\la}^{(k+1)}) = \begin{cases} 1 & \text{if } \mu \uparrow_k \la, \\ 0 & \text{otherwise}, \end{cases}
\end{equation*}
where $\mu\uparrow_k \la$ means that $\mu$ is an extremal weight of the level-$k$ integrable highest-weight module of $\tg$ whose highest weight is $\la$.
\end{fcor}

\begin{fex}
Assume that $\g = \mathfrak{sl}(2)$, whose fundamental weight is denoted by $\varpi$. The first few branching rules are:
\begin{align*}
\gch W_{2 \varpi}^{(1)} &= \gch W_{2 \varpi}^{(2)} + q \, \gch W_{0}^{(2)}, & \gch W_{3 \varpi}^{(1)} &= \gch W_{3 \varpi}^{(2)} + q^2 \, \gch W_{\varpi}^{(2)}, \\
\gch W_{3 \varpi}^{(2)} &= \gch W_{3 \varpi}^{(3)} + q \, \gch W_{\varpi}^{(3)}, & \gch W_{4 \varpi}^{(2)} &= \gch W_{4 \varpi}^{(3)} + q^2 \, \gch W_{0}^{(3)}, \\
\gch W_{4 \varpi}^{(3)} &= \gch W_{4 \varpi}^{(4)} + q \, \gch W_{2 \varpi}^{(4)}, & &\cdots.
\end{align*}
By iterating these rules, we obtain expansions such as
\begin{align*}
\gch W_{4 \varpi}^{(1)} &= \gch W_{4 \varpi}^{(2)} + (q^2 + q^3) \gch W_{2 \varpi}^{(2)} + q^4 \gch W_{0}^{(2)}, \\
\gch W_{6 \varpi}^{(1)} &= \gch W_{6 \varpi}^{(2)} + (q^3 + q^4 + q^5) \gch W_{4 \varpi}^{(2)} + (q^6 + q^7 + q^8) \gch W_{2 \varpi}^{(2)} + q^9 \gch W_{0}^{(2)}.
\end{align*}
The appearance of the highest degree terms $q^4$ and $q^9$ in these expressions implies that the module $W_0^{(2)}$ appears in the filtrations of $W_{4\varpi}^{(1)}$ and $W_{6 \varpi}^{(1)}$ only at their socles. Theorem~\ref{f-filt}(5) then asserts the dual relations:
\begin{align*}
\gch \bW_{0}^{(2)} &= \gch \bW_{0}^{(1)} + q \, \gch \bW_{2\varpi}^{(1)} + q^4 \gch \bW_{4\varpi}^{(1)} + q^9 \gch \bW_{6\varpi}^{(1)} + \cdots, \\
\gch \bW_{2 \varpi}^{(2)} &= \gch \bW_{2 \varpi}^{(1)} + (q^2 + q^3) \gch \bW_{4\varpi}^{(1)} + (q^6 + q^7 + q^8) \gch \bW_{6\varpi}^{(1)} + \cdots.
\end{align*}
\end{fex}

\medskip

\noindent\textit{An overview of our proof of Theorem~\ref{f-filt}.}

\smallskip

Theorem~\ref{f-filt} is the central result of the paper, so we briefly explain the structure of its proof. The base of the induction, the $k=1$ case, was established for $\g$ of certain types in~\cite{FKM}, and for general $\g$ in~\cite{CK18}. 
The main argument proceeds by induction on $k$, relying on the following two statements, which form the technical core of Theorem~\ref{f-filt}:
\begin{itemize}
    \item[$(\clubsuit)_k$] The modules $W^{(k)}_\la$ and $\bW^{(k)}_\mu$ and their variants satisfy the $\Ext$-orthogonality condition~\eqref{Intro-Extorth2}.
    \item[$(\spadesuit)_k$] An $\Ext$-criterion for a module to admit a filtration by $\{ W^{(k)}_\la \}_{\la \in P^+}$ or $\{ \bW^{(k)}_\mu \}_{\mu \in P^+}$, and their variants.
\end{itemize}
The precise statement of $(\clubsuit)_k$ is given by Theorem~\ref{thm:dual} at level $k$, and that of $(\spadesuit)_k$ is given by Theorem~\ref{thm:critfilt} at level $k$.  
The key idea behind the implication $(\clubsuit)_k \Rightarrow (\spadesuit)_k$ is an adaptation of the standard machinery of highest-weight categories (\cite{CPS88,BS24}), stated as Proposition~\ref{prop:Dhw}.  
It therefore remains to prove $(\clubsuit)_k$ using the induction hypothesis $(\clubsuit)_{k-1}$ and $(\spadesuit)_{k-1}$.

This is done in several steps, based on the following three ingredients:
\begin{itemize}[leftmargin=*, align=left]
    \item[\textit{Reduction.}] A structural property of special Demazure modules that holds particularly in the setting of twisted affinizations (Proposition~\ref{prop:Demext}).
    \item[\textit{Lifting.}] Morphisms of certain types between modules admit a lifting to a morphism from a larger module obtained by applying a Demazure functor (Theorem~\ref{thm:indsurj} and Corollary~\ref{cor:indsurj}).
    \item[\textit{Vanishing.}] $\Ext^1$-vanishing between modules arising from level $k$ and level $(k-1)$ (Theorem~\ref{thm:Lext}).
\end{itemize}
Here, \textit{Reduction} is a property that fails to hold for a general affine Lie algebra $\tg$. Consequently, if one could find a proper replacement for this property, it might be possible to generalize our results to the remaining types of affine Lie algebras. 
The proof of \textit{Lifting} uses a careful analysis, adapted from~\cite{Jos85}, of so-called ``string modules''. 
The proof of \textit{Reduction} is a case-by-case analysis using the structure of affine root systems. 
The proof of \textit{Vanishing} uses, in addition to \textit{Lifting} and \textit{Reduction}, a characterization of (global) Weyl modules (\cite{CK18,Kat18}), a rough estimate of the size of \(W^{(k)}_\la\) (Lemma~\ref{lem:we}), and the adjoint property of Demazure functors (\cite{FKM,Mat89c}).
The proofs of \textit{Reduction} and \textit{Lifting} do not, by themselves, require the induction hypotheses $(\clubsuit)$ or $(\spadesuit)$.

Assuming \((\spadesuit)_{k-1}\), these three ingredients yield a filtration of \(\bW^{(k)}_\la\) by modules of the form \(\{\bW^{(k-1)}_\mu\}_\mu\) (Proposition~\ref{prop:passage}).
This, in turn, leads to the vanishing theorem (Theorem~\ref{thm:orth}) needed to carry out the proof of $(\clubsuit)_{k}$.

\medskip

\noindent\textit{Organization of the paper.}

\smallskip

The remainder of the paper is organized as follows.
First, in Sections~\ref{sec:prelim} and \ref{sec:module}, we provide the necessary background material, including the notation used throughout the paper and foundational module-theoretic results.

Sections~\ref{sec:tmodule} and~\ref{sec:lift} establish two of the three crucial ingredients mentioned in the overview.
Section~\ref{sec:tmodule} establishes a structural property of special Demazure modules, namely the \textit{Reduction} property.
Section~\ref{sec:lift} proves our key theorem on the lifting of morphisms, referred to as the \textit{Lifting} property.

Section~\ref{sec:statement} gives the precise formulation of our main theorems (Theorems~\ref{thm:dual} and~\ref{thm:critfilt}), including a refined version of Theorem~\ref{f-filt}(3) and a discussion of the highest-weight structure (\S\ref{subsec:hw}).

Sections~\ref{sec:hw},~\ref{sec:ext1}, and~\ref{sec:filt} form the induction loop in the proofs of the main theorems.
Section~\ref{sec:hw} characterizes the defining equations of the modules $W^{(k)}_\la$ and $\bW^{(k)}_\mu$ and their variants.
Section~\ref{sec:ext1} utilizes the results on Lifting and Reduction to derive the key $\Ext^1$-vanishing theorem, which is the \textit{Vanishing} property.
Section~\ref{sec:filt} then combines the induction hypothesis with the Vanishing results from Section~\ref{sec:ext1} to establish the crucial filtration properties of our modules, a decisive step for advancing the induction.

Finally, Section~\ref{sec:proof} assembles all these ingredients to complete the inductive argument and finalize the proofs of the main theorems.

The subsequent sections are devoted to applications.
Section~\ref{sec:branch} proves the branching rules (Theorem~\ref{f-filt}(4)(5)) and
Section~\ref{sec:app} connects our results with level-restricted Kostka polynomials.

In Appendix~A, we exhibit higher-level analogues of the $q$-Cauchy kernel identities, which were studied in the literature for the level-one case in rather recent works (\cite{FMO23,FKMO}). These numerical identities highlight the richness of the theory developed here, suggesting deep connections with the emerging study of higher-level $q$-series.

\medskip

We end this introduction with a heuristic remark. While we have expressed our results by using ungraded $\Hom$s and $\Ext$s in the Introduction, we mainly employ the graded ones, denoted by $\hom$s and $\ext$s (\S\ref{subsubsec:ghom}). They are isomorphic as vector spaces in the arguments of this paper, but we prefer to keep the gradings in mind, since the latter retain crucial information.

\section{Preliminaries}\label{sec:prelim}

We work over the field $\C$ of complex numbers.
In this paper, all graded vector spaces are assumed to be $\Z$-graded. For such a vector space $V$, we denote its degree-$m$ component by $V_m$ and set
\[ \gdim V := \sum_{m \in \Z} q^{m} \dim V_m. \]

For background in this section, see Kac~\cite{Kac} and Kumar~\cite{Kum02}.

\subsection{Lie algebras and their root systems}
In this subsection, we fix our conventions for the finite and twisted affine root data.

Let $\g$ be a simple Lie algebra over $\C$, with a fixed Cartan subalgebra $\h$ and a Borel subalgebra $\gb \supset \h$. We set $\gn := [\gb, \gb]$. We have the root space decompositions
\[ \gn = \bigoplus_{\al \in \Delta^+} \g_\al, \quad \gn^- := \bigoplus_{\al \in \Delta^+} \g_{-\al}, \quad \text{and} \quad \g = \gn^- \oplus \h \oplus \gn, \]
where $\Delta^+ \subset \h^*$ is the set of positive roots, $\Delta := \Delta^+ \sqcup (- \Delta^+)$, and $\g_{\al}$ is the root space for a root $\al \in \Delta$. Let $r$ be the lacing number of $\g$, defined as the ratio of the squared length of a long root to that of a short root. Thus, $r = 1$ for types $\mathsf{A,D,E}$, $r = 2$ for types $\mathsf{B,C,F}$, and $r = 3$ for type $\mathsf{G}_2$.

\subsubsection{Finite root systems and weights}
Let $\tI$ be the set of vertices of the Dynkin diagram of $\g$. Let $\Pi = \{\al_i\}_{i \in \tI} \subset \Delta^+$ be the set of simple roots, and let $\vartheta \in \Delta^+$ be the highest short root. We set $\Delta_s := W \vartheta \subset \Delta$, and $\Delta_l := \Delta \setminus \Delta_s$.
For each $\al \in \Delta$, we have its coroot $\al^{\vee} \in \h$. We set $\Pi^{\vee} := \{ \al_i^{\vee} \}_{i \in \tI}$.

The Weyl group $W$ of $\g$ is a subgroup of $\mathop{GL}( \h )$ generated by reflections $s_\al$ with respect to $\al^{\vee}$ for all $\al \in \Delta^+$. It is known that $W$ is generated by the simple reflections $\{s_i\}_{i \in \tI}$, where $s_i := s_{\al_i}$. The length function with respect to simple reflections is denoted by $\ell$.
Let $w_0 \in W$ be the longest element. We have $-w_0 \vartheta = \vartheta$, since $\vartheta$ is the unique dominant short root.

Let $Q := \sum_{i \in \tI} \Z \al_i \subset \h^*$ be the root lattice, and let $Q_+ := \sum_{i \in \tI} \Z_{\ge 0} \al_i \subset Q$ be its positive submonoid.
We set $Q^{\vee} := \sum_{i \in \tI} \Z \al_i^{\vee} \subset \h$ to be the coroot lattice.
The set of integral weights $P$ and the set of dominant integral weights $P^+$ are defined as follows:
\begin{align*}
    P &:= \{ \la \in \h^* \mid \langle \la, \al_i^\vee \rangle \in \Z \text{ for all } i \in \tI \}, \\
    P^+ &:= \{ \la \in P \mid \langle \la, \al_i^\vee \rangle \in \Z_{\ge 0} \text{ for all } i \in \tI \}.
\end{align*}
Let $\{\varpi_i\}_{i \in \tI} \subset P^+$ be the set of fundamental weights of $\g$, defined by the condition $\langle \varpi_j, \al_i^\vee \rangle = \delta_{ij}$ for all $i,j \in \tI$.

Each $\la \in P^+$ determines an irreducible finite-dimensional $\g$-module $V_\la$ with highest weight $\la$. We fix a highest-weight vector $\bv_\la$ in $V_\la$, unique up to scalar, which is $\gn$-invariant.

\subsubsection{Twisted affine root systems and their Weyl groups}

The twisted affine root system associated with $\g$ is
\[
\Delta_{\af} := \{ \alpha + n\delta \mid \alpha \in \Delta_s,\ n \in \Z \} \sqcup \{ \beta + nr\delta \mid \beta \in \Delta_l,\ n \in \Z \},
\]
where $\delta$ denotes the primitive imaginary root.
We fix the standard set of positive affine roots $\Delta_{\af}^{+} \subset \Delta_{\af}$, so that $\Delta^+ \subset \Delta_{\af}^{+}$. We define the simple affine roots by $\Pi_{\af} := \Pi \cup \{ \alpha_0 \}$, where $\alpha_0 := - \vartheta + \delta$, and index them by $\tI_{\af} := \tI \cup \{ 0 \}$.

The affine Weyl group is the semidirect product $W_{\af} := W \ltimes Q$. It is a reflection group generated by $\{s_i \mid i \in \tI_{\af} \}$, where $s_0$ is the reflection with respect to $\alpha_0$. This equips $W_{\af}$ with a length function $\ell \colon W_\af \to \Z_{\ge 0}$ that extends the one on $W$.
For an element $w \in W_\af$, a sequence of indices $i_1, \dots, i_k \in \tI_\af$ is a \emph{reduced expression} for $w$ if $w = s_{i_1} \cdots s_{i_k}$ and $k=\ell(w)$.

The Bruhat order on $W_\af$ is characterized by $w \le v$ if and only if $w$ can be written as a subexpression of a reduced expression of $v$. The embedding $Q \hookrightarrow W_\af$ associates each $\gamma \in Q$ with a translation element $t_\gamma$. A key relation involving these elements is $t_{- \vartheta} = s_{\vartheta} s_0$.

\subsection{Twisted affine Lie algebras}\label{aff-Lie}

Let $\tg$ be the twisted affinization of $\g$. When $\g$ is of type $\mathsf{A}$, $\mathsf{D}$, or $\mathsf{E}$, this is just the untwisted affinization. Then $\tg$ is the affine Kac--Moody algebra associated with the root system $\Delta_\af$. The resulting affine Kac--Moody algebras are of the following types:
\[
\mathsf A_{\ell}^{(1)}, \mathsf D_{\ell}^{(1)}, \mathsf E_{6}^{(1)}, \mathsf E_7^{(1)}, \mathsf E_8^{(1)} \quad \text{and} \quad \mathsf A_{2\ell-1}^{(2)}, \mathsf D_{\ell+1}^{(2)}, \mathsf E_6^{(2)}, \mathsf D_4^{(3)}.
\]
The set $\tI$ of vertices of the Dynkin diagram of the finite-dimensional Lie algebra $\g$ is obtained from the set $\tI_\af$ of vertices of the Dynkin diagram of $\tg$ by removing the vertex $0$; see \cite[Chap.~4, Table~Aff]{Kac}.

\subsubsection{Triangular decomposition}
We have the triangular decomposition
\[
\tg = \tn^- \oplus \wth \oplus \tn,
\]
and the finite-dimensional decomposition $\g = \gn^- \oplus \h \oplus \gn$ embeds into it component-wise. The components are given by
\begin{equation*}
    \tn = \bigoplus_{\al \in \Delta_\af^+} \tg_{\al}, \quad
    \wth = \h \oplus \C K \oplus \C d, \quad \text{and} \quad
    \tn^- = \bigoplus_{\al \in \Delta_\af^+} \tg_{-\al}.
\end{equation*}
Here, $K$ is the central element ($[K, \tg] = 0$), $d$ is the derivation, and the enlarged Cartan subalgebra $\wth$ is abelian. For $\al \in \Delta_\af$, the root space $\tg_\al$ is the corresponding eigenspace for the adjoint action of $\wth$. Here we view $\al$ as an element of $\wth^*$ via the inclusion $\Delta_\af \subset (\h \oplus \C d)^* \subset \wth^*$.

We extend the definition of simple coroots from the finite case by setting $\al_0^{\vee} = K - \vartheta^{\vee}$. We write $\wg := [\tg,\tg]$ for the derived subalgebra, and set $\tb := \wth \oplus \tn$ and $\tb^- := \wth \oplus \tn^-$.

Finally, let $\tg_{\ge 0}$ and $\tg_{\le 0}$ denote the direct sums of the non-negative and non-positive $d$-eigenspaces of $\tg$, respectively. These are Lie subalgebras of $\tg$.

\subsubsection{Simple roots and fundamental weights}\label{subsubsec:sf}
For each $i \in \tI_\af$, let $E_i \in \tg_{\al_i}$ and $F_i \in \tg_{-\al_i}$ be the standard Kac--Moody generators, satisfying $[E_i,F_i] = \al_i^{\vee} \in \wth$ (see \cite[\S 1.2]{Kac}). The Chevalley involution $\theta$ on $\tg$ is an automorphism that satisfies
\[ \theta ( E_i ) = F_i, \quad \theta ( F_i ) = E_i \quad (i \in \tI_\af) \quad \text{and}\quad \theta ( d ) = -d. \]
Let $\{\La_i\}_{i \in \tI_\af} \subset \wth^*$ be the set of fundamental weights of $\tg$, and let $\rho \in \wth^*$ be the Weyl vector. These are characterized by
\[ \langle \La_i, \al_j^{\vee} \rangle = \delta_{ij}, \quad \langle \rho, \al_i^{\vee} \rangle = 1, \quad \text{and} \quad \langle \La_i, d \rangle = \langle \rho, d \rangle = 0. \]
We set $P_\af := \bigoplus_{i \in \tI_\af} \Z \La_i \oplus \Z \delta$ and $P_\af^+ := \bigoplus_{i \in \tI_\af} \Z_{\ge 0} \La_i \oplus \Z \delta$. The inclusion $\h \hookrightarrow \wth$ induces a classical projection map $\lambda \mapsto \overline{\lambda}$ from $P_\af$ onto $P$, which satisfies
\[ \overline{\La}_i = \varpi_i \quad (i \in \tI) \quad \text{and} \quad \overline{\La}_0 = \overline{\delta} = 0. \]
We also fix a section $P \hookrightarrow P_\af$ of the classical projection by extending the assignment $\varpi_i \mapsto \La_i - \langle \varpi_i, \vartheta^{\vee} \rangle \La_0$ for $i \in \tI$ linearly.

The affine Weyl group $W_\af$ acts on $P_\af$ such that each simple reflection $s_i$ ($i \in \tI_\af$) acts as the reflection with respect to the coroot $\al_i^\vee$.

\subsubsection{Affine weights of level $k$}
For each $k \in \Z$, we define the set of \emph{level $k$ weights} as
\[ P_{\af,k} := \{\La \in P_\af \mid \langle \La, K \rangle = k \}. \]
The set of \emph{level $k$ dominant weights}, denoted $P_{\af, k}^+$, consists of weights $\La \in P_{\af,k}$ that also satisfy $\langle \La, \al_i^{\vee} \rangle \ge 0$ for all $i \in \tI_\af$.

For a positive level $k \in \Z_{>0}$, the set $P_{\af, k}^+$ has only finitely many orbits under the additive action of $\Z \delta$. For each dominant weight $\La \in P_{\af, k}^+$, there exists a unique (up to isomorphism) integrable highest-weight representation $L(\La)$ of level $k$. This module contains a highest-weight vector $\bv_\La$, and for any $w \in W_\af$, the $w\La$-weight space of $L(\La)$ is one-dimensional. We denote a vector spanning this space by $\bv_{w\La}$, which is unique up to scalar.

By~\cite[Corollary~10.1]{Kac}, each $W_\af$-orbit in $P_{\af,k}$ meets $P_{\af,k}^+$ in exactly one point.
In particular, for any $\la \in P$, the weight $\la + k\La_0$ can be written in the form
\begin{equation}
\la + k\La_0 = w\La\label{eqn:wt-conv}
\end{equation}
for a unique $\La \in P_{\af,k}^+$ and some $w \in W_\af$ which is unique up to $\mathsf{Stab}_{W_\af} ( \La )$.

\subsubsection{Definition of the modules $D^{(k)}_\la$ and $\widehat{\bD}^{(k)}_{\la}$}
Using the correspondence from~\eqref{eqn:wt-conv}, we define the \emph{level $k$ thin Demazure module} as
\[ D^{(k)}_\la := U(\tb) \bv_{w\La} \subset L(\La). \]
We remark that $D^{(k)}_\la$ is finite-dimensional; indeed, its character is obtained from the highest-weight character by finitely many applications of Demazure operators \cite[Proposition~8.1.17]{Kum02}.
Similarly, we define the \emph{level $k$ thick Demazure module} as
\[ L(\La)^w := U(\tb^-)\bv_{w\La} \subset L(\La). \]
We then twist the $\tb^-$-module structure on $L(\La)^w$ by the Chevalley involution $\theta$ to obtain a $\tb$-module, which we denote by $\widehat{\bD}^{(k)}_{-\la}$. When convenient, we use the notation ${}^\theta M$ for a $\tg$-module $M$ with the twisted action; thus, $\widehat{\bD}^{(k)}_{-\la}$ is the module ${}^\theta(L(\La)^w)$, which we often denote by ${}^\theta L(\La)^w$.

Both modules $D^{(k)}_\la$ and $\widehat{\bD}^{(k)}_{\la}$, when regarded as $\tb$-modules, are $\wth$-semisimple. Each is generated by a cyclic vector whose $\h$-weight is $\la$.

As a $\tg$-module, ${}^\theta L(\La)$ has level $-k$. Consequently, $\widehat{\bD}^{(k)}_{-\la}$ can be viewed as a submodule of the irreducible lowest-weight integrable representation of $\tg$ with lowest weight $-\La$.

\subsubsection{Order relations on $P$ and $W_\af$}

\begin{defn}\label{def:ordering}
For each $\la \in P$, let $\la_+$ and $\la_-$ denote the unique elements in $W \la \cap P^+$ and $W \la \cap (-P^+)$, respectively. Let $\preceq$ be the partial order on $P$ defined by the following conditions:
\begin{align*}
\la \preceq \mu \quad \Longleftrightarrow \quad & (\mu_+ - \la_+ \in Q_+ \text{ and } \la_+ \neq \mu_+), \quad \text{or} \\
& (\la_+ = \mu_+ \text{ and } \la - \mu \in Q_+).
\end{align*}
We write $\la \prec \mu$ if $\la \preceq \mu$ and $\la \neq \mu$. We set
\[ \Sigma (\la) := \{ \mu \in P \mid \mu \preceq \la\} \quad \text{and} \quad \Sigma_* (\la) := \{ \mu \in P \mid \mu \prec \la\}. \]
The partial order $\preceq$ equips $P_\af$ with a preorder induced by the projection map $\la \mapsto \overline{\la}$. We define the shifted action of $W_\af$ on $\R \otimes_\Z P$ as
\[ w(\!(\la)\!) := \overline{w(\la + \La_0) - \La_0} \quad \text{for} \quad \la \in \R \otimes_\Z P, w \in W_\af. \]
\end{defn}

Note that we have the relation $\la_+ = w_0 \la_-$.

\begin{rem}
The order $\preceq$ originates in the work of Heckman~\cite{Hec87}, and is also recorded in Macdonald's Bourbaki lecture~\cite[\S4.10]{Mac95B}.
\end{rem}

\begin{lem}[\cite{CK18}~\S2.2.1]\label{lem:shiftconv}
For each $\la \in P$, we have
\[ \Sigma ( \la ) = \left( \mathrm{Conv} \,\Sigma ( \la )\right)\cap ( \la + Q ), \]
where $\mathrm{Conv}$ denotes the convex hull in $\R \otimes_\Z P$.
\end{lem}

\begin{prop}[\cite{CK18} \S4.1.1]\label{prop:ordercomp}
Let $k \in \Z_{>0}$. For each $\La \in P_{\af, k}^+$ and $w,v \in W_\af$, if $w \le v$ in the Bruhat order, then $\overline{w \La} \preceq \overline{v \La}$.
\end{prop}

\subsubsection{Definition of the module $\bD^{(k)}_{\la}$}

\begin{thm}[\cite{Kas93}; see also {\cite[Theorem~C]{Kat18b}}]\label{thm:D-incl}
Let $k \in \Z_{>0}$. For each $\La \in P_{\af, k}^+$ and $w,v \in W_\af$ such that $w \le v$, we have an inclusion of $\tb^-$-modules $L(\La)^w \supset L(\La)^v$. This inclusion is an equality if and only if $w \La = v \La$.
\end{thm}

\begin{cor}[{\cite[Theorem~C]{Kat18b}}; cf.\ \cite{Kas93}]\label{cor:D-dist}
Let $k \in \Z_{>0}$. For each $\La \in P_{\af, k}^+$ and $w,v \in W_\af$, there exists a subset $S(w,v) \subset W_\af$ such that
\[ L(\La)^w \cap L(\La)^v = \sum_{u \in S(w,v)} L(\La)^u \subset L(\La). \]
\end{cor}

Let $k \in \Z_{>0}$. For each $\La \in P_{\af, k}^+$ and $w \in W_\af$, we define
\[ \gr^w L(\La) := \frac{L(\La)^w}{\sum_{v > w, \, w\La \neq v\La} L(\La)^v}. \]
Assuming the correspondence $\la + k\La_0 = w\La$ from~\eqref{eqn:wt-conv}, we denote by $\bD^{(k)}_{-\la}$ the $\tb$-module obtained by twisting the action on $\gr^w L(\La)$ by the Chevalley involution.
Since $\gr^w L(\La)$ is a quotient of $L(\La)^w$ as a $\tb^-$-module, it follows that $\bD^{(k)}_{-\la}$ is a quotient of $\widehat{\bD}^{(k)}_{-\la}$ as a $\tb$-module. In particular, $\bD^{(k)}_{-\la}$ is a cyclic $\tb$-module generated by a vector of $\h$-weight $-\la$.

Thus, for each $\la \in P$, we obtain three related $\tb$-modules: the thin Demazure module $D^{(k)}_\la$, the thick Demazure module $\widehat{\bD}^{(k)}_{\la}$, and its graded quotient $\bD^{(k)}_{\la}$.

\begin{cor}\label{cor:LtoDfilt}
Let $k \in \Z_{>0}$. For each $\La \in P_{\af, k}^+$ and $w \in W_\af$, the module ${}^\theta L(\La)^w$ admits a decreasing, separable filtration whose associated graded components are direct sums of modules from the set $\{ \bD^{(k)}_\mu \mid \mu \in P \}$.
\end{cor}

\begin{proof}
Applying the involution $\theta$, the assertion is equivalent to the claim that $L(\La)^w$ has a decreasing, separable filtration whose successive quotients are isomorphic to modules of the form $\gr^v L(\La)$. This is a direct consequence of Theorem~\ref{thm:D-incl} and Corollary~\ref{cor:D-dist}.
\end{proof}

\subsection{Categories of representations}\label{rep-current}

We set $Q_\af^+ := \sum_{\al \in \Delta_\af^+} \Z_{\ge 0} \al \subset P_\af$. For an inclusion of Lie algebras $\mathfrak{q} \subset \mathfrak{p}$ with $\mathfrak{q}$ abelian, we write $(\mathfrak{p},\mathfrak{q})$ for the category of $\mathfrak{p}$-modules on which $\mathfrak{q}$ acts semisimply.

A $\tg_{\ge 0}$-module $M$ is said to be \emph{$\g$-integrable} if its restriction to $\g$ decomposes as a direct sum of finite-dimensional $\g$-modules.

\begin{defn}[Category $\fB$]
A $\tb$-module $M$ is called \emph{graded} if it is $\wth$-semisimple, all $\wth$-eigenvalues belong to $P_\af$, and every $\wth$-weight space has at most countable dimension. Let
\[
\Psi(M) := \{ \La \in P_\af \mid \Hom_{\wth} ( \C_{\La}, M ) \neq 0 \}
\]
denote the set of $\wth$-weights of $M$.
The module $M$ is said to be \emph{bounded} if its $d$-degrees are bounded from below and each $d$-weight space is finite-dimensional.
Let $\fB$ and $\fB_{\bd}$ be the full subcategories of the category of $\tb$-modules consisting of graded modules and bounded graded modules, respectively.
\end{defn}

\begin{defn}[Category $\fC$]
A $\tg_{\ge 0}$-module $M$ is called \emph{graded} if its restriction to $\tb$ is a graded module. Let $\fC$ be the full subcategory of graded $\tg_{\ge 0}$-modules whose objects are also $\g$-integrable. Let $\fC_{\bd}$ denote the full subcategory of $\fC$ consisting of modules whose restrictions to $\tb$ belong to $\fB_\bd$.
\end{defn}

Note that any module $M$ in $\fB_\bd$ or $\fC_\bd$ is generated by (any $\wth$-splitting of) its head, $\hd M$. Indeed, boundedness with respect to the $d$-grading implies that $M$ is generated by any $\wth$-splitting of its head.

\subsubsection{Simple and projective objects in $\fB$}\label{subsubsec:Bbasic}

A simple object in $\fB$ is a one-dimensional representation corresponding to a character of $\wth$.

For each $\La \in P_\af$, we define the standard module
\[ Q_\La := U(\tb) \otimes_{U(\wth)} \C_{\La}. \]
The module $Q_\La$ is the projective cover of the simple module $\C_{\La}$ in the category $\fB$.
The projective cover in $\fB$ of an object of $\fB_{\bd}$ is a direct sum of modules of the form $Q_\Lambda$, with multiplicities determined by the $\wth$-weights occurring in its head.

\subsubsection{Simple and projective objects in $\fC$}\label{subsubsec:Cbasic}

For each $\la \in P^+$, the simple $\g$-module $V_\la$ can be regarded as a graded $\tg_{\ge 0}$-module via the natural surjection $\tg_{\ge 0} \to \g$ that annihilates summands with positive $d$-degree. A simple object in $\fC$ is obtained from $V_\la$ for some $\la \in P^+$ by tensoring with a one-dimensional character. The parabolic Verma module
\[ P_\la := U(\tg_{\ge 0}) \otimes_{U(\g + \wth)} V_\la \]
is the projective cover of $V_\la$ in the category $\fC$ (\cite{CG07}).

For any affine weight $\La \in P_\af$, we define the corresponding projective module $P_\La$ and simple module $V_\La$ by twisting their finite counterparts. Let $\la := \overline{\La}$ be the classical projection and let $\la_+$ be the dominant weight in its $W$-orbit (the same one as in Definition~\ref{def:ordering}). We then define
\[ P_\La := P_{\la_+} \otimes \C_{\La - \la} \quad \text{and} \quad V_\La := V_{\la_+} \otimes \C_{\La - \la}. \]
Here, we identify $\la$ with its canonical image in $P_\af$ to form the difference $\La - \la$. The module $V_\La$ is regarded as a $(\g + \wth)$-module.
Every projective object in $\fC_\bd$ is a direct sum of modules of the form $P_\La$ $(\La \in P_\af)$. This follows from the fact that each object of $\fC_\bd$ is generated by any $\wth$-splitting of its head.

\subsubsection{Graded Homs and Exts}\label{subsubsec:ghom}

We now record the grading conventions and the derived functors used throughout the paper.

The four categories defined above are abelian. They are equipped with shift functors $\bq^m$ for each $m \in \Z$, which send a module to its tensor product with the one-dimensional representation $\C_{m\delta}$. Namely, if we are given a module $M$ equipped with a graded decomposition $M = \bigoplus_{s \in \Z} M_s$, then for each $m, s \in \Z$, we define
\[
( \bq^m M )_s := M_{s - m}.
\]
This functor commutes with natural functors among our categories.
Using the shift functors $\bq^m$, we define the graded Hom-space between objects $M$ and $N$ in one of these categories $\mathcal{C}$ as
\begin{equation}
\shom_{\mathcal{C}}(M, N) := \bigoplus_{m, l \in \Z} \mathrm{Hom}_{\mathcal{C}}(\bq^m M, N \otimes \C_{l\La_0}).\label{eqn:defhom}
\end{equation}
We similarly define the graded Ext-groups using a projective resolution of $M$:
\begin{equation}
\ext^i_{\mathcal{C}}(M, N) := \bigoplus_{m, l \in \Z} \mathrm{Ext}^i_{\mathcal{C}}(M, \bq^{-m} N \otimes \C_{l\La_0}).\label{eqn:defext}
\end{equation}
We consider these as graded vector spaces, where the grading corresponds to the $d$-degree (i.e., the $\bq$-shifts). The extra $\Lambda_0$-twists are included in order to keep track of the $K$-eigenvalue separately from the $d$-grading.

Despite the conventional isomorphism of graded vector spaces
\begin{equation}
\ext^\bullet_{\mathcal C} ( M, N ) \cong \ext^\bullet_{\mathcal C} ( M, N \otimes \C_{\La_0} ),\label{eqn:Lneglect}
\end{equation}
we generally keep track of the grading by $\La_0$-twists.

\subsubsection{Restricted duals}\label{subsubsec:rd}

For a module $M$ in $\fB_\bd$ or $\fC_\bd$, its {\it restricted dual} is defined as
\[ M^{\vee} := \bigoplus_{\La \in P_\af} \mathrm{Hom}_{\wth}(\C_{\La}, M)^*. \]
The dual module $M^{\vee}$ is an object in $\fB$ or $\fC$, respectively. In case $\dim M < \infty$, the full dual $M^*$ coincides with $M^{\vee}$.

Moreover, for any $M,N\in \fB_{\bd}$, there are natural isomorphisms
\begin{equation}\label{eqn:transpose}
    \shom_{\fB}(M, N^{\vee}) \cong \shom_{\fB}(\C, M^{\vee} \otimes N^{\vee} ) \cong \shom_{\fB}(N, M^{\vee}).
\end{equation}
These isomorphisms are functorial. Standard arguments in homological algebra show that they extend to the higher Ext-groups as well; see \cite[Chap.~2~\S1]{Gro57}.

\begin{rem}
The extension of~\eqref{eqn:transpose} to higher Ext-groups is obtained by
appealing to the formalism of $\delta$-functors, in particular through the
universal property of derived functors regarded as $\delta$-functors%
~\cite[Proposition~2.2.1]{Gro57}.

If~\eqref{eqn:transpose} were available functorially for all $\wth$-semisimple
$\tb$-modules, then the passage to higher
Ext-groups would follow directly from the formalism of derived categories
(cf.~\cite{Ver96}). However, since $\shom$ does not in general commute with
direct sums, such a functorial extension is not valid in complete generality.
To address this issue, we restrict to situations where~\eqref{eqn:transpose}
holds under certain boundedness assumptions, and we promote it to a natural
transformation by fixing one of the arguments.

This technical constraint explains the asymmetric role of $\fB_\bd$ in our
setup. Indeed, objects in $\fB_\bd$ are compact in the derived category of $\fB$ (cf.~\eqref{eqn:cpt}), and this compactness
ensures the desired control over direct sums.

In the sequel, we shall refrain from further use of the $\delta$-functor
formalism. All arguments will instead be carried out within the more concrete
setting of standard homological algebra (e.g.,~\cite{CE56}).
\end{rem}

\subsubsection{Filtrations}

Consider a decreasing, separable filtration of an object $M$ in one of the categories defined above:
\[
M = M^0 \supset M^1 \supset M^2 \supset \cdots, \qquad \text{with} \quad \bigcap_{r \ge 0} M^r = \{0\}.
\]
Assume that this filtration is exhaustive with respect to the $d$-grading; that is, for any degree bound $m \in \Z$, there exists an integer $r$ such that all weights in the support $\Psi(M^r)$ have $d$-degree at least $m$.

Under this assumption, the object $M$ is complete with respect to the filtration, i.e.,
\[
\varprojlim_{r} M/M^r \cong M
\]
as an object in the respective category.

In particular, this implies that $M$ is complete with respect to the topology induced by its $d$-grading. It is important to note, however, that this notion of completion—taken within specific categories such as $\fB$ or $\fC$—may differ from the standard completion in the larger category of all modules. This distinction will be used implicitly at various places in the sequel.

\subsubsection{Yoneda interpretation of $\ext^1$}\label{subsubsec:Yoneda}

We include the following only to justify our later use of $\ext^1$ as a space of extension classes.

\begin{prop}[See e.g., {\cite[Theorem~3.4.3]{Wei94}}]\label{prop:Yoneda}
Let $\mathcal{A}$ be an abelian category, and let $M, N \in \mathcal{A}$. Then there is a natural bijection
\[
\Ext^1_{\mathcal{A}}(M, N) \;\cong\; \left\{
\begin{array}{c}
\text{equivalence classes of short exact sequences of the form} \\
0 \to N \to E \to M \to 0
\end{array}
\right\}.
\]
Two extensions are considered equivalent if there exists a commutative diagram
\[
\begin{tikzcd}
0 \ar[r] & N \ar[r] \ar[equal]{d} & E \ar[r] \ar[d,"\sim"] & M \ar[r] \ar[equal]{d} & 0 \\
0 \ar[r] & N \ar[r] & E' \ar[r] & M \ar[r] & 0
\end{tikzcd}
\]
in which the middle vertical map is an isomorphism.

In particular, there exists a nontrivial extension of $M$ by $N$ in $\mathcal{A}$ if and only if $\Ext^1_{\mathcal{A}}(M, N) \ne 0$.
\end{prop}

We apply Proposition~\ref{prop:Yoneda} below for $\mathcal{A} = \fB_\bd$ or $\fC_\bd$, replacing $\Ext^1_{\mathcal{A}}(M, N)$ with $\ext^1_{\mathcal{A}}(M, N)$. This amounts to allowing nontrivial extensions of $M$ by $\bq^m N$ for some $m \in \Z$, in accordance with our convention~\eqref{eqn:defext}.

\subsubsection{$\ext$ and inverse limits}

The following result is a variant of~\cite[Theorem~3.5.8]{Wei94}. 

\begin{lem}\label{lem:comm}
Let $M,N \in \fB_{\bd}$, and suppose that
\[
M \;\cong\; \varprojlim_r M(r)
\]
as a projective system $\{M(r)\}_{r\in \Z_{\ge 0}}$ with surjective transition maps. Then, for each $i \in \Z$,
\begin{equation}\label{eqn:Kcptc}
\varinjlim_r \ext^{i}_{\fB}\big(M(r),N^{\vee}\big)
\;\cong\;
\ext^{i}_{\fB}\big(\varprojlim_r M(r),N^{\vee}\big).
\end{equation}
\end{lem}

\begin{proof}
Let $X(r):=M(r)\otimes N$ and
\[
X:=\varprojlim_r X(r) \cong \bigl( \varprojlim_r M(r) \bigr)\otimes N \qquad \text{(since $M,N\in\fB_{\bd}$)}.
\]

By the duality from \S\ref{subsubsec:rd}, we have
\begin{equation}
\ext^{\bullet}_{\fB}\big(M(r),N^{\vee}\big)\;\cong\;
\ext^{\bullet}_{\fB}\big(\C, X(r)^{\vee}\big),\quad
\ext^{\bullet}_{\fB}\big(M,N^{\vee}\big)\;\cong\;
\ext^{\bullet}_{\fB}\big(\C, X^{\vee}\big).\label{eqn:CMN}
\end{equation}

Since all $\wth$-weights of $\tn$ are (strictly) positive, the Koszul resolution of $\C$ has finite-dimensional $\wth$-weight components (and therefore, for each fixed weight, the corresponding component complex is bounded). In particular, we have
\begin{equation}
\ext^{\bullet}_{\fB}\big(\C, \bigoplus_r X(r)^{\vee}\big) \cong \bigoplus_r \ext^{\bullet}_{\fB}\big(\C, X(r)^{\vee}\big).\label{eqn:cpt}
\end{equation}

The transition maps \(\iota_r: X(r)^{\vee}\hookrightarrow X(r{+}1)^{\vee}\) are injective. Hence we have a standard short exact sequence
\[
0 \;\longrightarrow\; \bigoplus_r X(r)^{\vee}
\xrightarrow{\,f\,}
\bigoplus_r X(r)^{\vee}
\;\longrightarrow\; \varinjlim_r X(r)^{\vee}
\;\longrightarrow\; 0,
\]
where $f(\{a_r\}_r) = \{\iota_{r-1}(a_{r-1}) - a_{r}\}_r$ with $\iota_{-1} ( a_{-1} ) := 0$, and $\varinjlim_r X(r)^{\vee} \cong \bigl( \varprojlim_r X(r) \bigr)^{\vee}$.  

Applying $\ext^{\bullet}_{\fB}(\C,-)$ yields the long exact sequence
\begin{align*}
0 &\
   \to \hom_{\fB}(\C, \bigoplus_r X(r)^{\vee})
   \to \hom_{\fB}(\C, \bigoplus_r X(r)^{\vee}) \to \hom_{\fB}(\C,\varinjlim_r X(r)^{\vee})\\
&\to \ext^1_{\fB}(\C,\bigoplus_r X(r)^{\vee}) \to \cdots .
\end{align*}

On the other hand, taking inductive limits after applying $\ext^{\bullet}_{\fB}(\C,-)$ to the system $\{X(r)^{\vee}\}_r$ gives another long exact sequence
\begin{align*}
0 &
   \to \bigoplus_r \hom_{\fB}(\C,X(r)^{\vee})
   \to \bigoplus_r \hom_{\fB}(\C,X(r)^{\vee}) \to \varinjlim_r \hom_{\fB}(\C,X(r)^{\vee})\\
&\to \bigoplus_r \ext^1_{\fB}(\C,X(r)^{\vee}) \to \cdots .
\end{align*}

By~\eqref{eqn:cpt}, these two long exact sequences fit into a commutative diagram
\begin{equation}\label{eqn:comm-amp}
\xymatrix{
\bigoplus_r \ext^{i}_{\fB}(\C,X(r)^{\vee})\ar[r]\ar[d] &\bigoplus_r \ext^{i}_{\fB}(\C,X(r)^{\vee})\ar[d]
\\
\ext^{i}_{\fB}(\C, \bigoplus_r X(r)^{\vee}) \ar[r] &
\ext^{i}_{\fB}(\C, \bigoplus_r X(r)^{\vee})
}
\end{equation}
for each $i \in \Z$. In addition, we have natural maps
\[
\varinjlim_r \ext^{i}_{\fB}(\C,X(r)^{\vee}) \longrightarrow \ext^{i}_{\fB}(\C,X^{\vee}) \qquad (i \in \Z)
\]
which commute with~\eqref{eqn:comm-amp}. By the five lemma, we obtain
\[
\varinjlim_r \ext^{\bullet}_{\fB}(\C,X(r)^{\vee})
\;\cong\;
\ext^{\bullet}_{\fB} (\C,\varinjlim_r X(r)^{\vee} ).
\]

Combining with~\eqref{eqn:CMN} yields~\eqref{eqn:Kcptc}.
\end{proof}

\begin{cor}\label{cor:comm}
Let $M, N \in \fB_{\bd}$, and suppose that $N$ is presented as an inverse limit
\[
N \;\cong\; \varprojlim_r N(r),
\]
where $\{N(r)\}_r$ has surjective transition maps. Then there are natural isomorphisms
\[
\varinjlim_r \ext^{\bullet}_{\fB}\!\big(M, N(r)^{\vee}\big)
\;\cong\;
\ext^{\bullet}_{\fB}\!\Big(M, (\varprojlim_r N(r))^{\vee}\Big)
\;\cong\;
\ext^{\bullet}_{\fB}\!\Big(M, \varinjlim_r N(r)^{\vee}\Big).
\]
\end{cor}

\begin{proof}
By the duality from~\S\ref{subsubsec:rd}, we have
\[
\ext^{\bullet}_{\fB}\!\big(M, N(r)^{\vee}\big) \cong \ext^{\bullet}_{\fB}\!\big( N(r), M^{\vee}\big), \qquad
\ext^{\bullet}_{\fB}\!\big(M, N^{\vee}\big) \cong \ext^{\bullet}_{\fB}\!\big( N, M^{\vee}\big). 
\]
Applying Lemma~\ref{lem:comm}, we obtain
\[
\varinjlim_r \ext^{\bullet}_{\fB}\!\big(M, N(r)^{\vee}\big)
\;\cong\;
\ext^{\bullet}_{\fB}\!\Big(M, (\varprojlim_r N(r))^{\vee}\Big).
\]
Finally, in $\fB_{\bd}$ the dual interchanges inductive limits and inverse limits
(since weight spaces are finite-dimensional); hence $\varinjlim_r N(r)^{\vee} \cong (\varprojlim_r N(r))^{\vee}$, and this yields the last isomorphism.
\end{proof}

\subsection{BGG resolution and Demazure functors}\label{subsec:BGGD}

\subsubsection{BGG resolutions}
For each $\La \in P_\af$, we define the \emph{Verma module} $M(\La)$ and its Chevalley-twisted analogue ${}^{\theta}M(\La)$ by
\[ M(\La) := U(\tg) \otimes_{U(\tb)} \C_{\La}, \quad \text{and} \quad {}^{\theta}M(\La) := U(\tg) \otimes_{U(\tb^-)} \C_{-\La}. \]
Since ${}^\theta M(\La)$ is isomorphic to the standard projective module $Q_{-\La}$ from the definition, it is a projective object in the category $\fB$.

For any dominant integral weight $\La \in P^+_\af$, the simple module $L(\La)$ admits a \emph{BGG resolution}; namely, there is an exact sequence~\cite[Theorem~9.1.3]{Kum02}
\begin{equation}\label{eqn:BGG-res}
	\cdots \to \bigoplus_{\ell(w) = 2} M(w \cdot \La) \to \bigoplus_{\ell(w) = 1} M(w \cdot \La) \to M(\La) \to L(\La) \to 0,
\end{equation}
where $w \cdot \La := w(\La + \rho) - \rho$ denotes the shifted ``dot'' action of the affine Weyl group.

Thus, applying the Chevalley involution $\theta$ to the BGG resolution~\eqref{eqn:BGG-res} yields a projective resolution of the module ${}^{\theta}L(\La)$ in $\fB$.

\subsubsection{Demazure functors}
For each $i \in \tI_\af$, we define the parabolic subalgebra $\tp_i := \tb \oplus \tg_{-\al_i}$. This contains the Lie subalgebra
\[ \mathfrak{sl}(2, i) := \C E_i \oplus \C \al_i^{\vee} \oplus \C F_i. \]
A $\tp_i$-module is called {\it $\mathfrak{sl}(2,i)$-integrable} if its restriction to $\mathfrak{sl}(2,i)$ decomposes into a direct sum of finite-dimensional $\mathfrak{sl}(2,i)$-modules.

The {\it Demazure functor} $\mathscr{D}_i$ is an endofunctor on $\fB$ that sends a module $M$ to the maximal $\mathfrak{sl}(2,i)$-integrable quotient of the induced module $U(\tp_i) \otimes_{U(\tb)} M$.

\begin{thm}[Joseph~\cite{Jos85}]\label{thm:DJ}
Let $i,j \in \tI_\af$. The Demazure functors have the following properties:
\begin{enumerate}
    \item There is a natural transformation of functors $\mathrm{Id} \to \mathscr{D}_i$.
    \item There is a natural isomorphism of functors $\mathscr{D}_i \cong \mathscr{D}_i \circ \mathscr{D}_i$.
    \item The functor $\mathscr{D}_i$ is right exact, and its left derived functors satisfy $\mathbb{L}^m \mathscr{D}_i = 0$ for $m \neq 0, -1$.
    \item If $M \in \fB$ is the restriction of an $\mathfrak{sl}(2,i)$-integrable $\tp_i$-module, then there is a natural isomorphism of functors
    \[ M \otimes \mathscr{D}_i(\bullet) \cong \mathscr{D}_i(M \otimes \bullet). \]
    \item (Braid Relations) If $s_i$ and $s_j$ satisfy $\overbrace{s_i s_j s_i \cdots}^m = \overbrace{s_j s_i s_j \cdots}^m$, then
    \[ \overbrace{\mathscr{D}_i \mathscr{D}_j \mathscr{D}_i \cdots}^m \cong \overbrace{\mathscr{D}_j \mathscr{D}_i \mathscr{D}_j \cdots}^m.\]
\end{enumerate}
\end{thm}

We choose a reduced expression $s_{i_1} \cdots s_{i_\ell}$ of $w \in W_\af$.
As a consequence of the braid relations in Theorem~\ref{thm:DJ}(5), the composite functor
\begin{equation}\label{Dwdef}
    \mathscr{D}_w := \mathscr{D}_{i_1} \circ \mathscr{D}_{i_2} \circ \cdots \circ \mathscr{D}_{i_\ell}
\end{equation}
is well-defined and does not depend on the choice of reduced expression for $w$.

The following observation shows that the Demazure functor depends only on the relevant $\mathfrak{sl}(2)$-structure, and not on the ambient Lie algebra used to define it.

\begin{cor}\label{cor:eDindep}
Let $i \in \tI_\af$, and consider a sequence of inclusions of $\wth$-stable Lie algebras
\begin{equation}
\mathfrak{sl}(2, i) \subset \mathfrak{q} = \mathfrak{sl}(2, i) \oplus \mathfrak{q}_{>0} \subset \tp_i,\label{eqn:Qinclstr}
\end{equation}
with $\wth \subset \mathfrak{q}$ and $\mathfrak{q}_{>0} \subset \tb$.

Then, for any $\wth$-semisimple $(\mathfrak{q} \cap \tb)$-module $M$, the Demazure module $\mathscr{D}_i(M)$ computed using $(\mathfrak{sl}(2,i) \cap \tb)$-structure coincides with that computed using $(\mathfrak{q} \cap \tb)$-structure.
\end{cor}

\begin{proof}
Let $\gb_i := \mathfrak{sl}(2, i) \cap \tb = \C \alpha_i^{\vee} \oplus \C E_i$, and set $\mathfrak{q}_{\ge 0} := \mathfrak{q} \cap \tb$.

The decomposition $\mathfrak{q} = \mathfrak{sl}(2, i) \oplus \mathfrak{q}_{>0}$ is compatible with the $\wth$-action, and $\mathfrak{q}_{>0}$ decomposes into finite-dimensional $\mathfrak{sl}(2, i)$-modules. We may write
\begin{equation}
\mathfrak{q}_{>0} = \wth' \oplus \bigoplus_{\beta \in \Delta'} \g_\beta,\label{eqn:gqstr}
\end{equation}
where $\wth' \subset \wth$ and $\Delta' \subset \Delta^+_\af \setminus \{\alpha_i\}$ is $s_i$-stable, again using the assumptions in~\eqref{eqn:Qinclstr}.

Let $\widetilde{M} := U(\mathfrak{q}) \otimes_{U(\mathfrak{q}_{\ge 0})} M \cong U(\mathfrak{sl}(2, i)) \otimes_{U(\gb_i)} M$.
Let $\mathscr{D}_i(M)$ and $\mathscr{D}_i'(M)$ denote the Demazure modules obtained by imposing maximality with respect to $\mathfrak{q}_{\ge 0}$ and $\gb_i$, respectively. That is,
\[
\mathscr{D}_i(M) = \widetilde{M}/R, \qquad \mathscr{D}_i'(M) = \widetilde{M}/R',
\]
where $R$ (resp. $R'$) is the largest submodule of $\widetilde M$ such that every simple $\mathfrak q$-module quotient (resp. simple $\mathfrak{sl}(2,i)$-module quotient) contains an infinite-dimensional irreducible $\mathfrak{sl}(2,i)$-constituent.

Here $R'\subset R$, because being stable under the larger Lie algebra imposes a stronger condition. Suppose $R' \subsetneq R$. Then the quotient $R / R'$ has a simple $\mathfrak{sl}(2,i)$-module quotient which, by construction, must be finite-dimensional. Since the $\wth$-action on $R/R'$ is semisimple, it follows that $R/R'$ admits a nonzero simple finite-dimensional $(\mathfrak{sl}(2,i) + \wth)$-module quotient.

However, any such simple $(\mathfrak{sl}(2,i) + \wth)$-module extends to a $\mathfrak{q}$-module via~\eqref{eqn:gqstr}. Indeed, each $\g_\beta$ (for $\beta \in \Delta'$) strictly raises $\wth$-weights, so every nonzero action by $\g_\beta$ maps into a proper submodule. This contradicts the defining property of $R$, whose simple $\mathfrak{q}$-module quotients contain an infinite-dimensional irreducible $(\mathfrak{sl}(2,i) + \wth)$-constituent.

We conclude that $R = R'$, and hence $\mathscr{D}_i(M) = \mathscr{D}_i'(M)$ as claimed.
\end{proof}

By Corollary~\ref{cor:eDindep}, the Demazure functors $\mathscr{D}_i$ and $\mathscr{D}_w$ are independent of the ambient Lie algebra module structures whenever they are well-defined.

A module $M \in \fB$ extends to an integrable $\g$-module if and only if $M \cong \mathscr{D}_{w_0}(M)$, equivalently, if and only if $\mathscr{D}_i(M) \cong M$ for all $i \in \tI$.

\subsubsection{The adjoint property}

The following result is stated under the assumption $\dim N < \infty$ in~\cite{FKM}, but holds whenever~\eqref{eqn:transpose} holds:

\begin{thm}[{\cite[Proposition~5.7]{FKM}}; cf.\ {\cite[Lemma~8]{Mat89c}}]\label{thm:adj}
For any $M, N \in \fB_\bd$, there is a functorial isomorphism between hypercohomology groups:
\[ \ext^m_{\fB}(\mathbb{L}^{\bullet} \mathscr{D}_i(M), N^{\vee}) \cong \ext^m_{\fB}(M, \mathbb{R}^{\bullet} ( \mathscr{D}_i (N) )^{\vee}) \quad (m \in \Z). \]
\end{thm}

\section{Basic module-theoretic results}\label{sec:module}

Throughout this section, we keep the notation and assumptions of the previous section.

\subsection{Comparison of $\Ext$s}

\begin{prop}[Cartan–Eilenberg~\cite{CE56}]\label{prop:DCext}
For each $M, N \in \fC$ such that $\dim N < \infty$, we have
\[
\Ext^{\bullet}_{\fC}(M,N)=\Ext^{\bullet}_{\fB}(M,N).
\]
In particular, we have $\ext^{\bullet}_{\fC} ( M, N ) = \ext^{\bullet}_{\fB} ( M, N )$.
\end{prop}

\begin{proof}
Let $\tg_{>0} := \ker ( \tg_{\ge 0} \to \g + \wth )$. The inclusions $\tg_{>0} \subset \tb$ and $\tg_{>0} \subset \tg_{\ge 0}$ form Lie algebra ideals with their quotients $\gb + \wth$ and $\g + \wth$ admit Lie algebra splittings.

For $(\tb,\wth)$-modules $M$ and $N$, the Hochschild--Serre spectral sequence \cite[Chapter XVI, \S6]{CE56} for the ideal $\tg_{>0}\subset \tb$ yields
\[
E_2^{q,p} := \mathrm{Ext}^q_{(\gb+\wth,\wth)} ( \C, \mathrm{Ext}^p_{\tg_{>0}}(M,N) ) \Rightarrow \mathrm{Ext}^{q+p}_{(\tb, \wth)} (M, N).
\]

Since $M, N \in \fC$, the Hochschild–Serre spectral sequence for the ideal $\tg_{>0} \subset \tg_{\ge 0}$ yields
\[
{}' E_2^{q,p} := \mathrm{Ext}^q_{(\g+\wth,\wth)} ( \C, \mathrm{Ext}^p_{\tg_{>0}}(M,N) ) \Rightarrow \mathrm{Ext}^{q+p}_{(\tg_{\ge 0}, \wth)} (M, N).
\]

Hence, the assertion reduces to the analysis of $\mathrm{Ext}^p_{\tg_{>0}}(M,N)$. The module $\mathrm{Ext}^p_{\tg_{>0}}(M,N)$ inherits a $\g$-action, which is $\g$-integrable because $U(\tg_{>0})$ is itself $\g$-integrable under the adjoint action. Indeed, the usual cochain complex computing $\Ext^\bullet_{\tg_{>0}}(M,N)$ is built from $\g$-integrable modules. Owing to the complete reducibility of integrable $\g$-modules in $(\gb,\h)$, we have
\begin{align*}
\mathrm{Ext}^{>0}_{(\gb+\wth,\wth)} ( \C, \mathrm{Ext}^p_{\tg_{>0}}(M,N) ) & = 0, \\
\mathrm{Hom}_{(\gb+\wth,\wth)} ( \C, \mathrm{Ext}^p_{\tg_{>0}}(M,N) ) & = \mathrm{Hom}_{(\g + \wth,\wth)} ( \C, \mathrm{Ext}^p_{\tg_{>0}}(M,N) ).
\end{align*}

The vanishing of $E_2^{>0,\bullet}$ and ${}' E_2^{>0,\bullet}$ implies that the spectral sequences collapse, and the second equality identifies the $E_2$-pages of the Ext-groups. This yields the desired isomorphism
\[
\mathrm{Ext}^{\bullet}_{\fC} ( M, N ) = \mathrm{Ext}^{\bullet}_{\fB} ( M, N ).
\]
From this, we also deduce $\ext^{\bullet}_{\fC} ( M, N ) = \ext^{\bullet}_{\fB} ( M, N )$ as claimed.
\end{proof}

\subsection{Properties of affine Demazure modules}

We next recall the behavior of thin and thick Demazure modules under Demazure functors.

\begin{thm}[see Kumar~{\cite[Chapter~VIII]{Kum02}}, {\cite[Proposition~4.16]{CK18}}]\label{thm:thinD}
Let $k \in \Z_{>0}$. For each $\la \in P$ and $i \in \tI_\af$, we have
\[
\mathbb L^{\bullet} \mathscr D_i ( D_\la^{(k)} ) =
\begin{cases}
\bq ^{\delta _{i0} ( \langle \vartheta^{\vee}, \la \rangle - k )} D_{\overline{s_i ( \la + k \La_0 )}}^{(k)} & \text{if } \overline{s_i ( \la + k \La_0 )} \succ \la, \\
D_{\la}^{(k)} & \text{if } \overline{s_i ( \la + k \La_0 )} \preceq \la.
\end{cases}
\]
\end{thm}

\begin{thm}[Kashiwara~\cite{Kas93}, see also~{\cite[Theorem~C]{Kat18b}}]\label{thm:tD}
Fix $k \in \Z_{>0}$, and let $\La \in P_{\af, k}^+$, $w \in W_\af$, and $i \in \tI_\af$. Then we have
\[
\mathbb L^{\bullet} \mathscr D_i ( {}^{\theta} L ( \La )^w ) =
\begin{cases}
{}^{\theta} L ( \La )^{s_iw} & \text{if } s_i w < w, \\
{}^{\theta} L ( \La )^w & \text{if } s_i w > w.
\end{cases}
\]
\end{thm}

\begin{cor}\label{cor:tD}
Fix $k \in \Z_{>0}$, and let $\La \in P_{\af, k}^+$. For a subset $S \subset W_\af$, define
\[
M := \sum_{w \in S} {}^{\theta} L ( \La )^w, \quad N := \bigcap_{w \in S} {}^{\theta} L ( \La )^w.
\]
Then for each $i \in \tI_\af$, the lower derived functors vanish:
\[
\mathbb{L}^{<0} \mathscr{D}_i ( M ) = \mathbb{L}^{<0} \mathscr{D}_i ( N ) = 0,
\]
and both $M$ and $N$ are preserved under $\mathscr{D}_i$:
\[
M \subset \mathscr{D}_i ( M ), \quad N \subset \mathscr{D}_i ( N ).
\]
\end{cor}

\begin{proof}
By~\cite[Proposition~4.1]{Kas93} (specialized to $q = 1$), the module $L(\La)$ admits a basis such that, for each $w \in W_\af$, a subset of this basis yields a basis of $L(\La)^w$. These bases consist of $\wth$-eigenvectors and satisfy the so-called \emph{string property} of Demazure crystals~\cite[(0.4)(iii), §4]{Kas93}, which is preserved under both sums and intersections of these subspaces.

This property determines decreasing separable filtrations $F_i M$ and $F_i N$ with respect to $E_i$-stable subspaces, such that their associated graded modules decompose as direct sums of indecomposable $\wth$-semisimple $\C[E_i]$-modules (cf.~\cite[Proof of Corollary~4.8]{Kat18}).

By~\cite[(0.5)]{Kas93}, for each $j \ge 0$, we have
\begin{align*}
& \mathbb{L}^{<0} \mathscr{D}_i ( F_j M / F_{j+1} M ) = 0 = \mathbb{L}^{<0} \mathscr{D}_i ( F_j N / F_{j+1} N ), \qquad \text{and}\\
& F_j M / F_{j+1} M \subset \mathscr{D}_i ( F_j M / F_{j+1} M ), \quad F_j N / F_{j+1} N \subset \mathscr{D}_i ( F_j N / F_{j+1} N ).
\end{align*}
Thus, by the d\'evissage applied to the induced filtrations on $M / F_{j+1} M$ and $N / F_{j+1} N$ from $F_\bullet M$ and $F_\bullet N$ for $j \ge 0$, we conclude
\begin{align*}
& \mathbb{L}^{<0} \mathscr{D}_i ( M / F_{j+1} M ) = 0 = \mathbb{L}^{<0} \mathscr{D}_i ( N / F_{j+1} N ), \qquad \text{and}\\
& M / F_{j+1} M \subset \mathscr{D}_i ( M / F_{j+1} M ), \quad N / F_{j+1} N \subset \mathscr{D}_i ( N / F_{j+1} N ).
\end{align*}
Therefore, taking inverse limits, for each $j \ge 0$, we obtain
\begin{align*}
& \varprojlim_j \mathbb{L}^{<0} \mathscr{D}_i ( M / F_{j+1} M ) = 0 = \varprojlim_j \mathbb{L}^{<0} \mathscr{D}_i ( N / F_{j+1} N ), \qquad \text{and}\\
& \varprojlim_j M / F_{j+1} M \subset \mathscr{D}_i ( \varprojlim_j M / F_{j+1} M ), \quad \varprojlim_j N / F_{j+1} N \subset \mathscr{D}_i ( \varprojlim_j N / F_{j+1} N ).
\end{align*}
Here, the functor $\mathscr{D}_i$ commutes with $\varprojlim$ since $M$ and $N$ are direct sums of finitely generated $\wth$-semisimple $\C[E_i]$-modules, and the derived functors $\mathbb{L}^{<0} \mathscr{D}_i$ vanish on each direct summand as well as at each stage of the induced filtrations (cf.~\cite[Proof of Lemma~4.4]{Kat18}). This proves the assertions.
\end{proof}

\begin{cor}\label{cor:dominant}
The setting is as in Theorem~\ref{thm:tD}. The module ${}^{\theta} L ( \La )^{w}$ is $\g$-stable if and only if $\overline{w \La}$ is dominant.
\end{cor}

\begin{proof}
For the “only if” direction, note that ${}^{\theta} L ( \La )^{w}$ being a $\g$-module requires the equality $U(\gb^-)\bv_{w\La} = U(\g)\bv_{w\La}$, which can occur only if $\overline{w\La} \in P^+$. Indeed, $\g$-stability means stability under the negative simple root spaces as well.

For the “if” direction, assume $\overline{w\La} \in P^+$. By combining Proposition~\ref{prop:ordercomp}, Theorem~\ref{thm:D-incl}, and Theorem~\ref{thm:tD}, we see that $\mathscr D_i ( {}^{\theta} L ( \La )^w ) = {}^{\theta} L ( \La )^w$ for each $i \in \tI$. This implies that ${}^{\theta} L ( \La )^w$ is $\g$-stable.
\end{proof}

Let $\la \in P^+$ be such that $- \la + k \La_0 = w \La$ for some $w \in W_\af$ and $\La \in P_{\af, k}^+$. We define
\[
W_{\la}^{(k)} := D_{\la_-}^{(k)}.
\]
We also set
\begin{equation}
\bW_{\la}^{(k)} := \frac{{}^{\theta}L ( \La )^{w_0w}}{\sum_{v \in W_\af,\ \la_- \prec (- \overline{v \La})_-} {}^{\theta}L ( \La )^v}. \label{eqn:defW}
\end{equation}

These modules will serve as the integrable analogues of the thin and thick Demazure modules.

\begin{lem}
Let $k \in \Z_{>0}$. For each $\la \in P^+$, $\La \in P^+_{\af,k}$, and $w \in W_\af$ such that $- \la + k \La_0 = w \La$, the module ${}^{\theta}L ( \La )^{w_0w}$ admits a decreasing separable filtration by grading shifts of the modules $\{ \bW_{\mu}^{(k)} \}_{\mu}$, where $\mu$ runs over the set
\[
\left\{ \mu \in \overline{W_\af \La} \cap P^+ \mid - \mu = \overline{v \La} \ \text{for some } v \ge w \right\}.
\]
\end{lem}

\begin{proof}
By Theorem~\ref{thm:D-incl} and Corollary~\ref{cor:D-dist}, the assertion follows immediately from the definition of $\bW_{\la}^{(k)}$ (see also Corollary~\ref{cor:LtoDfilt}).
\end{proof}

\begin{lem}\label{lem:headW}
Let $k \in \Z_{>0}$. For each $\la \in P^+$, the modules $W_\la^{(k)}$ and $\bW_\la^{(k)}$ are $\tg_{\ge 0}$-modules generated by $V_\la$.
\end{lem}

\begin{proof}
Both $W_\la^{(k)}$ and $\bW_\la^{(k)}$ are quotients of Demazure modules, and hence each admits a unique simple quotient as a $\wth$-semisimple $\tb$-module. Moreover, they are $\g$-stable by construction. The assertion follows from the fact that both $W_\la^{(k)}$ and $\bW_\la^{(k)}$ are cyclic $\tb$-modules generated by a vector of $\h$-weight $\la_-$. In either case, this vector generates an irreducible $\g$-submodule isomorphic to $V_\la$, which is the unique simple quotient of the corresponding $\tg_{\ge 0}$-module.
\end{proof}

\subsection{Filtration of $\bW_\la^{(k)}$ by thick Demazure modules}

The following filtration statement is a straightforward generalization of~\cite{CK18}.

\begin{prop}\label{prop:WDfilt}
Let $k \in \Z_{>0}$. For each $\la \in P^+$, the module $\bW_\la^{(k)}$ admits a filtration as a graded $\tb$-module whose associated graded components are the modules $\{ \bD^{(k)}_{\mu} \}_{\mu \in W \la}$, each occurring with multiplicity one.
\end{prop}

\begin{proof}
The existence of the desired filtration follows from the fact that the class of Demazure submodules is closed under sums and intersections, and these operations are compatible with taking quotients (cf.~Corollary~\ref{cor:D-dist}). The module $\bW_\la^{(k)}$ inherits this lattice structure because it is constructed as a quotient: its parent module, $\widehat{\bD}_{\la_-}^{(k)}$, is identified with a thick Demazure module via the twist by $\theta$, and the quotient is taken by a submodule that is a sum of other thick Demazure submodules.

Having established the existence of such a filtration by thick Demazure modules, we next identify its associated graded components. By the definition of $\bW_\la^{(k)}$ and Proposition~\ref{prop:ordercomp}, the only possible factors $\bD^{(k)}_\mu$ are those for which the weight $\mu \in P$ satisfies
\[
\la_+ \preceq \mu \preceq \la_-.
\]
These are precisely the weights in the $W$-orbit of $\la$. Furthermore, the corresponding $\h$-weight vectors lie in the submodule $V_\la \subset (\bW_\la^{(k)})_0$ (which is in fact an equality), where each such weight occurs with multiplicity one because it is extremal in the irreducible $\g$-module $V_\la$. These vectors serve as cyclic generators of the corresponding thick Demazure subquotients. This shows that each $\bD_\mu^{(k)}$ appears exactly once in the filtration, completing the proof.
\end{proof}

\begin{lem}\label{lem:DtoWbyD}
Let $k \in \Z_{>0}$ and $\la \in P^+$. Then we have an isomorphism
\[
\bW_\la^{(k)} \cong \mathbb L^{\bullet} \mathscr D_{w_0} ( \bD^{(k)}_{\la} )
\]
of graded $\tb$-modules.
\end{lem}

\begin{proof}
We first note from the definition of $\bW_\la^{(k)}$ and Corollary~\ref{cor:tD} that $\bW_\la^{(k)}$ is a quotient of the dual of the thick Demazure module $\mathscr D_{w_0} ( \widehat{\bD}^{(k)}_{\la} )$. In addition, Corollary~\ref{cor:tD} implies $\mathbb L^{<0} \mathscr D_{w_0} ( \widehat{\bD}^{(k)}_{\la} ) \equiv 0$. This gives a natural surjection:
\begin{equation}
\mathbb L^{\bullet} \mathscr D_{w_0} ( \widehat{\bD}^{(k)}_{\la} ) \equiv \mathscr D_{w_0} ( \widehat{\bD}^{(k)}_{\la} ) \twoheadrightarrow \bW_\la^{(k)}.\label{eqn:ssurj}
\end{equation}
Both sides of the surjection~\eqref{eqn:ssurj} admit filtrations by modules in the family $\{ \bD^{(k)}_{\mu} \}_{\mu \in P}$. A comparison of their degree-zero parts, along with Proposition~\ref{prop:WDfilt}, shows that the surjection must be compatible with these filtrations. Consequently, its kernel must also be filtered by modules of the form $\bD^{(k)}_{\mu}$ with $\la_- \prec \mu_-$.

To identify the kernel precisely, we first consider a related map. The $\theta$-twist of $\widehat{\bD}^{(k)}_{\la}$ has a $\tb^-$-cyclic vector of $\h$-weight $-\la$. By Proposition~\ref{prop:ordercomp} and Theorem~\ref{thm:D-incl}, the $\tb^-$-cyclic vector of any proper thick Demazure submodule has $\h$-weight strictly greater than $-\la$ with respect to $\prec$. It follows that the kernel
\[
\ker_0 := \ker ( \widehat{\bD}^{(k)}_\la \to \bD^{(k)}_\la )
\]
is the sum of all modules $\widehat{\bD}^{(k)}_\mu$ with $\la_- \prec \mu_-$ that are contained in $\widehat{\bD}^{(k)}_\la$.

By applying Corollaries~\ref{cor:D-dist} and~\ref{cor:tD}, we deduce that $\mathscr D_{w_0} (\ker_0)$ is the sum of all modules $\widehat{\bD}^{(k)}_\mu$ with $\la_- \prec \mu_-$ that are contained in $\widehat{\bD}^{(k)}_{\la_-}$ (and $\mathbb L^{<0} \mathscr D_{w_0} (\ker_0) \equiv 0$). This coincides with the definition of the submodule used in the construction of $\bW_\la^{(k)}$, showing that the surjection~\eqref{eqn:ssurj} is indeed an isomorphism.
\end{proof}

\subsection{The module $\bW_\la$ and its properties}

By definition, $W_\la^{(1)}$ and $\bW_\la^{(1)}$ (for $\la \in P^+$) are level-one modules on which the central element $K \in \tg$ acts by $+1$ and $-1$, respectively.

\begin{thm}[Fourier–Littelmann–Manning–Senesi~\cite{FL07,FMS13}]\label{thm:CKmain1}
For each $\la \in P^+$, the module $\bW_\la^{(1)} \otimes \C_{\La_0}$ is a nontrivial self-extension of $W_\la^{(1)} \otimes \C_{-\La_0}$ as a $\tg_{\ge 0}$-module. Moreover, its endomorphism ring
\[
\mathrm{End}_\fC ( \bW_\la^{(1)} \otimes \C_{\La_0} )
\]
is isomorphic to a graded polynomial ring generated by elements of strictly positive degree. In addition, $\bW_\la^{(1)} \otimes \C_{\La_0}$ is a projective object in the full subcategory of $\fC$ consisting of those modules $M$ with
\[
\overline{\Psi ( M )} \le \la,
\]
where $\le$ is the dominance order.
\end{thm}

The following result is a reformulation of~\cite[Theorem~5.2]{CK18}, adapted to the setting of~\cite[Theorem~4.13]{Kat18}.  
Note that there exists a $W$-equivariant inclusion $Q \cap P \subset Q^{\vee}$, which identifies short roots with long coroots.
This identification is specific to the twisted affinizations.

\begin{thm}[\cite{CK18, Kat18}]\label{thm:CKmain2}
Let $\la \in P^+$ and $\beta \in Q \cap P^+ \subset Q^{\vee}$. Then there is an isomorphism of graded $\tb$-modules:
\[
\mathbb L^{\bullet} \mathscr D_{t_{- \beta}} \left( \bW_\la ^{(1)} \otimes \C_{\La_0} \right)
\cong \bq^{- \langle \beta, \la \rangle} \, \bW_\la ^{(1)} \otimes \C_{\La_0}.
\]
\end{thm}

Let us extend the definition of $\bW_\la^{(1)} \otimes \C_{\La_0}$ to arbitrary $\la \in P$ in a slightly twisted manner:

\begin{prop}[{\cite[\S4.2.2]{CK18}}]\label{prop:Wzero}
For each $\la \in P$, there exists a $\wth$-semisimple $\tb$-module $\bW_\la$ that carries a $\tb$-cyclic vector of $\wth$-weight $\la$ with the following properties:
\begin{enumerate}
\item We have an isomorphism $\bW^{(1)}_{\la_+} \otimes \C_{\La_0} \cong \bW_{\la_-}$ as $\tb$-modules;
\item We have an isomorphism $\bD^{(1)}_{\la_+} \otimes \C_{\La_0} \cong \bW_{\la_+}$ as $\tb$-modules;
\item The module $\bW_\la$ admits a finite filtration with successive quotients isomorphic to $\bD^{(1)}_\mu \otimes \C_{\La_0}$, where $\mu$ ranges over
\[
\{ \mu \in W \la \mid \la \le \mu \le \la_+ \},
\]
each occurring exactly once;
\item For each $i \in \tI_\af$, we have
\[
\mathbb{L}^{\bullet} \mathscr{D}_i ( \bW_\la ) \cong \begin{cases} 
\bW_\la & \text{if } s_i (\!(\la)\!) \preceq \la, \\
\bq^{-\langle \vartheta^{\vee}, \la \rangle \delta_{i,0}} \bW_{\overline{s_i \la}} & \text{if } s_i (\!( \la)\!) \succ \la.
\end{cases}
\]
\item For each $i \in \tI_\af$ such that $s_i (\!(\la)\!) \preceq \la$, we have $\mathbb{L}^{\bullet} \mathscr{D}_i ( \bD_\la^{(1)} ) \equiv 0$. For each $i \in \tI$ such that $s_i (\!(\la)\!) \succ \la$, we have a short exact sequence
\[
0 \rightarrow \bD_\la^{(1)} \longrightarrow \mathscr{D}_i ( \bD_\la^{(1)} ) \longrightarrow \bD_{s_i (\!(\la)\!)}^{(1)} \longrightarrow 0.
\]
\end{enumerate}
\end{prop}

\begin{proof}
The cyclic vector is provided by the construction at the beginning of~\cite[\S4.2.2]{CK18}; the first item appears at the beginning of~\cite[\S4.2.2]{CK18}; the second is~\cite[Corollary~4.15]{CK18}; the third is~\cite[Lemma~4.13]{CK18}; the fourth ~\cite[Proposition~4.14]{CK18}, while the case $i = 0$ follows from~\cite[Theorem~5.2]{CK18} and its proof; the fifth is~\cite[Corollary~4.15]{CK18}.
\end{proof}

\begin{rem}
The family of modules $\{\bW_\la\}_{\la \in P}$ was first realized both as the duals of the global sections of line bundles on semi-infinite flag varieties~\cite{Kat18} and, simultaneously, as the non-symmetric counterpart of $q$-Whittaker functions for untwisted affinizations~\cite{FMO18}. Extensions to general twisted affinizations are given in~\cite{CK18}, and the case of type $\mathsf{A}^{(2)}_{2\ell}$ is treated in~\cite{Chi22}.
\end{rem}

Let $\Theta$ be an automorphism of the affine Dynkin diagram that permutes the vertices $\tI_\af$. Then $\Theta$ defines an automorphism of $\tb$ and $\tg$, which we denote by the same symbol $\Theta$, such that
\[
\Theta ( E_i ) = E_{\Theta(i)}, \Theta ( F_i ) = F_{\Theta(i)}, \Theta ( d ) = d \qquad i \in \tI_\af.
\]

\begin{cor}\label{cor:CKmain2}
Let $\Theta$ be an automorphism of $\tb$ induced by an affine Dynkin diagram automorphism. For each $\la \in P$, there exists $\la' \in P$ such that
\[
\Theta^* ( \bW_\la ) \cong \bW_{\la'}.
\]
\end{cor}

\begin{proof}
Combining Theorems~\ref{thm:CKmain1} and~\ref{thm:CKmain2} and Proposition~\ref{prop:Wzero}, we observe that the module
\[
\bigcup_{w \in W_\af} \mathscr D_w ( \bW_\la )
\]
satisfies the universal property of a level-zero extremal weight module of $\tg$ (see~\cite[\S8]{Kas93}).  
It follows that $\bW_\la$ must be a Demazure module of such a level-zero extremal weight module (see~\cite{Kas05}).  
Thus, up to grading shifts, we obtain a bijection between the set $\{\bW_{\la}\}_{\la \in P}$ and the set of Demazure submodules of level-zero extremal weight modules.

Since the structure of level-zero extremal weight modules and their Demazure submodules is preserved under diagram automorphisms, the same conclusion holds for the transformed module $\bW$.
\end{proof}

\begin{rem}
Let $\Theta$ be the automorphism of the affine Dynkin diagram that permutes the vertices $\tI_\af$. Then the weights $\La$ appearing in Corollary~\ref{cor:CKmain2} all satisfy the level zero condition $\La(K) = 0$. The automorphism $\Theta$ induces linear automorphisms of the spaces
\[
\{ \La \in \wth^* \mid \La(K) = 0 \} \subset \wth^*
\]
satisfying
\[
\Theta(\alpha_i) = \alpha_{\Theta(i)}, \quad \Theta(\Lambda_i) = \Lambda_{\Theta(i)}, \quad \Theta(\delta) = \delta \qquad (i \in \tI_\af).
\]

In~\S\ref{subsubsec:sf}, each fundamental weight $\varpi_i \in P$ ($i \in \tI$) is identified with
\[
\varpi_i = \Lambda_i - \langle \varpi_i, \vartheta^{\vee} \rangle \Lambda_0 = \Lambda_i - \Lambda_i(K) \Lambda_0 \in P_\af.
\]
By convention, we formally set $\varpi_0 := \Lambda_0 - \Lambda_0(K) \Lambda_0 = 0$, corresponding to the vertex $0 \in \tI_\af \setminus \tI$.

It then follows that
\[
\Theta(\varpi_i) = \Lambda_{\Theta(i)} - \Lambda_i(K) \Lambda_{\Theta(0)} \qquad (i \in \tI_\af),
\]
which satisfies the consistency condition
\[
\Theta(\varpi_0) = \Lambda_{\Theta(0)} - \Lambda_0(K) \Lambda_{\Theta(0)} = \Lambda_{\Theta(0)} - \Lambda_{\Theta(0)} = 0.
\]
Hence, applying $\Theta$ linearly yields $\la' = \Theta (\la)$ in Corollary~\ref{cor:CKmain2}.
\end{rem}

\subsection{Spaces of characters $\mathbf{B}$ and $\mathbf{C}$}

We now introduce the ambient spaces in which the graded characters of our modules naturally live.

We define the $\Q(\!(q)\!)$-vector space
\[
\mathbf B := \Q(\!(q)\!) \otimes_{\Z[\![q]\!]} \varprojlim_{m} \bigoplus_{\la \in P} \left( \Z[q] / (q^{m+1}) \right) e^\la.
\]
We define $\mathbf C$ to be the subspace of $W$-invariants in $\mathbf B$, where the $W$-action is given by $w \cdot (q^m e^\la) := q^m e^{w\la}$ for $w \in W$, $\la \in P$, and $m \in \Z$.  
Both vector spaces are equipped with the inverse limit topology induced by the $q$-adic filtration.

For each module $M$ in $\fB_\bd$ or $\fC_\bd$, we define its graded character as
\[
\gch M := \sum_{\la \in P,\ m \in \Z} q^m e^\la \cdot \dim \mathrm{Hom}_{(\h + \C d)} ( \C_{\la + m \delta},\ M ).
\]
In particular, this character is invariant under level shifts:
\[
\gch M = \gch ( M \otimes \C_{k \La_0} ) \quad \text{for any } k \in \Z.
\]

\begin{lem}
For each $M$ in $\fB_\bd$ or $\fC_\bd$, we have $\gch M \in \mathbf B$ or $\gch M \in \mathbf C$, respectively. 	
\end{lem}

\begin{proof}
This follows directly from the definitions of the categories and the grading structure of $M$.
\end{proof}

\begin{prop}[{\cite[Appendix]{FKM}}]\label{lem:ep}
For each $M, N \in \fB_{\bd}$, define the Euler–Poincaré pairing by
\[
\langle M, N \rangle_{\mathtt{EP}} := \sum_{i \in \Z} (-1)^i \gdim \mathrm{Ext}_{\fB}^i ( M, N^\vee ).
\]
This pairing extends uniquely to a $\Z[q^{\pm 1}]$-antibilinear map
\[
[\fB_{\bd}] \times [\fB_{\bd}] \longrightarrow \Q(\!(q^{-1})\!),
\]
where $[\mathcal C]$ denotes the Grothendieck group of the abelian category $\mathcal C$.
\end{prop}

\begin{lem}\label{lem:top}
Let $k \in \Z_{>0}$. We have the following:
\begin{enumerate}
\item The sets of characters $\{\gch D^{(k)}_\la \}_{\la \in P}$ and $\{\gch \bD^{(k)}_\la \}_{\la \in P}$ are topological $\Q(\!(q)\!)$-bases of $\mathbf B$.
\item The sets of characters $\{\gch W^{(k)}_\la \}_{\la \in P^+}$ and $\{\gch \bW^{(k)}_\la \}_{\la \in P^+}$ are topological $\Q(\!(q)\!)$-bases of $\mathbf C$.
\end{enumerate}
\end{lem}

\begin{proof}
Any integrable lowest-weight module of $\tg$ is a quotient of a parabolic Verma module $P_\la$ for some $\la \in P^+$ (see~\cite[\S9.2]{Kac}). The graded character of such a module belongs to $\mathbf C$, as we have
\[
\gch P_\la = \gch \bigl( U(\tn / \gn) \bigr) \cdot \ch V_\la \in \mathbf C.
\]
The modules $\bD^{(k)}_\mu$ (for $\mu \in W\la$) and $\bW^{(k)}_\la$ are subquotients of such modules. It follows that their graded characters are termwise bounded above by that of a parabolic Verma module, and hence belong to $\mathbf B$ or $\mathbf C$, respectively.

By Lemma~\ref{lem:headW}, both $W^{(k)}_\la$ and $\bW^{(k)}_\la$ contain $V_\la$ in degree zero and are generated by it as $\tg_{\ge 0}$-modules. Therefore, we have the congruences
\[
\gch W^{(k)}_\la \equiv \gch \bW^{(k)}_\la \equiv \gch V_\la \pmod{q}
\quad\text{for each $\la \in P^+$}.
\]

Furthermore, for each $\mu \in P$, the module $D^{(k)}_\mu$ is contained in $W^{(k)}_{\mu_+}$, and $\bD^{(k)}_{u \mu_+}$ appears in a filtration of $\bW^{(k)}_{\mu_+}$. Thus, we have
\[
\left. \gch D^{(k)}_\mu \right|_{q=0} \in \sum_{\la \in \Sigma(\mu)} \Z e^\la \ni \left. \gch \bD^{(k)}_\mu \right|_{q=0}.
\]
Here, the specialization $\left. \gch M \right|_{q=0}$ denotes the image of $\gch M$ under the reduction map $\Z[q] \to \Z[q]/(q) \cong \Z$ applied coefficientwise.

These specializations belong to $\bigoplus_{\la \in P} \Z e^\la$, and the resulting sets
\[
\Bigl\{\left. \gch D^{(k)}_\mu \right|_{q=0}\Bigr\}_{\mu \in P} \qquad \text{and} \qquad \Bigl\{\left. \gch \bD^{(k)}_\mu \right|_{q=0} \Bigr\}_{\mu \in P}
\]
form $\Z$-bases of this space. The same argument applies to $\{\gch W^{(k)}_\la\}_{\la \in P^+}$ and $\{\gch \bW^{(k)}_\la\}_{\la \in P^+}$.

Since all the characters involved are $q$-series with integer coefficients whose degrees are bounded below, it follows that any element of $\mathbf B$ or $\mathbf C$ can be uniquely expressed as a formal $q$-series in terms of these characters. This proves that the families listed in the statement form topological $\Q(\!(q)\!)$-bases.
\end{proof}

\subsection{$\Ext$-orthogonality in the level-one case}

\begin{thm}[{\cite[Theorem~4.19]{CK18}}]\label{thm:Dorth}
For each $\la, \mu \in P$, we have
\[
\ext_{\fB}^{i} \left( \bD^{(1)}_\la \otimes \C_{\La_0}, ( D^{(1)}_\mu \otimes \C_{-\La_0} )^* \right)
= \begin{cases} \C & \text{if } i = 0 \text{ and } \la + \mu = 0, \\
0 & \text{otherwise}. \end{cases}
\]
Moreover, the set $\Psi ( D^{(1)}_\la )$ of weights of $D^{(1)}_\la$ satisfies
\[
\overline{\Psi ( D^{(1)}_\la )} \subset \Sigma(\la).
\]
\end{thm}

\begin{lem}\label{lem:extDC}
For $\la, \mu \in P$, we have
\begin{equation}
\ext_{\fB}^{\bullet} \left( \bD^{(1)}_\la \otimes \C_{\La_0}, (\C_{\mu})^* \right) = 0 \quad \text{for all } - \la \not\preceq \mu.\label{eqn:ext-ch}
\end{equation}
\end{lem}

\begin{proof}
We prove the assertion by induction on $\mu$ with respect to $\prec$, starting from minimal elements.
In case $\mu$ is minimal, we have $\dim D^{(1)}_\mu = 1$ since $\mu + \La_0$ must be a level one dominant weight.

Assume~\eqref{eqn:ext-ch} holds for all $\nu \prec \mu$. Consider the short exact sequence
\[
0 \to \ker f_{\mu} \longrightarrow D^{(1)}_\mu \stackrel{f_{\mu}}{\longrightarrow} \C_\mu \to 0.
\]
By the inductive hypothesis, all weights $\gamma$ occurring in $\ker f_{\mu}$ satisfy $\gamma \prec \mu$, so we have
\[
\ext_{\fB}^{\bullet} \left( \bD^{(1)}_\la \otimes \C_{\La_0}, (\ker f_\mu \otimes \C_{-\La_0})^* \right) = 0.
\]
Moreover, Theorem~\ref{thm:Dorth} gives
\[
\ext_{\fB}^{\bullet} \left( \bD^{(1)}_\la \otimes \C_{\La_0}, (D^{(1)}_\mu \otimes \C_{-\La_0})^* \right) = 0.
\]
Applying the long exact sequence for $\ext_{\fB}^\bullet$ to the dual of the above short exact sequence (after twisting), we conclude that
\[
\ext_{\fB}^{\bullet} \left( \bD^{(1)}_\la \otimes \C_{\La_0}, (\C_{\mu})^* \right) = 0,
\]
which completes the induction step and proves~\eqref{eqn:ext-ch}.
\end{proof}

\subsection{Two technical lemmas}

\begin{prop}[{\cite[Lemma~4.4]{Kat18}; see also~\cite{Jos85,Mat88}}]\label{prop:string}
Assume that $\g = \mathfrak{sl}(2,\C)$. Let $\mathscr D$ denote the Demazure functor with respect to $\gb$, and let $M$ be a $\gb$-module with semisimple $\h$-action. Suppose that:
\begin{itemize}
  \item[(i)] $\mathbb{L}^{-1} \mathscr D(M) = 0$, and
  \item[(ii)] there exists a $\gb$-module embedding $M \hookrightarrow \mathscr D(M)$.
\end{itemize}
Then $M$ admits a finite filtration whose associated graded is a direct sum of modules of the following two types:
\begin{enumerate}
  \item irreducible finite-dimensional $\mathfrak{sl}(2,\C)$-modules;
  \item one-dimensional $\gb$-modules $\C_\mu$ with $\langle \alpha^\vee, \mu \rangle \ge 0$, where $\alpha^{\vee}$ denotes the unique positive coroot of $\mathfrak{sl}(2,\C)$.
\end{enumerate}
\end{prop}

\begin{proof}
Since $\mathbb{L}^{-1} \mathscr D(M) = 0$, the derived Demazure functor coincides with its underived version in degree zero. In particular, we have
\[
\gch \mathscr D(M) = \partial \gch M,
\]
where $\partial$ denotes the Demazure operator associated with $\mathfrak{sl}(2,\C)$. Therefore, the assumptions of~\cite[Lemma~4.4]{Kat18} are satisfied, and the desired conclusion follows.
\end{proof}

\begin{lem}\label{lem:app}
Let $i \in \tI_\af$, and let $M$ be a finite-dimensional
$\wth$-semisimple indecomposable $\C E_i \oplus \C \al_i^{\vee}$-module which satisfies the
assumptions of Proposition~\ref{prop:string}.
Suppose that $v \in M$ is an $\al_i^{\vee}$-eigenvector with
eigenvalue $t$. Then
\[
v=\sum_{\substack{|t|\le t'\\ t'\equiv t\!\!\!\pmod 2}}v_{t'},
\]
where the sum ranges over the highest weights $t'\varpi$ of the
irreducible $\mathfrak{sl}(2,i)$-summands of $\mathscr D_i(M)$, and
each $v_{t'}$ is a nonzero $\al_i^{\vee}$-eigenvector with eigenvalue
$t$ in the corresponding summand.
\end{lem}

\begin{proof}
Since $M$ is finite-dimensional and indecomposable over
$\C E_i \oplus \C\al_i^\vee$, it is generated by a lowest-weight
vector $u$, and $v$ is a nonzero scalar multiple of $E_i^r u$ for
some $r\ge0$.

By the definition of the Demazure functor, $\mathscr D_i(M)$ is
generated by $M$. Hence the projection of $u$ to every irreducible
$\mathfrak{sl}(2,i)$-summand of $\mathscr D_i(M)$ is nonzero.
In a summand of highest weight $t'\varpi$, inspection of the
$\mathfrak{sl}(2)$-weight string shows that the corresponding
component of $E_i^r u$ is nonzero precisely when
\[
|t|\le t'
\qquad\text{and}\qquad
t'\equiv t\pmod 2;
\]
once the resulting weight exceeds the highest weight $t'$, that
component is zero. The assertion follows.
\end{proof}

\subsection{An estimate of weights in $D_\la^{(k)}$}

\begin{lem}\label{lem:we}
For each $\la \in P$ and $k \in \Z_{\ge 2}$, the module $D_\la^{(k)}$ is a quotient of $D_\la^{(k-1)} \otimes \C_{\La_0}$ in $\fB$. Moreover, we have
\[
\overline{\Psi ( D_\la ^{(k)})} \subset \overline{\Psi ( D_{\la}^{(1)} )} \subset \Sigma (\la).
\]
\end{lem}

\begin{proof}
By~\cite[3.4~Theorem]{Jos85}, the defining relations of $D_{\la}^{(k)}$ are given by
\begin{align}\nonumber
e_{\al}^{\max \{ - \langle \al^{\vee}, \la \rangle, 0 \} + 1} \bv &= 0 && \text{for } \al \in \Delta_+,\\
e_{\al + m \delta}^{\max\{ - \langle \al^{\vee}, \la \rangle - 2 mk / r_\al, 0\} + 1} \bv &= 0 && \text{for } \al \in \Delta,\ m \in \Z_{>0},\label{eqn:defD}
\end{align}
where $\bv$ is the cyclic vector of $\h$-weight $\la$, and $r_\al$ is the squared root length.

From this description, it is clear that $D_\la^{(k)}$ is a quotient of $D_\la^{(k-1)} \otimes \C_{\La_0}$. In particular, we obtain the inclusion
\[
\Psi ( D_\la^{(k)} ) \subset \Psi ( D_{\la}^{(1)} ).
\]

The further inclusion $\Psi ( D_{\la}^{(1)} ) \subset \Sigma(\la)$ is a part of Theorem~\ref{thm:Dorth}, and it also follows from the identification of the graded character of $D_\la^{(1)}$ with the non-symmetric Macdonald polynomial $E_\la(q,0)$ (see~\cite{Ion03}), together with the weight support estimate in~\cite[(4.4)]{Che95}, which is part of the characterization of these polynomials.
\end{proof}

\begin{rem}\label{rem:thin}
Joseph's result~\cite{Jos85} is stated only for semisimple Lie algebras $\g$, but the construction and arguments extend without modification to thin Demazure modules for general symmetrizable Kac–Moody algebras. This is because the proof in~\cite[\S4]{Jos85} proceeds by induction on the length of the sequence of simple reflections, starting from the one-dimensional case corresponding to the empty sequence. In contrast, the argument does not generalize to thick Demazure modules.
\end{rem}

\section{On the structure of affine Demazure modules}\label{sec:tmodule}

We retain the setting of the previous section.

In this section, we establish a technical but crucial structural result concerning affine Demazure modules, stated as Proposition~\ref{prop:Demext}, referred to as {\it Reduction} in the Introduction. Although the proof is necessarily case-by-case and reflects features specific to twisted affinizations, this result plays an essential role in the inductive proof of the $\ext^1$-vanishing theorem (Theorem~\ref{thm:Lext}) by enabling a reduction from level $k$ to level $k{-}1$.

More precisely, in the argument for Theorem~\ref{thm:Lext}, the $\ext^1$-vanishing for weights $\La \in P^+_{\af, k}$ is formally reduced to that for weights of the form $\La - \La_0 \in P_{\af, (k-1)}$. Proposition~\ref{prop:Demext} enables a further key reduction by replacing this with a vanishing condition for another dominant weight $\La' \in P^+_{\af, (k-1)}$, in a way that is compatible with the Demazure module structure.

\smallskip
For $k \in \Z_{> 0}$ and $\La \in P^+_{\af, k}$, we set
\[
S^+ ( \La ) := \{ j \in \tI_\af \mid \langle \al_j^{\vee}, \La \rangle = 0 \},
\]
and let $S ( \La ) \subset S^+ ( \La )$ be the maximal subset that forms a connected subdiagram of the Dynkin diagram of $\tg$ containing $0$. (We have $S (\La) = \emptyset$ if and only if $\langle \vartheta^{\vee}, \La \rangle < k$.) 
Since $k > 0$, we have $S(\La) \subsetneq \tI_\af$.  
Let $W_{S(\La)} := \left< s_j \mid j \in S ( \La ) \right>$, and let $\g_{S(\La)}$ be the simple Lie subalgebra of $\tg$ generated by $\{ E_j, \al_j^{\vee}, F_j \}_{j \in S(\La)}$.

\begin{lem}\label{lem:mini}
Let $k \in \Z_{\ge 2}$ and $\La \in P_{\af, k}^+$. Then the weight $\La - \La_0$ is a minuscule antidominant weight of $\g_{S(\La)}$. In particular, $\g_{S(\La)}$ is not of types $\mathsf{E}_8$, $\mathsf{F}_4$, or $\mathsf{G}_2$.
\end{lem}

\begin{proof}
The weight $\La_0$, viewed as a weight of $\g_{S(\La)}$, is the fundamental weight corresponding to the Dynkin vertex $0$.  
Since $S(\La) \subsetneq \tI_\af$ and $\tg$ is a twisted affinization, a case-by-case inspection using the classification of twisted affine Lie algebras (\cite[Chap.~4, Table~Aff]{Kac}) and that of minuscule weights (\cite[V\!I\!I\!I, \S7.3]{Bou08}) shows that $\La_0$ is a minuscule dominant weight of $\g_{S(\La)}$.  
Since $\La$ is trivial on the coroots of $\g_{S(\La)}$, the restriction of $\La-\La_0$ to $\g_{S(\La)}$ is the negative of this minuscule dominant weight. This establishes the first assertion.

The second assertion follows from the fact that the simple Lie algebras of types $\mathsf{E}_8$, $\mathsf{F}_4$, and $\mathsf{G}_2$ admit no nonzero minuscule weights (see the above classification).
\end{proof}

\begin{lem}\label{lem:coreDem}
Let $k \in \Z_{\ge 2}$. Suppose that $\La \in P_{\af, k}^+$ and $w \in W_{S(\La)}$ satisfy $\langle \vartheta^{\vee}, \overline{\La} \rangle = k$ and $w^{-1} ( \La - \La_0 ) \in P^+_{\af, (k-1)}$. In addition, assume that $w$ is the minimal element with this property. Then $\ell ( s_0 w ) < \ell ( w )$, and we have
\[
\mathscr D_{w} ( \C_{w^{-1} (\La - \La_0)} ) \big/ \mathscr D_{s_0 w} ( \C_{w^{-1} (\La - \La_0)} ) \cong \C_{\La - \La_0}.
\]
\end{lem}

\begin{proof}
Since $w^{-1} ( \La - \La_0 ) \in P^+_{\af,(k-1)}$ and
\[
s_0 w w^{-1} ( \La - \La_0 ) 
= (\La - \La_0) + \al_0 
= w w^{-1} ( \La - \La_0 ) + \al_0,
\]
we deduce that $s_0 w < w$, which implies $\ell ( s_0 w ) < \ell ( w )$. This proves the first assertion.

Let $S:=S(\La)$. By Lemma~\ref{lem:mini}, the module $\mathscr D_w(\C_{w^{-1}(\La-\La_0)})$ is isomorphic to the irreducible highest-weight $\mathfrak g_S$-module of minuscule highest weight $w^{-1}(\La-\La_0)$, so all of its weight spaces are one-dimensional and indexed by the $W_S$-orbit of $\La-\La_0$.
Thus, the second assertion holds if every such $v ( \La - \La_0 ) \neq \La - \La_0$ yields a nonzero weight space in $\mathscr D_{s_0 w} ( \C_{w^{-1} (\La - \La_0)} )$.

Since $\La - \La_0$ is antidominant for $\g_S$, there exists a minimal sequence
\[
i_1, i_2, \ldots, i_t \in S
\]
such that
\[
(s_{i_1} s_{i_2} \cdots s_{i_t}) v ( \La - \La_0 ) = \La - \La_0,
\quad\text{provided } v ( \La - \La_0 ) \neq \La - \La_0.
\]
By the minimality of $t$, we must have $i_1 = 0$, so that
\[
(s_{i_2} \cdots s_{i_t}) v ( \La - \La_0 ) = s_0( \La - \La_0 ).
\]

Moreover, we may replace $v$ with $(s_{i_1} s_{i_2} \cdots s_{i_t})^{-1}$ without loss of generality since only the weight $v(\Lambda-\Lambda_0)$ matters.  
Then, by the minimality assumption on $w$, we have $\ell ( v w ) = \ell ( w ) - \ell ( v )$, and
\[
s_{i_2} \cdots s_{i_t} \cdot v w = s_0 w.
\]
Therefore,
\[
\C_{v ( \La - \La_0 )} \subset 
( \mathscr D_{i_2} \circ \cdots \circ \mathscr D_{i_t} ) 
\bigl( \mathscr D_{v w} ( \C_{w^{-1} ( \La - \La_0 )} ) \bigr)  
= \mathscr D_{s_0 w} ( \C_{w^{-1} ( \La - \La_0 )} ).
\]
Hence, each such $v ( \La - \La_0 )$ appears as a weight of $\mathscr D_{s_0 w} ( \C_{w^{-1} (\La - \La_0)} )$, as desired.
\end{proof}

\begin{prop}\label{prop:Demext}
Let $k \in \Z_{\ge 2}$. For each $\La \in P_{\af, k}^+$ satisfying $\langle \vartheta^{\vee}, \overline{\La} \rangle = k$, there exist $w \in W_\af$ and $\La' \in P_{\af,(k-1)}^+$ such that
\begin{equation}
\mathscr D_w ( \C_{\La'} ) \big/ \mathscr D_{s_0 w} ( \C_{\La'} ) \cong \C_{ \La - \La_0 }.
\label{eqn:special}
\end{equation}
\end{prop}

\begin{proof}
We have
\[
s_0 ( \La - \La_0 ) = \La - \La_0 + \al_0,
\]
since $\langle \al_0^{\vee}, \La - \La_0 \rangle = 0 - 1 = -1$.  
If $s_0 ( \La - \La_0 ) \in P_{\af,(k-1)}^+$, then
\[
\mathscr D_0 ( \C_{s_0 ( \La - \La_0 )} ) \cong \C_{\La - \La_0} \oplus \C_{s_0 ( \La - \La_0 )}.
\]
Thus, the assertion holds in this case.

This includes, in particular, the case when $\g$ is of type $\mathsf{A}_1$.  
In the remainder of the proof, we assume that $\g$ is not of type $\mathsf{A}_1$.

\medskip

We now assume that $s_0 ( \La - \La_0 ) \not\in P_{\af,(k-1)}^+$.  
Thus, the possible values of $j \in \tI_\af$ satisfying $\langle \al_j^{\vee}, s_0 ( \La - \La_0 ) \rangle < 0$ are as follows.

If $\g$ is of type $\mathsf{A}_n$, then $j = 1$ or $j = n$ (cf.~\cite[p.~206 (VI)]{Bou02}, with $n = \ell$).  
In all other types, $j = i$ for a uniquely determined $i$, by inspection (cf.~\cite[Table~Aff]{Kac}).

Suppose that $\langle \al_j^{\vee}, \La \rangle > 0$ for all such possible $j$.  
Then necessarily $s_0 ( \La - \La_0 ) \in P_{\af,(k-1)}^+$, contradicting our assumption.

Hence, we must have $\langle \al_i^{\vee}, \La \rangle = 0$ for at least one such $i$:  
specifically, $i = 1$ or $i = n$ in type $\mathsf{A}_n$, and a uniquely determined $i$ in all other types.

We write $S$ for $S(\La)$ for the remainder of the proof.

By Lemma~\ref{lem:mini}, the weight $\La - \La_0$ is an antidominant minuscule weight (corresponding to the vertex $0$) of $\g_S$, and the subdiagram that appears in this construction is of type $\mathsf{A}$, $\mathsf{B}$, $\mathsf{C}$, $\mathsf{D}$, $\mathsf{E}_6$, or $\mathsf{E}_7$.

Hence, there exists a minimal length element $w \in W_S$ such that $w^{-1} ( \La - \La_0 )$ is a dominant minuscule (fundamental) weight of $\g_S$.

If $w^{-1} ( \La - \La_0 ) \in P_{\af,(k-1)}^+$, then the assertion follows from Lemma~\ref{lem:coreDem}.

\medskip

We now verify whether $w^{-1} ( \La - \La_0 ) \in P_{\af,(k-1)}^+$ case by case, according to the type of the Dynkin subdiagram determined by $S$.
In each case, we compute $w^{-1}(\Lambda-\Lambda_0)-(\Lambda-\Lambda_0)$ and check that its pairing with every simple coroot outside $S$ remains nonnegative because $k\ge 2$.

We first remark that if $j \in \tI_\af \setminus S$ is adjacent to a vertex of $S$ in the extended Dynkin diagram, then $\langle \al_j^{\vee}, \La \rangle \ge 1$.  
Moreover, if $|\tI_\af| - |S| = 1$ and
\begin{equation}
\langle \vartheta^{\vee}, \varpi_j \rangle = 1 
\quad \Leftrightarrow \quad 
\vartheta^{\vee} \in \al_j^{\vee} + \sum_{r \in \tI,\, r \neq j} \Z_{\ge 0} \al_r^{\vee},
\label{eqn:hrmult}
\end{equation}
then we obtain
\[
\langle \al_j^{\vee}, \La \rangle = \langle \vartheta^{\vee}, \La \rangle = k \ge 2.
\]

The condition~\eqref{eqn:hrmult} holds automatically when $\g$ is of type $\mathsf{A}$, or when $\g$ is of type $\mathsf{B}$, $\mathsf{C}$, or $\mathsf{D}$ and $j \in \tI$ corresponds to a short root connected to a unique other vertex in the (non-extended) Dynkin diagram, by inspection.

To distinguish between the original Dynkin indices and their rearrangement to match the conventions of~\cite{Bou02}, we may denote the key original Dynkin indices—namely $0$ (the affine vertex) and $i$ (a vertex adjacent to $0$)—by $\underline{0}$ and $\underline{i}$, respectively.

\paragraph{Type $\mathsf{A}_l$: $\underline{0}$ has degree one in $S$}  
We then label the vertices as
\[
S = \{1, 2, \ldots, l \} \quad \text{with} \quad 1 = \underline{0}, \quad 2 = \underline{i},
\]
so that $j$ and $j{-}1$ are adjacent in $S$.

In this case, the weight $\La - \La_0$ is the lowest weight of the vector representation of type $\mathsf{A}_l$.  
We then have
\[
w^{-1} ( \La - \La_0 ) - ( \La - \La_0 ) = \sum_{j = 1}^{l} \al_j.
\]

For each $j \in S$ and $t \in \tI_\af \setminus S$, we have $\langle \al_t^{\vee}, \al_j \rangle \in \{0, -1\}$.  
If $t \in \tI_\af \setminus S$ is adjacent to two simple roots in $S$ (which is the maximal possible), then necessarily $\g$ is of type $\mathsf{A}_n$ and $|\tI_\af| - |S| = 1$.

Therefore, since $k \ge 2$, we conclude that $w^{-1} ( \La - \La_0 ) \in P_{\af,(k-1)}^+$.

\paragraph{Type $\mathsf{A}_l$: $\underline{0}$ has degree two in $S$}  
We then label the vertices as
\[
S = \{-r_1, -1, 0, 1, 2, \ldots, r_2\},
\]
so that $j$ and $j{-}1$ are adjacent for each $j \in S$.

In this case, the weight $\La - \La_0$ corresponds to the minuscule representation of type $\mathsf{A}_{r_1 + r_2 + 1}$, and $\g$ is necessarily of type $\mathsf{A}_n$.  
We have
\[
w^{-1} ( \La - \La_0 ) - ( \La - \La_0 ) = \al_{-r_1} + 2 \al_{-r_1 + 1} + \cdots + \al_{r_2}.
\]

Each $j \in \tI_\af \setminus S$ that is adjacent to $S$ in the extended Dynkin diagram is adjacent either to $-r_1$ or to $r_2$.  
Therefore, since $k \ge 2$, we conclude that $w^{-1} ( \La - \La_0 ) \in P_{\af,(k-1)}^+$.

\paragraph{Type $\mathsf{D}_l$: $\underline{i}$ has degree three in $S$}  
We label the vertices as
\[
S = \{1, 2, \ldots, l\} \quad \text{with} \quad l = \underline{0}, \quad l{-}2 = \underline{i},
\]
so that $j$ and $j{-}1$ are adjacent for $j < l$, and $l$ is adjacent to $l{-}2$.

This situation occurs only when $\g$ is of type $\mathsf{C}$ or $\mathsf{D}$.  
In this case, the weight $\La - \La_0$ is the lowest weight of a half-spin representation.  
According to \cite[p.~209~(VI)]{Bou02}, we have
\begin{equation}
w^{-1} ( \La - \La_0 ) - ( \La - \La_0 ) = \varpi_l + \varpi_{l-t} = \al_{1} + 2 \al_{2} + \cdots,
\label{eqn:spin}
\end{equation}
where $t = 1$ when $l$ is odd and $t =0$ when $l$ is even.

A vertex in $\tI_\af \setminus S$ can be adjacent only to $1$ or $2$ in the extended Dynkin diagram.  
The latter case occurs only when $|\tI_\af| - |S| = 1$ and $\g$ is of type $\mathsf{D}$.  
Moreover, $\al_1$ must be a short root.

Therefore, since $k \ge 2$, we conclude that $w^{-1} ( \La - \La_0 ) \in P_{\af,(k-1)}^+$.

\paragraph{Type $\mathsf{D}_l$: $\underline{i}$ has degree two in $S$}  
We then label the vertices as
\[
S = \{1, 2, \ldots, l\} \quad \text{with} \quad 1 = \underline{0}, \quad 2 = \underline{i},
\]
so that $j$ and $j{-}1$ are adjacent for $j < l$, and $(l{-}2)$ is adjacent to $l$.

This situation occurs only when $\g$ is of type $\mathsf{D}$ or $\mathsf{E}$.  
In this case, the weight $\La - \La_0$ is the lowest weight of a minuscule representation.  
According to \cite[p.~209~(VI)]{Bou02}, we have
\[
w^{-1} ( \La - \La_0 ) - ( \La - \La_0 ) = 2 \varpi_1 = 2 \al_{1} + 2 \al_{2} + \cdots + \al_{l-1} + \al_{l}.
\]

If $\g$ is of type $\mathsf{D}_n$, then $|\tI_\af| - |S| = 1$, and the unique vertex in $\tI_\af \setminus S$ is connected to vertex $2$ in the extended Dynkin diagram; the corresponding root lengths are equal.  
If $\g$ is of type $\mathsf{E}$, then all external vertices are connected only to $l{-}1$ or $l$.

In either case, we conclude that $w^{-1} ( \La - \La_0 ) \in P_{\af,(k-1)}^+$.

\paragraph{Type $\mathsf{B}_l$ case}
We label the vertices as
\[
S = \{1, 2, \ldots, l\} \quad \text{with} \quad l = \underline{0}, \quad l{-}1 = \underline{i},
\]
so that the vertices $j$ and $j{-}1$ are adjacent in $S$ for each $j > 1$, in accordance with \cite[p.~202~(IV)]{Bou02}.

In this case, the vertex $l$ corresponds to the short root, and this situation occurs only when $\g$ is of type $\mathsf{B}_n$.  
Moreover, the weight $\La - \La_0$ is the lowest weight of a minuscule representation of $\g_S$.  
Thus, we have
\[
w^{-1} ( \La - \La_0 ) - ( \La - \La_0 ) = 2 \varpi_l = \al_1 + 2 \al_2 + \cdots.
\]

A vertex in $\tI_\af \setminus S$ can be connected only to vertex $1$ in the extended Dynkin diagram.  
The corresponding roots have different lengths only when $|\tI_\af| - |S| = 1$, in which case $\al_1$ is the long root.

Therefore, since $k \ge 2$, we conclude that $w^{-1} ( \La - \La_0 ) \in P_{\af,(k-1)}^+$.

\paragraph{Type $\mathsf{C}_l$ case}  
We label the vertices as
\[
S = \{1, 2, \ldots, l\} \quad \text{with} \quad 1 = \underline{0}, \quad 2 = \underline{i},
\]
so that the vertices $j$ and $j{-}1$ are adjacent in $S$ for all $j > 1$, and $l$ corresponds to the long root.  
(This occurs only when $\g$ is of type $\mathsf{C}$ or $\mathsf{F}$.)

In this case, the weight $\La - \La_0$ is the lowest weight of the vector representation of $\g_S$.  
Thus, we have
\[
w^{-1} ( \La - \La_0 ) - ( \La - \La_0 ) = \al_l + \sum_{j = 1}^{l - 1} 2 \al_j.
\]

If $l$ is connected to some $t \in \tI_\af \setminus S$ in the extended Dynkin diagram, then $\al_l$ and $\al_t$ have the same length.  
If another vertex $j \in S$, with $j \ne l$, is connected to $t \in \tI_\af \setminus S$, then $\g$ must be of type $\mathsf{C}_n$, $|\tI_\af| - |S| = 1$, and $\al_j$ and $\al_t$ again have the same length.

Therefore, since $k \ge 2$, we conclude that $w^{-1} ( \La - \La_0 ) \in P_{\af,(k-1)}^+$.

\paragraph{Type $\mathsf{E}_6$ case}  
Then $\g$ must also be of type $\mathsf{E}_6$, and in particular we have $|\tI_\af| - |S| = 1$.

We label the vertices as
\[
S = \{1, 2, \ldots, 6\} \quad \text{with} \quad 1 = \underline{0}, \quad 3 = \underline{i},
\]
in accordance with \cite[p.~218~(IV)]{Bou02}, where vertex $4$ is the trivalent vertex.

In this case, the weight $\La - \La_0$ is the lowest weight of the minuscule representation.  
Thus, we have
\[
w^{-1} ( \La - \La_0 ) - ( \La - \La_0 ) = \varpi_1 + \varpi_6 = 2 \al_1 + 2 \al_2 + \cdots.
\]

Moreover, vertex $2$ is the only vertex in $S$ that is adjacent to the unique vertex $j \in \tI_\af \setminus S$ in the extended Dynkin diagram.  
We have
\[
\langle \vartheta^{\vee}, \varpi_j \rangle = 1
\]
by \cite[p.~218~(IV)]{Bou02}.

Therefore, since $k \ge 2$, we conclude that $w^{-1} ( \La - \La_0 ) \in P_{\af,(k-1)}^+$.

\paragraph{Type $\mathsf{E}_7$ case}  
Then $\g$ must also be of type $\mathsf{E}_7$, and in particular we have $|\tI_\af| - |S| = 1$.

We label the vertices as
\[
S = \{1, 2, \ldots, 7\} \quad \text{with} \quad 7 = \underline{0}, \quad 6 = \underline{i},
\]
in accordance with \cite[p.~216~(IV)]{Bou02}. Here $4$ is again the trivalent vertex.

In this case, the weight $\La - \La_0$ is the lowest weight of the minuscule representation.  
We have
\[
w^{-1} ( \La - \La_0 ) - ( \La - \La_0 ) = 2 \varpi_7 = 2 \al_1 + 3 \al_2 + \cdots.
\]

Moreover, vertex $1$ is the only vertex in $S$ that is adjacent to the unique vertex in $\tI_\af \setminus S$ in the extended Dynkin diagram.  
We have
\[
\langle \vartheta^{\vee}, \varpi_j \rangle = 1,
\]
as recorded in \cite[p.~216~(IV)]{Bou02}.

Therefore, since $k \ge 2$, we conclude that $w^{-1} ( \La - \La_0 ) \in P_{\af,(k-1)}^+$.

\smallskip

This completes the proof.
\end{proof}

\section{A lifting theorem}\label{sec:lift}
We continue to work in the setting of the previous section.

\begin{thm}\label{thm:indsurj}
Let $M, N \in \fB_\bd$ be modules on which $K$ acts by $-k$ and $k$, respectively. Suppose that for some $i \in \tI_\af$, we have
\[
M \subset \mathscr D_i(M), \qquad N \subset \mathscr D_i(N), \qquad \text{and} \qquad \mathbb L^{<0} \mathscr D_i(M) = \mathbb L^{<0} \mathscr D_i(N) = 0.
\]
Assume further that $N$ admits a $\tb$-cyclic $\wth$-eigenvector $\bv$. Then the natural map
\[
\shom_{\fB} ( \mathscr D_i(M), N^{\vee} ) \longrightarrow \shom_{\fB} ( M, N^{\vee} )
\]
is surjective.
\end{thm}

\begin{rem}
Theorem~\ref{thm:indsurj} extends naturally to the setting of an arbitrary symmetrizable Kac–Moody algebra.
\end{rem}

\begin{proof}[Proof of Theorem~\ref{thm:indsurj}]
If either $\mathscr D_i(M) \cong M$ or $\mathscr D_i(N) \cong N$, then the map is an isomorphism by Theorem~\ref{thm:adj} together with the idempotence $\mathscr D_i^2 \cong \mathscr D_i$. Thus, the assertion follows in this case.

By the transpose isomorphism~\eqref{eqn:transpose}, it suffices to prove the surjectivity of the dual map
\begin{equation}
\shom_{\fB} ( N, \mathscr D_i(M)^{\vee} ) \longrightarrow \shom_{\fB} ( N, M^{\vee} ).\label{eqn:Lsurj}
\end{equation}

By Proposition~\ref{prop:string} and the isomorphism $U(\g_{\al_i}) \cong \C[E_i]$, the module $N$ decomposes as a direct sum of indecomposable $(\C E_i + \wth)$-modules, each of which is isomorphic to an irreducible $\mathfrak{sl}(2,i)$-module tensored with a one-dimensional $\wth$-module of weight $\la$ satisfying $\langle \al_i^{\vee}, \la \rangle \ge 0$. (This is essentially a Smith normal form argument for the action of $E_i$.) 
Similarly, $M$ admits a decomposition into a direct sum of modules of the same type.

Consider a nonzero morphism of $\tb$-modules
\[
f : N \longrightarrow M^{\vee}.
\]
Let $B$ be an indecomposable direct summand of $N$ as a $(\C E_i + \wth)$-module, viewed as a module over $\C E_i \oplus \C \al_i^{\vee}$. Then $B$ has head $\C_{r\varpi}$ and socle $\C_{s\varpi}$ for some integers $r \le s$.

Let $L$ be an indecomposable direct summand of $M$ as a $(\C E_i + \wth)$-module, regarded also as a $(\C E_i \oplus \C \al_i^{\vee})$-module. Then $L$ is isomorphic to an irreducible $\mathfrak{sl}(2,i)$-module tensored with a one-dimensional $(\wth + \C E_i)$-module $\C_{\la}$, on which $\wth$ acts by $\la$ with $\langle \al_i^{\vee}, \la \rangle \ge 0$. Consequently, $L^* \subset M^{\vee}$ is obtained by tensoring its dual $\C_{-\la}$ with an irreducible $\mathfrak{sl}(2,i)$-module.

Let $f_L$ denote the composition
\[
N \longrightarrow M^{\vee} \longrightarrow L^{*}.
\]

The image $f_L(B)$ is either zero or a quotient of $B$ as a $(\C E_i \oplus \C \al_i^{\vee})$-module. In the latter case, it has head $\C_{r\varpi}$ and socle $\C_{s'\varpi}$ for some $s' \le s$.

Using the natural embedding $N \subset \mathscr D_i(N)$, we consider the $\mathfrak{sl}(2,i)$-submodule $\widetilde{B} \subset \mathscr D_i(N)$ generated by $B$.  
By the defining property of $\mathscr D_i$, we have $\mathscr D_i(B) = \widetilde{B}$.  
It follows that $\widetilde{B}$ decomposes as a direct sum of irreducible $\mathfrak{sl}(2,i)$-modules with highest weights
\begin{equation}
s\varpi,\ (s{-}2)\varpi,\ \ldots,\ (|r|{+}2)\varpi,\ |r|\varpi. \label{eqn:listhw}
\end{equation}

We now examine whether the map $f_L$ extends from $B$ to $\widetilde{B}$.

\smallskip

If $s' > |r|$, then there exists a unique irreducible $\mathfrak{sl}(2,i)$-direct summand of $\widetilde{B}$ with highest weight $s'\varpi$ that admits an injection into $L^*$. This summand contains the image $f_L(B)$, since $s'\varpi$ (and hence $-s'\varpi$) is a weight of $L^*$.

If $s' \le |r|$, then there exists a unique irreducible $\mathfrak{sl}(2,i)$-direct summand of $\widetilde{B}$ with highest weight $|r|\varpi$ that surjects onto $f_L(B) \subset L^*$.

\smallskip

In either case, we obtain a nonzero map $\widetilde{f}_L$ from $\widetilde{B}$ to $L^*$.

By Lemma~\ref{lem:app}, the map $\widetilde{f}_L$ can be chosen uniquely to satisfy $\widetilde{f}_L |_B = f_L$.

Now set $B_0 := U(\C E_i)\, \bv$. This is a direct summand of $N$ as a $(\C E_i + \wth)$-module, by $\wth$-weight considerations ($B_0$ is precisely the sum of all $\wth$-weight spaces whose weights belong to $\La + \Z_{\ge 0} \al_i$ if $\bv$ has $\wth$-weight $\La$). In particular, its $\mathfrak{sl}(2,i)$-span $\widetilde{B}_0 \subset \mathscr D_i(N)$ admits a lift
\[
\widetilde{f} : \widetilde{B}_0 \longrightarrow M^{\vee}
\]
of $(\C E_i + \wth)$-modules such that its restriction to $B_0$ coincides with $f|_{B_0}$.

Let
\[
\mathfrak{u} := \bigoplus_{\beta \in \Delta_\af^+ \setminus \{\al_i\}} \tg_\beta \subset \tn.
\]
Since $U(\tn) \cong U(\mathfrak{u}) \otimes_{\C} U(\C E_i)$, we can regard
\[
U(\mathfrak{u}) \otimes_\C B_0 \qquad \text{and} \qquad U(\mathfrak{u}) \otimes_\C \widetilde{B}_0
\]
as $\tb$-modules, which are projective as $U(\mathfrak{u})$-modules. Moreover, $\mathfrak{u}$ is a $\mathfrak{sl}(2,i)$-integrable $\tp_i$-module by inspection.

Therefore, we have
\[
U(\mathfrak{u}) \otimes_\C \widetilde{B}_0 \;\cong\; U(\mathfrak{u}) \otimes_\C \mathscr D_i(B_0) \;\cong\; \mathscr D_i(U(\mathfrak{u}) \otimes_\C B_0)
\]
by Theorem~\ref{thm:DJ}(4).
In particular, we obtain surjections of $\tb$-modules
\[
U(\mathfrak{u}) \otimes_\C B_0 \longrightarrow N, \qquad \text{and} \qquad U(\mathfrak{u}) \otimes_\C \widetilde{B}_0 \longrightarrow \mathscr D_i(N).
\]

Set
\[
I := \ker\bigl( U(\mathfrak{u}) \otimes_\C B_0 \longrightarrow N \bigr).
\]
By assumption, we have $\mathbb L^{-1} \mathscr D_i(N) = 0$. It follows that the exact sequence
\[
0 = \mathbb L^{-1} \mathscr D_i(N) \longrightarrow \mathscr D_i(I) \longrightarrow U(\mathfrak{u}) \otimes_\C \widetilde{B}_0 \longrightarrow \mathscr D_i(N) \longrightarrow 0
\]
is short exact.

By Theorem~\ref{thm:adj}, giving a $\tb$-module map
\[
\widetilde{f} : \mathscr D_i(N) \longrightarrow M^{\vee}
\]
is equivalent to giving a $\tb$-module map
\[
\widetilde{f}' : \mathscr D_i(N) \longrightarrow \mathscr D_i(M)^{\vee}.
\]

Since $\mathscr D_i(M)^{\vee} \rightarrow M^{\vee}$ is surjective by assumption, it suffices to construct a map $\widetilde{f}'$ lifting $f$, that is, a map
\[
\widetilde{\psi} : U(\mathfrak{u}) \otimes_\C \widetilde{B}_0 \longrightarrow \mathscr D_i(M)^{\vee}
\]
such that $\widetilde{\psi}\bigl( \mathscr D_i(I) \bigr) = 0$.

Now, the map $f$ lifts to a $U (\tb)$-module map
\[
\psi : U(\mathfrak{u}) \otimes_\C B_0 \longrightarrow M^{\vee}
\]
satisfying $\psi(I) = 0$. Let
\[
\widetilde{\psi}' : U(\mathfrak{u}) \otimes_\C \widetilde{B}_0 \longrightarrow M^{\vee}
\]
be the unique lift of $\psi$ as a $U (\tb)$-module morphism, induced by the universal property of the induction functor $U ( \tn ) \otimes_{\C [E_i]} \bullet$ from the category of $(\C E_i + \wth)$-modules to the category of $\tb$-modules applied to the $(\C E_i + \wth)$-module morphism lift $\widetilde{f}$ from $f$. Then $\widetilde{\psi}'(I) = 0$ holds as $\psi$ is written as the composition map
\[
U(\mathfrak{u}) \otimes_\C B_0 \hookrightarrow U(\mathfrak{u}) \otimes_\C \widetilde{B}_0 \stackrel{\widetilde{\psi}'}{\longrightarrow} M^{\vee}.
\]

By applying the correspondence given by Theorem~\ref{thm:adj}, we obtain a map
\[
\widetilde{\psi}'' : U(\mathfrak{u}) \otimes_\C \widetilde{B}_0 \longrightarrow \mathscr D_i(M)^{\vee}
\]
corresponding to $\widetilde{\psi}'$. Because $\widetilde{\psi}'$ vanishes on $I$, functoriality of $\mathscr D_i$ implies that $\widetilde{\psi}''$ vanishes on
\[
\mathscr D_i(I) \subset \mathscr D_i\bigl(U(\mathfrak{u}) \otimes_\C B_0\bigr)
= U(\mathfrak{u}) \otimes_\C \widetilde{B}_0.
\]
Hence
\[
\widetilde{\psi}''\bigl(\mathscr D_i(I)\bigr)=0.
\]

Therefore, we may take $\widetilde{\psi} := \widetilde{\psi}''$ as the desired lift. This defines the required map $\widetilde{f}'$, and hence the original map $\widetilde{f}$.

This completes the proof.
\end{proof}

\begin{cor}\label{cor:indsurj}
Let $k, l \in \Z_{\ge 0}$ with $k \ge l$, and let $M \in \fB_\bd$ be a module on which $K$ acts by $-k$. Suppose that for some $i \in \tI_\af$, we have
\[
M \subset \mathscr D_i(M), \qquad \text{and} \qquad \mathbb L^{<0} \mathscr D_i(M) = 0.
\]
Then, for each $\mu \in P$, there is a surjection:
\[
\shom_{\fB} \bigl( \mathscr D_i(M), ( D^{(l)}_{\mu} \otimes \C_{(k-l)\La_0} )^{*} \bigr) \longrightarrow \shom_{\fB} \bigl( M, ( D^{(l)}_{\mu} \otimes \C_{(k-l)\La_0} )^{*} \bigr).
\]
\end{cor}

\begin{proof}
We apply Theorem~\ref{thm:indsurj} with
\[
N = D^{(l)}_{\mu} \otimes \C_{(k-l)\La_0}.
\]
This module is $\tb$-cyclic by definition, so it remains to verify that
\[
N \subset \mathscr D_i(N)
\qquad\text{and}\qquad
\mathbb L^{<0}\mathscr D_i(N)=0.
\]

If $i \neq 0$, then both assertions follow directly from Theorem~\ref{thm:thinD}.

It remains to treat the case $i=0$. By Theorem~\ref{thm:thinD}, the module $D^{(l)}_{\mu}$ satisfies the hypotheses of Proposition~\ref{prop:string}. Hence $D^{(l)}_{\mu}$ admits a finite filtration whose associated graded pieces are of the two types described there. After twisting by $\C_{(k-l)\La_0}$, the finite-dimensional $\mathfrak{sl}(2,0)$-constituents remain finite-dimensional $\mathfrak{sl}(2,0)$-modules, and the one-dimensional constituents still satisfy the required non-negativity condition with respect to $\alpha_0^\vee$. Therefore the same conclusion applies to $D^{(l)}_{\mu} \otimes \C_{(k-l)\La_0}$. For each such graded piece, the required inclusion into its Demazure image and the vanishing of the lower derived Demazure functors follow by inspection. Applying d\'evissage to this filtration, we obtain
\[
D^{(l)}_{\mu} \otimes \C_{(k-l)\La_0}
\subset
\mathscr D_0\bigl(D^{(l)}_{\mu} \otimes \C_{(k-l)\La_0}\bigr)
\]
and
\[
\mathbb L^{<0}\mathscr D_0\bigl(D^{(l)}_{\mu} \otimes \C_{(k-l)\La_0}\bigr)=0.
\]

Therefore the hypotheses of Theorem~\ref{thm:indsurj} are satisfied, and the desired surjection follows.
\end{proof}

\section{Statement of the Main Theorems}\label{sec:statement}

We retain the notation and assumptions of the previous section. This section presents the main results of the paper and discusses some of their consequences. The proofs will be given in the subsequent sections.

\subsection{\texorpdfstring{$\ext$}{Ext}-orthogonality theorem}

\begin{thm}[Level $k$-duality]\label{thm:dual}
Let $k \in \Z_{>0}$. Then the following statements hold:
\begin{enumerate}
\item For each $\la, \mu \in P$, we have
\begin{equation}
\ext^{i}_\fB ( \bD _\la^{(k)}, ( D_\mu ^{(k)} )^{*}) \cong 
\begin{cases} 
\C & \text{if } i = 0 \text{ and } \la + \mu = 0, \\
0 & \text{otherwise}.
\end{cases}
\end{equation}

\item For each $\la, \mu \in P^+$, we have
\begin{equation}
\ext^{i}_\fC ( \bW _\la^{(k)}, ( W_\mu ^{(k)} )^{*}) \cong 
\begin{cases} 
\C & \text{if } i = 0 \text{ and } \la + \mu_- = 0, \\
0 & \text{otherwise}.
\end{cases}
\end{equation}
\end{enumerate}	
\end{thm}

The case $k = 1$ of Theorem~\ref{thm:dual}(1) corresponds to~\cite[Theorem~4.19]{CK18} (Theorem~\ref{thm:Dorth}). The case $k = 1$ of Theorem~\ref{thm:dual}(2) will be deduced from Theorem~\ref{thm:dual}(1) in the proof of Theorem~\ref{thm:dual}(2).

\subsection{Filtration criteria}\label{subsec:filt}

\begin{defn}
Let $k \in \Z_{> 0}$.

A module $M \in \fB_\bd$ is said to admit a $D^{(k)}$-filtration (resp. a $\bD^{(k)}$-filtration) if there exists a decreasing separable filtration
\[
M = M^0 \supset M^1 \supset M^2 \supset \cdots, \qquad \bigcap_{i \ge 0} M^i = 0
\]
by $\tb$-submodules, such that for each $i \in \Z_{\ge 0}$, the successive quotient $M^i / M^{i+1}$ is isomorphic to $\bq^{m_i} D^{(k)}_{\mu_i}$ (resp. $\bq^{m_i} \bD^{(k)}_{\mu_i}$), for some $m_i \in \Z$ and $\mu_i \in P$, and possibly after twisting by $\C_{l\Lambda_0}$ for some $l \in \Z$.

Similarly, a module $M \in \fC_\bd$ is said to admit a $W^{(k)}$-filtration (resp. a $\bW^{(k)}$-filtration) if there exists a decreasing separable filtration
\[
M = M^0 \supset M^1 \supset M^2 \supset \cdots, \qquad \bigcap_{i \ge 0} M^i = 0
\]
by $\tg_{\ge 0}$-submodules, such that for each $i \in \Z_{\ge 0}$, the successive quotient $M^i / M^{i+1}$ is isomorphic to $\bq^{m_i} W^{(k)}_{\mu_i}$ (resp. $\bq^{m_i} \bW^{(k)}_{\mu_i}$), for some $m_i \in \Z$ and $\mu_i \in P^+$, and possibly after twisting by $\C_{l\La_0}$ for some $l \in \Z$.
\end{defn}

Suppose that $M$ admits an $X$-filtration for some $X \in \{ \bD^{(k)}, D^{(k)}, \bW^{(k)}, W^{(k)} \}$. Then we define the collection of elements $(M : X_\la)_q \in \Z(\!(q)\!)$ by the expansion
\begin{equation}
\gch M = \sum_{\la} (M : X_\la)_q \cdot \gch X_\la,\label{eqn:()q}
\end{equation}
where the sum is taken over $\la \in P$ or $\la \in P^+$, depending on whether $X = \bD^{(k)}, D^{(k)}$ or $X = \bW^{(k)}, W^{(k)}$, respectively. By Lemma~\ref{lem:top}, these coefficients are uniquely determined. These coefficients record the graded multiplicities of the modules $X_\la$ in $M$ with respect to the $X$-filtration.

Similarly, for each $M \in \fC_\bd$, we write $[M : V_\la]_q \in \Z_{\ge 0}(\!(q)\!)$ for the graded multiplicity of $V_\la$ in $M$. Namely, we have
\[
\gch M = \sum_{\la} [M : V_\la]_q \cdot \ch V_\la.
\]
The subscript $q$ may be replaced by another variable (e.g., $q^{-1}$) to indicate a change of grading.  
In particular, we have
\begin{equation}
(M : X_\la)_{q^{-1}} = \Bigl. (M : X_\la)_q \Bigr|_{q \mapsto q^{-1}}, \quad
[M : V_\la]_{q^{-1}} = \Bigl. [M : V_\la]_q \Bigr|_{q \mapsto q^{-1}}.\label{eqn:qtoqinv}
\end{equation}

\begin{thm}[Level $k$ criterion for filtrations]\label{thm:critfilt}
Let $k \in \Z_{>0}$. Then the following statements hold:
\begin{enumerate}
\item A module $M \in \fB_\bd$ admits a $D^{(k)}$-filtration if and only if
 \begin{equation}
\ext^{1}_\fB ( \bD _\la^{(k)}, M^{\vee} ) = 0 \qquad \text{for all } \la \in P.
\end{equation}

\item A module $M \in \fB_\bd$ admits a $\bD^{(k)}$-filtration if and only if
 \begin{equation}
\ext^{1}_\fB ( M, ( D_\la ^{(k)} )^{*} ) = 0 \qquad \text{for all } \la \in P.
\end{equation}

\item A module $M \in \fC_\bd$ admits a $W^{(k)}$-filtration if and only if
\begin{equation}
\ext^{1}_\fC ( \bW _\la^{(k)}, M^{\vee} ) = 0 \qquad \text{for all } \la \in P^+.
\end{equation}

\item A module $M \in \fC_\bd$ admits a $\bW^{(k)}$-filtration if and only if
\begin{equation}
\ext^{1}_\fC ( M, ( W_\la ^{(k)} )^{*} ) = 0 \qquad \text{for all } \la \in P^+.
\end{equation}
\end{enumerate}
Moreover, in each of the above cases, all higher $\ext^i$ vanish for $i > 1$.
\end{thm}

The case $k = 1$ of Theorem~\ref{thm:critfilt}(2)(4) corresponds to~\cite[Theorem~5.9(ii)(i)]{CK18}, respectively.  
The proofs of Theorem~\ref{thm:critfilt}(3)(4) are analogous to those of Theorem~\ref{thm:critfilt}(1)(2), with the use of Theorem~\ref{thm:dual}(1) replaced by that of Theorem~\ref{thm:dual}(2).

Throughout the proofs of Theorems~\ref{thm:dual} and~\ref{thm:critfilt}, we assume their validity for all strictly smaller values of $k$, with the base case given by Theorem~\ref{thm:dual}(1) at $k = 1$ (Theorem~\ref{thm:Dorth}).

In fact, the implication from Theorem~\ref{thm:dual} at level $k$ to Theorem~\ref{thm:critfilt} at level $k$ is largely formal, once Proposition~\ref{prop:Dhw} is established at level $k$.  
Accordingly, in addition to the preparations in Sections~\ref{sec:prelim}--\ref{sec:lift}, we devote Sections~\ref{sec:hw}--\ref{sec:filt} to preparing several auxiliary results that are required to prove Theorem~\ref{thm:dual} at level $k$ from Theorems~\ref{thm:dual} and~\ref{thm:critfilt} at levels $< k$.

We also note that the induction closes within Theorems~\ref{thm:dual}(1) and~\ref{thm:critfilt}(2); Theorems~\ref{thm:dual}(2) and~\ref{thm:critfilt}(1)(3)(4) may then be regarded as their corollaries.

\begin{cor}\label{cor:direct}
Let $k \in \Z_{>0}$, and assume that Theorem~\ref{thm:critfilt} holds at level $k$. Then:
\begin{itemize}
\item If $M, N \in \fC_\bd$ and $M \oplus N$ admits a $\bW^{(k)}$-filtration, then $M$ also admits a $\bW^{(k)}$-filtration;
\item If $M, N \in \fB_\bd$ and $M \oplus N$ admits a $\bD^{(k)}$-filtration, then $M$ also admits a $\bD^{(k)}$-filtration.
\end{itemize}
In particular, the class of $\bW^{(k)}$-filtered (resp. $\bD^{(k)}$-filtered) modules in $\fC_\bd$ (resp. $\fB_\bd$) is closed under taking direct summands.
\end{cor}

\begin{proof}
This follows immediately from Theorem~\ref{thm:critfilt}, as the functors $\ext^1_\fC$ and $\ext^1_\fB$ commute with finite direct sums.
\end{proof}

\begin{lem}\label{lem:surjfilt}
Let $k \in \Z_{>0}$, and assume that Theorem~\ref{thm:critfilt} holds at level $k$. Let $M, N \in \fB_{\bd}$ be two $\bD^{(k)}$-filtered modules, and suppose that there is a surjective morphism $f : M \twoheadrightarrow N$. Then $\ker f$ also admits a $\bD^{(k)}$-filtration. The analogous statement holds for the filtrations of type $D^{(k)}$, $\bW^{(k)}$, or $W^{(k)}$.
\end{lem}

\begin{proof}
Since the proofs are similar in all cases, we treat only the case of $\bD^{(k)}$-filtrations. Consider the short exact sequence
\[
0 \rightarrow \ker f \rightarrow M \rightarrow N \rightarrow 0.
\]
Applying the functor $\ext^\bullet_\fB(\bullet, ( D_\la^{(k)} )^{*})$ yields a long exact sequence, from which we extract an exact sequence for each $i \ge 0$:
\[
\ext_{\fB}^i ( M, ( D_\la^{(k)} )^{*} ) 
\longrightarrow\ext_{\fB}^i ( \ker f, ( D_\la^{(k)} )^{*} ) 
\longrightarrow\ext_{\fB}^{i+1} ( N, ( D_\la^{(k)} )^{*} ).
\]
Since $M$ and $N$ both admit $\bD^{(k)}$-filtrations, Theorem~\ref{thm:critfilt} implies
\[
\ext_{\fB}^i ( M, ( D_\la^{(k)} )^{*} ) = 0 =\ext_{\fB}^{i+1} ( N, ( D_\la^{(k)} )^{*} ) \qquad \text{for all } i > 0.
\]
Therefore,
\[
\ext_{\fB}^{>0} ( \ker f, ( D_\la^{(k)} )^{*} ) = 0,
\]
and Theorem~\ref{thm:critfilt} implies that $\ker f$ admits a $\bD^{(k)}$-filtration.
\end{proof}

\begin{cor}\label{cor:multcrit}
Let $k \in \Z_{>0}$, and assume Theorem~\ref{thm:dual} and Theorem~\ref{thm:critfilt} hold at level~$k$. Then:
\begin{enumerate}
\item If $M \in \fB_{\bd}$ admits a $\bD^{(k)}$-filtration, then for each $\la \in P$, we have
\[
(M : \bD^{(k)}_{\la})_{q^{-1}} = \gdim \hom_{\fB} \big( M, ( D^{(k)}_{-\la} )^{*} \big).
\]
Similarly, if $M \in \fB_{\bd}$ admits a $D^{(k)}$-filtration, then
\[
(M : D^{(k)}_{\la})_{q^{-1}} = \gdim \hom_{\fB} \big( \bD^{(k)}_{-\la}, M^{\vee} \big).
\]

\item If $M \in \fC_{\bd}$ admits a $\bW^{(k)}$-filtration, then for each $\la \in P^+$, we have
\[
(M : \bW^{(k)}_{\la})_{q^{-1}} = \gdim \hom_{\fC} \big( M, ( W^{(k)}_{-\la_-} )^{*} \big).
\]
Similarly, if $M \in \fC_{\bd}$ admits a $W^{(k)}$-filtration, then
\[
(M : W^{(k)}_{\la})_{q^{-1}} = \gdim \hom_{\fC} \big( \bW^{(k)}_{-\la_-}, M^{\vee} \big).
\]
\end{enumerate}
\end{cor}

\begin{proof}
Since the proofs in all cases are similar (and simpler for the cases involving $W^{(k)}$ and $D^{(k)}$), we focus on the case of a $\bD^{(k)}$-filtration.

Assume $M \in \fB_\bd$. Then there exists a quotient map
\[
f : M \twoheadrightarrow \bq^m \bD^{(k)}_{\mu}
\]
for some $m \in \Z$ and $\mu \in P$. By Lemma~\ref{lem:surjfilt}, we obtain the following isomorphism for each $\la \in P$:
\begin{equation}
\hom_\fB \big( \ker f, ( D_{\la}^{(k)} )^{*} \big) \oplus \bq^{-m} \C^{\delta_{\la+\mu,0}} 
\cong
\hom_\fB \big( M, ( D_{\la}^{(k)} )^{*} \big). \label{eqn:onefilt}
\end{equation}
In addition, we have
\[
\ext^{>0}_\fB \big( \ker f, ( D_{\la}^{(k)} )^{*} \big) = 0 \qquad \text{for all } \la \in P.
\]

We examine a decreasing separable $\bD^{(k)}$-filtration
\[
M=M^0 \supset M^1 \supset M^2 \supset \cdots \qquad \bigcap_{i \ge 0} M^i = 0
\]
provided by Theorem~\ref{thm:critfilt}(2). Since $M \in \fB_\bd$, for each $m \in \Z$, there exists $i_m \gg 0$ such that $M^i$ is concentrated in degrees $> m$ for all $i > i_m$. This implies that
\[
M = \varprojlim_i M / M^i,
\]
and the claim follows by repeatedly applying~\eqref{eqn:onefilt}, as required.
\end{proof}

\subsection{Highest-weight structure}\label{subsec:hw}

Let $(\bullet,\bullet)_0$ be the $W$-invariant bilinear form on~$\h$ normalized by $(\vartheta^{\vee},\vartheta^{\vee})_0 = 2$. We extend it to a $W_\af$-invariant bilinear form $(\bullet,\bullet)$ on~$\wth$ by setting
\[
 (\la + k \La_0 + m \delta,\, \mu + k' \La_0 + m' \delta ) := (\la,\mu)_0 + km' + k' m \qquad (\la,\mu \in P;\; k,k',m,m' \in \Z).
\]

We define a $W_\af$-invariant quadratic form $\mathbf{q}$ on $P_\af$ by
\[
\mathbf{q}(\la + k \La_0 + m \delta) := \big( \la + k \La_0 + m \delta,\; \la + k \La_0 + m \delta \big) \qquad (\la \in P;\; m \in \Z).
\]
Since we have
\[
\mathbf{q}(\la + k \La_0 + m \delta) = (\la, \la)_0 + 2km,
\]
we obtain the inequality
\[
\mathbf{q}(\La - \delta) < \mathbf{q}(\La) \qquad \text{for all } \La \in \bigsqcup_{k > 0} P_{\af,k}.
\]

\begin{prop}[{\cite[Proposition~12.5(d)]{Kac}}]\label{prop:bound}
Let $k \in \Z_{>0}$. For each $\La \in P_{\af,k}^+$, every $\wth$-weight $\La'$ occurring in $L(\La)$ satisfies $\mathbf{q}(\La') \le \mathbf{q}(\La)$.
\end{prop}

\begin{defn}
We define a partial order $\blacktriangleleft$ on $P_\af$ by declaring:
\begin{align*}
\La' \blacktriangleleft \La \quad \Longleftrightarrow & \quad \mathbf{q}(\La') < \mathbf{q}(\La), \quad \text{and} \quad \La' (K) = \La (K) \quad \text{or} \\
& \quad \mathbf{q}(\La') = \mathbf{q}(\La),\quad W_\af \La' = W_\af \La, \quad \text{and} \quad \La' \in \La + \sum_{i \in \tI_\af} \Z_{\ge 0} \al_i.
\end{align*}
\end{defn}

\begin{rem}
Two elements $\La, \La' \in P_\af$ can be comparable with respect to $\blacktriangleleft$ only if $\La(K) = \La'(K)$. Note that $\mathbf{q}(\La') = \mathbf{q}(\La)$ and $W_\af \La' \neq W_\af \La$ implies that $\La'$ is not a $\wth$-weight of $L(\La)$.
\end{rem}

\begin{prop}\label{prop:Dhw}
For each $k \in \Z_{>0}$, $\la \in P$, and $m \in \Z$, we have:
\begin{enumerate}
\item The module $\bq^m D_\la^{(k)}$ is the largest $\tb$-module generated by a $\wth$-weight of the form $\la + k \La_0 + m \delta$, whose $\wth$-weights $\La$ satisfy
\begin{equation}
\La \blacktriangleleft \la + k \La_0 + m \delta.\label{eqn:exbound}
\end{equation}
In addition, a $\wth$-weight $\La$ of $\bq^m \hd \bigl( \ker ( Q_{\la + k \La_0} \to D_\la^{(k)} ) \bigr)$ satisfies
\begin{equation}
\la + k \La_0 + m \delta \blacktriangleleft \La.\label{eqn:defbound}
\end{equation}

\item The module $\bq^m ( \bD_\la^{(k)} )^{\vee}$ is the largest $\tb$-module cogenerated by a $\wth$-weight of the form $- \la + k \La_0 + m \delta$, whose $\wth$-weights $\La$ satisfy
\begin{equation}
\La \blacktriangleleft - \la + k \La_0 + m \delta.\label{eqn:exbound2}
\end{equation}
In addition, a $\wth$-weight $\La$ of $\bq^m \bigl( \hd \ker ( Q_{\la - k \La_0} \to \bD_\la^{(k)}  ) \bigr)^{\vee}$ satisfies
\begin{equation}
- \la + k \La_0 + m \delta \blacktriangleleft \La.\label{eqn:defbound2}
\end{equation}
\end{enumerate}
\end{prop}

\begin{proof}
The proof of Proposition~\ref{prop:Dhw} will be given later as part of the inductive proof of Theorem~\ref{thm:dual}(1). In particular, once Theorem~\ref{thm:dual}(1) is known at level $k$, Proposition~\ref{prop:Dhw} follows at level $k$. Thus, the assertion is established when Theorem~\ref{thm:dual}(1) is established for all levels.
\end{proof}

\begin{cor}\label{cor:Whw}
For each $k \in \Z_{>0}$, $\la \in P^+$, and $m \in \Z$, we have:
\begin{enumerate}
\item The module $\bq^m W_\la^{(k)}$ is the largest $\tg_{\ge 0}$-module generated by a $\wth$-weight of the form $\la_- + k \La_0 + m \delta$, whose $\wth$-weights $\La$ satisfy
\[
\La \blacktriangleleft \la_- + k \La_0 + m \delta.
\]

\item The module $\bq^m ( \bW_\la^{(k)} )^{\vee}$ is the largest $\tg_{\ge 0}$-module cogenerated by a $\wth$-weight of the form $- \la + k \La_0 + m \delta$, whose $\wth$-weights $\La$ satisfy
\[
\La \blacktriangleleft - \la + k \La_0 + m \delta.
\]
\end{enumerate}
\end{cor}

\begin{proof}
Since any $\tg_{\ge 0}$-module extension restricts to a $\tb$-module extension, we may restrict the module structures and reduce the assertions to Proposition~\ref{prop:Dhw}.

The first assertion follows from Proposition~\ref{prop:Dhw}(1), using the identification $W_\la^{(k)} = D_{\la_-}^{(k)}$ and the $W$-invariance of $\mathbf{q}$. The second assertion follows from Proposition~\ref{prop:WDfilt} and Proposition~\ref{prop:Dhw}(2).
\end{proof}

\begin{rem}\label{rem:hw}
Proposition~\ref{prop:Dhw} shows that the collection $\{ \bq^m D_\la^{(k)} \}_{\la, m}$ forms the set of standard modules, and $\{ \bq^m ( \bD_\la^{(k)} )^{\vee} \}_{\la, m}$ forms the set of costandard modules, with respect to the partial order $\blacktriangleleft$, in the sense of highest-weight theory~\cite{CPS88} (without the finiteness assumption on the poset). An axiomatic treatment of such a setting is given in~\cite{BS24}.

The same description applies to the pair $\{ \bq^m W_\la^{(k)} \}_{\la, m}$ and $\{ \bq^m ( \bW_\la^{(k)} )^{\vee} \}_{\la, m}$ by Corollary~\ref{cor:Whw}. Thus, Theorem~\ref{thm:dual}(1) establishes that the collection of modules under consideration satisfies the axioms required for a highest-weight category with respect to $\blacktriangleleft$.
\end{rem}

\section{Highest-weight structure: proofs}\label{sec:hw}

We retain the setting of the previous section.

\begin{prop}\label{prop:Dhw-ind}
For each $k \in \Z_{>0}$, $\la \in P$, and $m \in \Z$, Proposition~\ref{prop:Dhw}(1) holds unconditionally, and Proposition~\ref{prop:Dhw}(2) holds if we assume Theorem~\ref{thm:dual}(1) at level $k$.
\end{prop}

The first part follows directly from the defining relations of thin Demazure modules, while the second is deduced from level-$k$ duality by ruling out nontrivial one-dimensional extensions.

\begin{proof}[Proof of Proposition~\ref{prop:Dhw-ind}]
We first prove Proposition~\ref{prop:Dhw}(1). By definition, we have
\[
\bq^m D_{\la}^{(k)} \subset \bq^m L(w(\la + k \La_0)),
\]
where $w \in W_\af$ is chosen so that $w(\la + k \La_0) \in P_{\af}^+$. By Proposition~\ref{prop:bound}, all $\wth$-weights appearing in $\bq^m D_\la^{(k)}$ satisfy~\eqref{eqn:exbound}. Moreover, by definition, $\bq^m D_{\la}^{(k)}$ is generated by a vector of $\wth$-weight $\la + k \La_0 + m \delta$.

It remains to show the maximality of $\bq^m D_{\la}^{(k)}$. From the description in~\eqref{eqn:defD}, the $\wth$-weights of the defining relations of $D_\la^{(k)}$ as a $U(\tb)$-module are given by
\begin{equation}
\la + k \La_0 + m_\al \al, \qquad m_\al = \begin{cases} -
\langle \al^{\vee}, \la + k \La_0 \rangle + 1 & \text{if } \langle \al^{\vee}, \la + k \La_0 \rangle \le 0\\
1 & \text{otherwise}
\end{cases}
\label{eqn:exdefD}
\end{equation}
for each real root $\al \in \Delta_\af^+$. This implies that the defining relation of weight $\la + k \La_0 + m_\al \al$ satisfies~\eqref{eqn:defbound} when $s_\al ( \la + k \La_0) \ge \la + k \La_0$ with respect to the dominance order.

By convexity of the level sets of $\mathbf{q}$ and Proposition~\ref{prop:bound}, we deduce that the weights appearing in~\eqref{eqn:exdefD} satisfy~\eqref{eqn:defbound}. This confirms the maximality of $\bq^m D_{\la}^{(k)}$, proving the first assertion.

\medskip
We now prove Proposition~\ref{prop:Dhw}(2), assuming Theorem~\ref{thm:dual}(1). We begin by examining the dual module $\bq^{-m} \bD_\la^{(k)}$. By definition, $\bq^{-m} \bD_\la^{(k)}$ is generated by a vector of $\wth$-weight $\la - k \La_0 - m \delta$. Moreover, it appears as a subquotient of
\[
{}^\theta L(w(-\la + k \La_0 + m \delta)),
\]
where $w \in W_\af$ is chosen such that $w(-\la + k \La_0) \in P_\af^+$. Again by Proposition~\ref{prop:bound}, any $\wth$-weight space of ${}^\theta L(w(-\la + k \La_0 + m \delta))$ in the $W_\af$-orbit of $-\la + k \La_0 + m \delta$ (except for $-\la + k \La_0 + m \delta$ itself, noting that $\theta$ negates $\wth$-weights) is quotiented out in $\bq^{-m} \bD_\la^{(k)}$. Therefore, all $\wth$-weights of $\bq^m ( \bD_\la^{(k)} )^{\vee}$ satisfy~\eqref{eqn:exbound2}. This proves the weight estimate part of the second assertion.

Let $\gamma + s \delta + k \La_0 \in P_\af$ be a weight such that we have a non-split short exact sequence
\begin{equation}
0 \rightarrow (\bD_\la^{(k)} )^{\vee} \rightarrow E \rightarrow \C_{\gamma + k \La_0 + s \delta} \rightarrow 0\label{eqn:bDext-one}
\end{equation}
in the category $\fB$. If no such weight satisfies
\[
-\la + k\La_0 \not\blacktriangleleft \gamma + k\La_0 + s\delta,
\]
then there is nothing to prove. Therefore, we assume
\begin{equation}
-\la + k\La_0 \not\blacktriangleleft \gamma + k\La_0 + s\delta\label{eqn:wt-est-bD}
\end{equation}
to derive a contradiction. By the first part of the assertion, we have
\[
\dim \Hom_{\wth} (\C_{\nu + k \La_0 + r \delta}, \bq^s D_\gamma^{(k)}) = \begin{cases} 1 & (\nu + k \La_0 + r \delta = \gamma + s \delta + k \La_0)\\0 & (\nu + k \La_0 + r \delta \not\blacktriangleleft \gamma + s \delta + k \La_0)\end{cases}.
\]
It follows that every one-dimensional subquotient $L$ of $\bq^s D_\gamma^{(k)}$ has $\wth$-weight $\not\blacktriangleright - \la + k \La_0$ by~\eqref{eqn:wt-est-bD}. Thus, we have
\[
\Hom_{\fB} ( \bD_\la^{(k)}, ( \bq^s \ker ( D_\gamma^{(k)} \twoheadrightarrow \C_{\gamma + k \La_0}) )^{\vee}) = 0.
\] 
It follows that the long exact sequence associated to the short exact sequence
\[
0 \rightarrow \C_{\gamma + k \La_0 + s \delta}^{\vee} \rightarrow \bq^{-s} ( D_\gamma^{(k)} )^{\vee} \rightarrow \bq^{-s} \ker ( D_\gamma^{(k)} \twoheadrightarrow \C_{\gamma + k \La_0})^{\vee} \rightarrow 0
\]
yields an exact sequence
\[
0 = \Hom_{\fB} ( \bD_\la^{(k)}, \bq^{-s} \ker ( D_\gamma^{(k)} \twoheadrightarrow \C_{\gamma + k \La_0})^{\vee}) \rightarrow \Ext^1_{\fB} ( \bD_\la^{(k)}, \C_{\gamma + k \La_0 + s \delta}^{\vee}) \rightarrow \Ext^1_{\fB} ( \bD_\la^{(k)}, \bq^{-s} ( D_\gamma^{(k)} )^{\vee}).
\]
In particular, we have
\[
0 \neq \Ext^1_{\fB} ( \bD_\la^{(k)}, \bq^{-s} ( D_\gamma^{(k)} )^{\vee}) \subset \ext^1_{\fB} ( \bD_\la^{(k)}, \bq^{-s} ( D_\gamma^{(k)} )^{\vee}).
\] 
This contradicts Theorem~\ref{thm:dual}(1). Therefore, we cannot have a non-split short exact sequence~\eqref{eqn:bDext-one} with~\eqref{eqn:wt-est-bD}, confirming~\eqref{eqn:defbound2}.

Hence, any $\wth$-weight occurring in a potential one-dimensional extension as a $\tb$-module violates~\eqref{eqn:exbound2}. This establishes the maximality of $\bD_{\la}^{(k)}$ and completes the proof of the maximality part of the second assertion.

This completes the proof.
\end{proof}

\begin{cor}\label{cor:Dhwext}
Fix $k \in \Z_{> 0}$, $\la, \mu \in P$, and $m,n \in \Z$. We have
\begin{enumerate}
\item If we have $\la + k \La_0 + m \delta \not\blacktriangleleft \mu + k \La_0 + n \delta$, then we have
\[
\Ext_{\fB}^1 ( \bq^m D_\la^{(k)}, \bq^{n} D_\mu^{(k)} ) = 0;
\]
\item Assume that Theorem~\ref{thm:dual}(1) holds for $k$. If we have
\[
- \mu + k \La_0 - n \delta \not\blacktriangleleft - \la + k \La_0 - m \delta,
\]
then we have
\[
\Ext_{\fB}^1 ( \bq^m \bD_\la^{(k)}, \bq^{n} \bD_\mu^{(k)} ) = 0.
\]
\end{enumerate}
\end{cor}

\begin{proof}
We prove the first assertion. By the first part of Proposition~\ref{prop:Dhw-ind}, we find that the $\wth$-weights $\La$ of $\bq^{n} D_\mu^{(k)}$ satisfy $\la + k \La_0 + m \delta \not\blacktriangleleft \La$. Since extensions by one-dimensional modules of such weights are prohibited by~\eqref{eqn:defbound}, d\'evissage applied to a $\tb$-module filtration of $\bq^{n} D_\mu^{(k)}$ yields
\[
\Ext_{\fB}^1 ( \bq^m D_\la^{(k)}, \bq^{n} D_\mu^{(k)} ) = 0.
\]
This is the first assertion.

We prove the second assertion. We have
\[
\Ext_{\fB}^1 ( \bq^m \bD_\la^{(k)}, \bq^{n} \bD_\mu^{(k)} ) \cong \Ext_{\fB}^1 ( (\bq^{n} \bD_\mu^{(k)} )^{\vee}, ( \bq^m \bD_\la^{(k)})^{\vee} )
\]
by the Yoneda interpretation of $\Ext^1$ (Proposition~\ref{prop:Yoneda}). By the second part of Proposition~\ref{prop:Dhw-ind}, we find that the $\wth$-weight $\La$ of $( \bq^m \bD_\la^{(k)} )^{\vee}$ satisfies $- \mu + k \La_0 - n \delta \not\blacktriangleleft \La$. Since extensions by one-dimensional modules of such weights are prohibited by~\eqref{eqn:defbound2}, we have a surjective projective system $\{X(r)\}_{r \ge 0}$ such that
\[
\varprojlim X( r ) \cong \bq^m \bD_\la^{(k)}
\]
and each $X(r)$ ($r \ge 0$) satisfies
\[
\Ext_{\fB}^1 ( X(r), \bq^{n} \bD_\mu^{(k)} ) = 0
\]
by the d\'evissage applied to a $\tb$-module filtration of $X(r)$. By Lemma~\ref{lem:comm}, we deduce
\[
\Ext_{\fB}^1 ( \bq^m \bD_\la^{(k)}, \bq^{n} \bD_\mu^{(k)} ) = \varinjlim_r \Ext_{\fB}^1 ( X(r), \bq^{n} \bD_\mu^{(k)} )  = 0.
\]
This proves the second assertion and completes the proof.
\end{proof}

\begin{cor}\label{cor:Whwext}
Fix $k \in \Z_{> 0}$, $\la, \mu \in P^+$, and $m,n \in \Z$. We have
\begin{enumerate}
\item If we have $\la + k \La_0 + m \delta \not\blacktriangleleft \mu + k \La_0 + n \delta$, then we have
\[
\Ext_{\fB}^1 ( \bq^m W_\la^{(k)}, \bq^{n} W_\mu^{(k)} ) = 0;
\]
\item Assume that Theorem~\ref{thm:dual}(1) holds for $k$. If we have
\[
- \mu_- + k \La_0 - n \delta \not\blacktriangleleft - \la + k \La_0 - m \delta,
\]
then we have
\[
\Ext_{\fB}^1 ( \bq^m \bW_\la^{(k)}, \bq^{n} \bW_\mu^{(k)} ) = 0.
\]
\end{enumerate}
\end{cor}

\begin{proof}
The first assertion follows from Corollary~\ref{cor:Dhwext}(1) by $W_\la^{(k)} = D_{\la_-}^{(k)}$ and $W_\mu^{(k)} = D_{\mu_-}^{(k)}$, together with the $W$-invariance of the form $\mathbf{q}$.

We prove the second assertion. 
Since $-\mu_- + k\La_0 - n\delta$ is the smallest element in its $W$-orbit with respect to $\blacktriangleleft$, we have $- \mu_- + k \La_0 - n \delta \blacktriangleleft - u \mu + k \La_0 - n \delta$ for each $u \in W$ with $- \mu_- \neq - u \mu$. Combined with $- \mu_- + k \La_0 - n \delta \not\blacktriangleleft - \la + k \La_0 - m \delta$, we conclude
\[
-u \mu + k \La_0 - n \delta \not\blacktriangleleft - \la + k \La_0 - m \delta \qquad u \in W
\]
by transitivity. By using the fact that $-\la + k\La_0 - m\delta$ is the largest element in its $W$-orbit with respect to $\blacktriangleleft$, we further deduce
\[
-u \mu + k \La_0 - n \delta \not\blacktriangleleft - v \la + k \La_0 - m \delta \qquad u, v \in W.
\]

Hence Corollary~\ref{cor:Dhwext}(2) and the d\'evissage applied to Proposition~\ref{prop:WDfilt} yield
\[
\Ext_{\fB}^1 ( \bq^m \bW_\la^{(k)}, \bq^{n} \bW_\mu^{(k)} ) = 0.
\]
This proves the second assertion and completes the proof.
\end{proof}

\section{$\ext^1$-vanishing results}\label{sec:ext1}

We retain the setup from the previous section.

\begin{prop}\label{prop:Lext}
Let $k \in \Z_{> 0}$. For each $\La \in P_{\af, k}^+$ and $\mu \in P^+$, we have
\begin{equation}
\ext^1_\fB ( {}^{\theta} L ( \La ), ( D^{(k)}_{\mu} )^{*} ) = 0.\label{eqn:pc}
\end{equation}
\end{prop}

\begin{proof}
We first find $\La' \in P_{\af, k}^+$ and $w \in W_\af$ such that
\[
\mu + k \La_0 = w \La' \quad \text{and hence} \quad D^{(k)}_{\mu} = \mathscr D_{w} ( \C _{\La'} ).
\]

We rewrite~\eqref{eqn:pc} as
\begin{align*}
\ext^1_\fB \bigl( {}^{\theta} L ( \La ), ( D^{(k)}_{\mu} )^{*} \bigr)
&=\ext^1_\fB \bigl( {}^{\theta} L ( \La ), \bigl( \mathscr D_w ( \C_{\La'} ) \bigr)^{*} \bigr) \\
&\cong\ext^1_\fB \bigl( \mathscr D_{w^{-1}} ( {}^{\theta} L ( \La ) ), ( \C _{\La'} )^{*} \bigr) \\
&=\ext^1_\fB \bigl( {}^{\theta} L ( \La ), ( \C_{\La'} )^* \bigr),
\end{align*}
where the isomorphism follows from Theorems~\ref{thm:adj} and~\ref{thm:tD}.

By the BGG-resolution~\eqref{eqn:BGG-res}, the last term vanishes:
\[
\ext^1_\fB \bigl( {}^{\theta} L ( \La ), \C_{-\La'} \bigr) = 0,
\]
since $\La'$ does not lie in the set $\{\, s_i \cdot  \La + m \delta \mid i \in \tI_\af,\ m \in \Z \,\}$. This completes the proof.
\end{proof}

\begin{lem}\label{lem:WWdual}
Assume Theorem~\ref{thm:dual}(2) for $k = 1$.
For $\la, \mu \in P$ with $\la_+ + \mu_- \neq 0$, we have
\[
\ext_{\fB}^{\bullet} ( \bW_\la, ( \bW_\mu )^{\vee} ) = 0.
\]
In the case $\la_+ + \mu_- = 0$, we have
\[
\ext_{\fB}^{\ge 1} ( \bW_\la, ( \bW_\mu )^{\vee} ) = 0.
\]
\end{lem}

\begin{proof}
Suppose $\la = u \la_-$. Then the length formula~\cite[(2.4.1)]{Mac03} asserts that for any $\beta \in Q \cap P^+$ with $\langle \alpha_i^{\vee}, \beta \rangle > 0$ for all $i \in \tI$, we have
\[
\ell ( u t_{-\beta} ) = - \ell ( u ) + \ell ( t_{-\beta} ).
\]
It follows that
\[
\bW_\la \cong \bq^{\langle \beta, \la_- \rangle} \mathscr{D}_{u t_{-\beta}} ( \bW_{\la_-} ),
\]
which implies
\begin{align*}
\ext^\bullet_{\fB} ( \bW_\la, ( \bW_\mu )^{\vee} ) 
& \cong \ext^\bullet_{\fB} ( \bq^{\langle \beta,\la_- \rangle} \mathscr{D}_{u t_{-\beta}} ( \bW_{\la_-} ), ( \bW_\mu )^{\vee} ) \\
& \cong \bq^{-\langle \beta, \la_- \rangle} \ext^\bullet_{\fB} ( \bW_{\la_-}, ( \mathscr{D}_{t_{\beta} u^{-1}} ( \bW_\mu ) )^{\vee} ).
\end{align*}
Since $\mathscr{D}_{w_0} ( \bW_{\la_-} ) \cong \bW_{\la_-}$, we further obtain
\[
\ext^\bullet_{\fB} ( \bW_{\la_-}, ( \mathscr{D}_{t_{\beta} u^{-1}} ( \bW_\mu ) )^{\vee} )
\cong \ext^\bullet_{\fB} ( \bW_{\la_-}, ( ( \mathscr{D}_{w_0} \circ \mathscr{D}_{t_{\beta} u^{-1}} ) ( \bW_\mu ) )^{\vee} ).
\]
Here, $( \mathscr{D}_{w_0} \circ \mathscr{D}_{t_{\beta} u^{-1}} ) ( \bW_\mu )$ is $\g$-stable, and is therefore isomorphic to a grading shift of $\bW_{\mu_-}$.

Let $R := \mathrm{End}_{\fC} ( \bW_{\mu_+}^{(1)} )$, which is a polynomial ring by Theorem~\ref{thm:CKmain1}. For each $m$, denote by $R_{\le m}$ the degree-$m$ truncation. Then
\begin{equation}\label{eqn:bWinv}
\bW_{\mu_+}^{(1)} \cong \varprojlim_m \bigl( \bW_{\mu_+}^{(1)} \otimes_{R} R_{\le m} \bigr).
\end{equation}

Each term in the inverse system~\eqref{eqn:bWinv} is a finite self-extension of
\[
W_{\mu_+}^{(1)} \otimes \C_{-2\La_0}.
\]

By Proposition~\ref{prop:DCext}, we deduce
\begin{align*}
\ext^\bullet_{\fB} \left( \bW_{\la_-}, \bW_{\mu_-}^{\vee} \right)
& \cong \varinjlim_m \ext^\bullet_{\fB} \left( \bW_{\la_-}, ( \bW_{\mu_+}^{(1)} \otimes_{R} R_{\le m} )^* \otimes \C_{-\La_0} \right) \\
& \cong \varinjlim_m \ext^\bullet_{\fC} \left( \bW_{\la_-}, ( \bW_{\mu_+}^{(1)} \otimes_{R} R_{\le m} )^* \otimes \C_{-\La_0} \right),
\end{align*}
where the first isomorphism follows from Corollary~\ref{cor:comm}.

By Proposition~\ref{prop:Wzero}(1), the first entry of $\ext$ in the last expression identifies with $\bW_{\la_+}^{(1)} \otimes \C_{\La_0}$. Hence the limit vanishes completely when $\la_+ + \mu_- \neq 0$, and vanishes in positive degrees when $\la_+ + \mu_- = 0$, by Theorem~\ref{thm:dual}(2) for $k = 1$.

This completes the proof.
\end{proof}

\begin{prop}\label{prop:WDfiltLv1}
Assume Theorem~\ref{thm:critfilt}(1) for $k=1$. For each $\la \in P$, the $\tb$-module $\bW_\la \otimes \C_{\La_0}$ admits a $D^{(1)}$-filtration.
\end{prop}

\begin{proof}
By Theorem~\ref{thm:critfilt}(1) for $k=1$, we need to verify
\begin{equation}
\ext^{1}_{\fB} ( \bW_\la, ( \bD_\mu^{(1)} \otimes \C_{\La_0} )^{\vee}) = 0 \qquad \la, \mu \in P.\label{eqn:extcrit}
\end{equation}
In case $\mu = \mu_+$, this is Proposition~\ref{prop:Wzero}(2) and Lemma~\ref{lem:WWdual}.

We prove the assertion for $\mu = u \mu_+$, where $u \in W$ is taken smallest possible with respect to the Bruhat order, by assuming the validity of~\eqref{eqn:extcrit} for the case of smaller $u$. By Proposition~\ref{prop:Wzero}(3), we have a short exact sequence
\begin{equation}
0 \longrightarrow \ker f \longrightarrow \bW_{u\mu} \stackrel{f}{\longrightarrow} \bD_{u\mu}^{(1)} \otimes \C_{\La_0} \longrightarrow 0,\label{eqn:indWD}
\end{equation}
where $\ker f$ is filtered by $\{\bD_{v \mu}^{(1)} \otimes \C_{\La_0}\} _{v < u}$.

We have
\[
\ext^{\ge 1}_{\fB} ( \bW_\la, ( \ker f )^{\vee}) = 0 = \ext^{\ge 1}_{\fB} ( \bW_\la, ( \bW_{u\mu} )^{\vee})
\]
by induction hypothesis and Lemma~\ref{lem:WWdual}. We have the long exact sequence obtained by applying $\hom_{\fB} (  \bW_\la,  \bullet )$ to~\eqref{eqn:indWD}:
\[
\shom_{\fB} (  \bW_\la,  ( \bW_{u\mu} )^{\vee} ) \stackrel{g}{\longrightarrow} \shom_{\fB} ( \bW_\la,  ( \ker f )^{\vee} ) \to \ext^{1}_{\fB} ( \bW_\la,  ( \bD_{u\mu}^{(1)}\otimes \C_{\La_0} )^{\vee} )  \to 0
\]
If $s_i u < u$, then applying $\mathscr{D}_i$ to~\eqref{eqn:indWD} yields
\[
0 = \mathbb{L}^{-1} \mathscr D_i ( \bD_{u\mu}^{(1)} \otimes \C_{\La_0} ) \to \mathscr{D}_i ( \ker f ) \longrightarrow \bW_{u\mu} \to \mathscr D_i ( \bD_{u\mu}^{(1)} \otimes \C_{\La_0} ) = 0
\]
by Proposition~\ref{prop:Wzero}(5) (and Theorem~\ref{thm:DJ}(4) for $M = \C_{\La_0}$), which forces
\[
\ker f \subset \mathscr{D}_i ( \ker f ) = \bW_{u\mu}.
\]
Thus, the dual form of Theorem~\ref{thm:indsurj} obtained via~\eqref{eqn:transpose} implies that $g$ is surjective.

This completes the induction and proves~\eqref{eqn:extcrit},
hence completing the proof.
\end{proof}

\begin{thm}\label{thm:Lext}
Let $k \in \Z_{\ge 2}$. For each $\La \in P_{\af,k}^+$ and $\mu \in P$, we have
\begin{equation}
	\ext^1_\fB ( {}^{\theta} L ( \La ), ( D^{(k-1)}_{\mu} \otimes \C_{\La_0})^{*} ) = 0.\label{eqn:B-ext}
\end{equation}
\end{thm}

We record several reduction results, which may be regarded as preliminary steps toward the proof of Theorem~\ref{thm:Lext}. The proof of Theorem~\ref{thm:Lext} is given after Lemma~\ref{lem:basecase}.

\begin{prop}\label{prop:Lext-tens}
Theorem~\ref{thm:Lext} holds for level $k$ if
\begin{equation}
\ext^1_\fB \bigl( {}^{\theta} L ( \La' ) \otimes \bW_\la, ( \C_{\La''})^{*} \bigr) = 0\label{eqn:extenh2}
\end{equation}
holds for every $\La', \La'' \in P_{\af, (k-1)}^+$, $\la \in P$.
\end{prop}

\begin{proof}
By~\cite{BM95}, there exist $\La' \in P_{\af,(k-1)}^+$ and $\La_i \in P_{\af,1}^+$ such that
\[
L ( \La ) \subset L ( \La' ) \otimes L ( \La_i ).
\]
Therefore, the assertion~\eqref{eqn:B-ext} follows from the vanishing of
\begin{equation}
\ext^1_\fB \bigl( {}^{\theta} L ( \La' ) \otimes {}^{\theta} L ( \La_i ), ( W^{(k-1)}_{\mu} \otimes \C_{\La_0} )^{*} \bigr) = 0, \label{eqn:extenh}
\end{equation}
for some $\La' \in P_{\af,(k-1)}^+$, $\La_i \in P_{\af,1}^+$, and $\mu \in P^+$.

Here, ${}^{\theta} L ( \La_i )$ admits a decreasing separable filtration whose associated graded is a direct sum of grading shifts of the modules $\{ \bW _\la^{(1)} \}_{\la \in P^+}$, by Proposition~\ref{prop:Lext} and Theorem~\ref{thm:critfilt}(4) for $k = 1$ using $D^{(k)}_{\gamma_-} = W^{(k)}_{\gamma}$ for $\gamma \in P^+$ (see also~\cite[Theorem~4.7]{CK18} for their constructions).

Thus, it suffices to prove
\[
\ext^1_\fB \bigl( {}^{\theta} L ( \La' ) \otimes \bW_\la^{(1)}, ( W^{(k-1)}_{\mu} \otimes \C_{\La_0})^{*} \bigr) = \ext^1_\fB \bigl( {}^{\theta} L ( \La' ) \otimes \bW_{\la_-}, ( W^{(k-1)}_{\mu} )^{*} \bigr) = 0
\]
for each $\La' \in P^+_{\af,(k-1)},\ \la, \mu \in P^+$ in order to establish Theorem~\ref{thm:Lext} for level $k$.

There exists $\La'' \in P_{\af,(k-1)}^+$ and $w \in W_\af$ such that $\mu + (k-1)\La_0 = w \La''$. We have $W^{(k-1)}_{\mu} \cong \mathscr{D}_w ( \C_{\La''} )$. By Theorem~\ref{thm:DJ}(4) and Proposition~\ref{prop:Wzero}(4), we have
\begin{align*}
\ext^1_\fB \bigl( {}^{\theta} L ( \La' ) \otimes \bW_{\la_-}, ( W^{(k-1)}_{\mu} )^{*} \bigr) &= \ext^1_\fB \bigl( {}^{\theta} L ( \La' ) \otimes \bW_{\la_-}, \mathscr D_{w} ( \C_{\La''})^{*} \bigr)\\
& \cong \ext^1_\fB \bigl( {}^{\theta} L ( \La' ) \otimes  \mathscr D_{w} ( \bW_{\la_-} ), (\C_{\La''})^{*} \bigr) \\
& = \ext^1_\fB \bigl( {}^{\theta} L ( \La' ) \otimes \bW_{u \la_-}, (\C_{\La''})^{*} \bigr)
\end{align*}
up to grading shift for some $u \in W$.

This proves the result.
\end{proof}

\begin{cor}\label{cor:Lext-tens}
Theorem~\ref{thm:Lext} holds for level $k$ if
\begin{equation}
\ext^1_\fB ( {}^{\theta} L ( \La' ) \otimes D_\la^{(1)}, ( \C_{\La'' - \La_0} )^{*} ) = 0\label{eqn:extenh3}
\end{equation}
holds for every $\La', \La'' \in P_{\af, (k-1)}^+$, $\la \in P$. Moreover, the weights $\La'$ and $\La''$ in~\eqref{eqn:extenh3} are the same as those in~\eqref{eqn:extenh2}. Hence, apart from the choice of $\la$, it suffices to verify~\eqref{eqn:extenh3} for the same pairs $(\La',\La'')$ as in~\eqref{eqn:extenh2}.
\end{cor}

\begin{proof}
By Proposition~\ref{prop:WDfiltLv1}, we have
\[
\bW_\la \otimes \C_{\La_0} \cong \varprojlim M(r),
\]
where $\{M(r)\}_r$ is an inverse system with surjective transition maps such that each $M(r)$ is finite-dimensional and admits $D^{(1)}$-filtration. Since $\bW_\la \in \fB_\bd$, the degree $\le m$-part of $\bW_\la \otimes \C_{\La_0}$ is the same as that of $M(r)$ for $r \gg 0$ for each $m \in \Z$.

Thus we have
\begin{equation}
\ext^1_\fB \bigl( {}^{\theta} L ( \La' ) \otimes \bW_\la, ( \C_{\La''})^{*} \bigr) \cong \varinjlim_r \ext^1_\fB \bigl( {}^{\theta} L ( \La' ) \otimes M(r), ( \C_{\La''})^{*} \bigr)\label{eqn:WMproj}
\end{equation}
by Lemma~\ref{lem:comm}. Here the right-hand side of~\eqref{eqn:WMproj} vanishes by assumption and Proposition~\ref{prop:Lext-tens}.

This completes the proof.
\end{proof}

\begin{lem}\label{lem:basecase}
Theorem~\ref{thm:Lext} holds for level $k=2$.
\end{lem}

\begin{proof}
It suffices to verify the condition~\eqref{eqn:extenh3} in Corollary~\ref{cor:Lext-tens} for $k=2$. To this end, we apply an affine Dynkin diagram automorphism which sends $\La''$ to $\La_0$ in~\eqref{eqn:extenh2}. This reduces the choice of $\La'' \in P_{\af,1}^+$ to $\La_0$. 

We have
\begin{equation}
\ext^1_\fB \bigl( {}^{\theta} L ( \La' ) \otimes D_\la^{(1)}, ( \C_{\La_0 - \La_0} )^{*} \bigr)
= \ext^1_\fB \bigl( {}^{\theta} L ( \La' ), ( D_\la^{(1)} )^* \bigr) \label{eqn:k2}
\end{equation}
for each $\La' \in P_{\af,1}^+$. This vanishes by Proposition~\ref{prop:Lext}.

This proves the assertion.
\end{proof}

\begin{proof}[Proof of Theorem~\ref{thm:Lext}]
We prove~\eqref{eqn:B-ext} by induction on $k$. 
The case $k = 2$ is proved in Lemma~\ref{lem:basecase}.

\medskip

We consider the case of $k > 2$ by assuming the validity of the assertion for lower values of $k$. We now find $x \in W_\af$ and $\La_j \in P^+_{\af,1}$ such that
\[
\la + \La_0 = x \La_j \quad \text{and hence} \quad D_{\la}^{(1)} \cong \mathscr D_{x} ( \C_{\La_j} ).
\]

We rewrite~\eqref{eqn:extenh3} as
\begin{align}
\ext^1_\fB \bigl( {}^{\theta} L ( \La' ) \otimes D_{\la}^{(1)}, ( \C_{\La'' - \La_0} )^{*} \bigr) &= \ext^1_\fB \bigl( {}^{\theta} L ( \La' ) \otimes \mathscr D_{x} ( \C_{\La_j} ), ( \C_{\La'' - \La_0} )^{*} \bigr) \label{eqn:lastone}  \\
&\cong \ext^1_\fB \bigl( {}^{\theta} L ( \La' ) \otimes \C_{\La_j}, ( \mathbb L^{\bullet} \mathscr D_{x^{-1}} ( \C_{\La'' - \La_0} ) )^{*} \bigr) \nonumber \\
&\cong \ext^1_\fB \bigl( {}^{\theta} L ( \La' ), ( \mathbb L^{\bullet} \mathscr D_{x^{-1}} ( \C_{\La'' - \La_0} ) \otimes \C_{\La_j} )^{*} \bigr). \nonumber
\end{align}

To the last expression, we may apply an affine Dynkin diagram automorphism that sends $\La_j$ to $\La_0$.

Therefore, if the derived Demazure functor is concentrated in degree $0$ on
\[
\mathbb L^{\bullet} \mathscr D_{x^{-1}} ( \C_{\La'' - \La_0} )
\]
and defines an actual Demazure module, then we obtain
\[
\ext^1_\fB \bigl( {}^{\theta} L ( \La' ), ( \mathscr D_{x^{-1}} ( \C_{\La'' - \La_0} ) \otimes \C_{\La_j} )^{*} \bigr) = 0
\]
from the level $(k{-}1)$ case of~\eqref{eqn:B-ext}.

This scenario applies when $\La'' - \La_0 \in P_{\af,(k-2)}^+$.

Thus, we are reduced to the case where $\La'' - \La_0 \not\in P_{\af,(k-2)}^+$ (but still $\La'' \in P_{\af, (k-1)}^+$). We have $\langle \al_0^{\vee}, \La'' - \La_0 \rangle = -1$, and hence
\[
s_0 ( \overline{\La''} + (k-2)\La_0  ) = ( \overline{\La''} + (k-2)\La_0  ) + \al_0.
\]
We take a minimal length $w \in W_\af$ such that $w ^{-1} s_0 ( \overline{\La''} + (k-2)\La_0  ) \in P_{\af,(k-2)}^+$ and apply Proposition~\ref{prop:Demext} to obtain a short exact sequence of $\tb$-modules:
\begin{equation}
0 \longrightarrow \bq D^{(k-2)}_{\overline{\La'' + \al_0}} \longrightarrow D^{(k-2)}_{\overline{\La''}} \longrightarrow \C_{\La'' - \La_0} \longrightarrow 0. \label{eqn:Demext}
\end{equation}

As in~\eqref{eqn:lastone}, we have
\begin{align*}
\ext^1_\fB \bigl( {}^{\theta} L ( \La' ) \otimes D_{\la}^{(1)}, ( D^{(k-2)}_{\overline{\La''}} )^{*} \bigr)
&\equiv \ext^1_\fB \bigl( {}^{\theta} L ( \La' ) \otimes \mathscr D_{x} ( \C_{\La_j} ), ( D^{(k-2)}_{\overline{\La''}} )^{*} \bigr) \\
&\cong \ext^1_\fB \bigl( {}^{\theta} L ( \La' ) \otimes \C_{\La_j}, ( \mathscr D_{x^{-1}} ( D^{(k-2)}_{\overline{\La''}} ) )^{*} \bigr) \\
&\cong \ext^1_\fB \bigl( {}^{\theta} L ( \La' ), ( \mathscr D_{x^{-1}} ( D^{(k-2)}_{\overline{\La''}} ) \otimes \C_{\La_j} )^{*} \bigr).
\end{align*}

We apply an affine Dynkin diagram automorphism that sends $\La_j$ to $\La_0$, using Theorem~\ref{thm:adj} and Corollary~\ref{cor:CKmain2}.  
Then, by applying~\eqref{eqn:B-ext} for level $(k{-}1)$, we deduce
\[
\ext^1_\fB \bigl( {}^{\theta} L ( \La' ) \otimes D_{\la}^{(1)}, ( D^{(k-2)}_{\overline{\La''}} )^{*} \bigr) = 0.
\]

We now consider the effect of applying the functor
\[
\shom_\fB \bigl( {}^{\theta} L ( \La' ) \otimes \mathscr D_{x} ( \C_{\La_j} ), ( \bullet )^{*} \bigr)
\]
to the short exact sequence~\eqref{eqn:Demext}. The vanishing of~\eqref{eqn:lastone}, namely
\[
\ext^1_\fB \bigl( {}^{\theta} L ( \La' ) \otimes D_{\la}^{(1)}, ( \C_{\La'' - \La_0} )^{*} \bigr) = 0,
\]
is then equivalent to the surjectivity of the map $\psi$:
\begin{align*}
&\shom_\fB \bigl( {}^{\theta} L ( \La' ) \otimes D_{\la}^{(1)}, ( D^{(k-2)}_{\overline{\La''}} )^{*} \bigr)
\stackrel{\psi}{\longrightarrow} \shom_\fB \bigl( {}^{\theta} L ( \La' ) \otimes D_{\la}^{(1)}, ( \bq D^{(k-2)}_{\overline{\La''+\al_0}} )^{*} \bigr) \\
&\hspace{10mm} \longrightarrow \ext^1_\fB \bigl( {}^{\theta} L ( \La' ) \otimes D_{\la}^{(1)}, ( \C_{\La'' - \La_0} )^{*} \bigr)
\longrightarrow \ext^1_\fB \bigl( {}^{\theta} L ( \La' ) \otimes D_{\la}^{(1)}, ( D^{(k-2)}_{\overline{\La''}} )^{*} \bigr) = 0.
\end{align*}

We have $D_\la^{(1)} \subset \mathscr D_0 ( D_\la^{(1)} )$ and $\mathbb L^{<0} \mathscr D_0 ( D_\la^{(1)} ) = 0$ as $D_\la^{(1)}$ is a level-one Demazure module. Theorem~\ref{thm:DJ} implies
\[
{}^{\theta} L ( \La' ) \otimes D_\la^{(1)} \subset \mathbb L^{\bullet} \mathscr D_0 \bigl( {}^{\theta} L ( \La' ) \otimes D_\la^{(1)} \bigr)
\cong {}^{\theta} L ( \La' ) \otimes \mathbb L^{\bullet} \mathscr D_0 ( D_\la^{(1)} ).
\]

Moreover, we have
\[
\mathbb L^{\bullet} \mathscr D_0 ( \bq D^{(k-2)}_{\overline{\La'' + \al_0}} ) = D^{(k-2)}_{\overline{\La''}}.
\]
Hence, by Theorem~\ref{thm:adj}, the map $\psi$ coincides with
\[
\shom_\fB \bigl( \mathscr D_0 ( {}^{\theta} L ( \La' ) \otimes D_{\la}^{(1)} ), ( D^{(k-2)}_{\overline{\La''+\al_0}} )^{*} \bigr)
\longrightarrow \shom_\fB \bigl( {}^{\theta} L ( \La' ) \otimes D_{\la}^{(1)}, ( D^{(k-2)}_{\overline{\La''+\al_0}} )^{*} \bigr).
\]

Hence, Corollary~\ref{cor:indsurj} asserts that $\psi$ must be surjective.

Therefore,~\eqref{eqn:B-ext} for level $(k-1)$ implies~\eqref{eqn:B-ext} for level $k$. This completes the induction and proves the assertion.
\end{proof}

\section{Some filtration results}\label{sec:filt}

We retain the setting of the previous section.

\begin{prop}\label{prop:passage}
Let $k, l \in \Z_{> 0}$ be such that $k \ge l$. Assume that Theorem~\ref{thm:critfilt} holds for all levels $< k$. Then for each $\la \in P$, the module $\bD_\la^{(k)} \otimes \C_{(k-l)\La_0}$ admits a $\bD^{(l)}$-filtration.
\end{prop}

\begin{proof}
The case $k = l$ is trivial by definition. Thus, we assume $k > l$ in the following.
By Theorem~\ref{thm:critfilt}(2) at level $l$, which is available since $l < k$, the assertion is equivalent to
\[
\ext^{1}_{\fB} \bigl( \bD^{(k)}_\la \otimes \C_{(k-l)\La_0}, ( D_{\mu}^{(l)} )^{*} \bigr) = 0
\qquad \text{for all } \la, \mu \in P.
\]
We first treat the case $k = l+1$, in which this condition reads
\begin{equation}
\ext^{1}_{\fB} \left( \bD^{(k)}_\la, ( D_{\mu}^{(k-1)} \otimes \C_{\La_0} )^{*} \right) = 0
\hskip 5mm \text{for all } \la, \mu \in P;\label{eqn:extdualone}
\end{equation}
the general case is deduced from this at the end of the proof.

For $\La \in P_{\af, k}^+$, let $W_\af^{\La} \subset W_\af$ denote the set of maximal length coset representatives for $W_\af / \mathsf{Stab}_{W_\af} ( \La )$. For each $\La \in P_{\af, k}^+$, $m \in \Z_{\ge 0}$, and $i \in \tI_\af$, define
\[
L^{m} ( \La ) := \sum_{\ell(v) \ge m} {}^{\theta} L ( \La )^v, \qquad
L_i^m ( \La ) := L^{m+1} ( \La ) + \sum_{\substack{\ell(v) = m \\ s_i v < v}} {}^{\theta} L ( \La )^v,
\]
where $v$ ranges over $W_\af^{\La}$ subject to the indicated conditions.

By Theorem~\ref{thm:tD} and Corollary~\ref{cor:tD}, we have
\[
\mathbb L^{\bullet} \mathscr D_i ( L^{m} ( \La ) )
= \sum_{\substack{w \in W_\af^{\La} \\ \ell(w) = m{-}1 \\ s_i w > w}} {}^{\theta} L( \La )^w + L^{m} ( \La ).
\]
In view of Theorem~\ref{thm:D-incl} (see also~\cite[Theorem~C]{Kat18b}), it follows that
\begin{equation}
\frac{\mathscr D_i ( L^{m} ( \La ) )}{L^{m} ( \La )}
\cong \bigoplus_{\substack{w \in W_\af^{\La} \\ \ell(w) = m{-}1 \\ s_i w > w}} \bq^{\langle d, \La - w \La \rangle} \bD^{(k)}_{- \overline{w \La}}.
\label{eqn:dec}
\end{equation}
In particular, we deduce that
\begin{equation}
\mathbb L^{\bullet} \mathscr D_i ( L^{m} ( \La ) )
= \mathbb L^{\bullet} \mathscr D_i ( L_i^{m} ( \La ) ),\label{eqn:DLeq}
\end{equation}
since the condition $s_i v \not< v$ implies $s_i v > v$ for each $v \in W_\af^{\La}$.

We prove the following assertion by induction on $m$:
\begin{itemize}
  \item[$(\bigstar)_1^m$] $L^m ( \La ) \otimes \C_{\La_0}$ admits a $\bD^{(k-1)}$-filtration;
  \item[$(\bigstar)_2^m$] For every $w \in W_\af^{\La}$ such that $\ell(w) < m$, $\bD^{(k)}_{-\overline{w \La}} \otimes \C_{\La_0}$ admits a $\bD^{(k-1)}$-filtration.
\end{itemize}

We may refer $(\bigstar)^m_1$ and $(\bigstar)_2^m$ as $(\bigstar)^m$ for each $m$. Here we note that $(\bigstar)_1^m$ is equivalent to
\[
\ext^1_{\fB} \left( L^{m} ( \La )\otimes \C_{\La_0}, ( D^{(k-1)}_{\mu} )^* \right) \cong \ext^1_{\fB} \left( L^{m} ( \La ), ( D^{(k-1)}_{\mu} \otimes \C_{\La_0} )^* \right) = 0,
\]
and $(\bigstar)_2^m$ is equivalent to
\[
\ext^1_{\fB} \left( \bD^{(k)}_{-\overline{w \La}} \otimes \C_{\La_0}, ( D^{(k-1)}_{\mu} )^* \right) \cong \ext^1_{\fB} \left( \bD^{(k)}_{-\overline{w \La}}, ( D^{(k-1)}_{\mu} \otimes \C_{\La_0} )^* \right) = 0
\]
by Theorem~\ref{thm:critfilt}(2).

The base case $(\bigstar)^0_1$ holds by Theorem~\ref{thm:Lext} and the base case $(\bigstar)^0_2$ holds by the definition of $\{ \bD^{(k)}_{-\overline{w \La}} \}_{w \in W_\af^{\La}}$ since $(\bigstar)_2^0$ is vacuous.

By Theorem~\ref{thm:tD} and Corollary~\ref{cor:tD}, we have a short exact sequence
\[
0 \longrightarrow L^m ( \La ) \longrightarrow L^{m-1} ( \La )
\longrightarrow \bigoplus_{\substack{w \in W_\af^{\La} \\ \ell(w) = m{-}1}} \bq^{\langle d,\, \La - w \La \rangle} \bD^{(k)}_{- \overline{w \La}} \longrightarrow 0.
\]
Applying Lemma~\ref{lem:surjfilt}, we see that $(\bigstar)_1^{m-1}$ and $(\bigstar)_2^m$ together imply $(\bigstar)_1^m$.

Thus, it remains to prove $(\bigstar)_2^m$ assuming $(\bigstar)^{m-1}$.

Assume now that $(\bigstar)^{m-1}$ holds. By~\eqref{eqn:dec}, and $(\bigstar)^{m-1}$, we have
\[
\ext^{1}_{\fB} ( L^{m-1} ( \La ), ( D^{(k-1)}_{\mu} \otimes \C_{\La_0} )^* ) = 0 = \bigoplus_{\substack{w \in W_\af^{\La} \\ \ell(w) = m{-}2}} \ext^{1}_{\fB} ( \bD^{(k)}_{- \overline{w \La}}, ( D^{(k-1)}_{\mu} \otimes \C_{\La_0} )^* )
\]
for every $\mu \in P$. By the long exact sequence of $\ext$ applied to
\[
0 \rightarrow L^{m-1} ( \La ) \rightarrow \mathscr{D}_i ( L^{m-1} ( \La ) ) \rightarrow \bigoplus_{\substack{w \in W_\af^{\La} \\ \ell(w) = m{-}2 \\ s_i w > w}} \bq^{\langle d, \La - w \La \rangle} \bD^{(k)}_{- \overline{w \La}}\rightarrow 0
\]
and~\eqref{eqn:DLeq}, we obtain
\begin{equation}
\ext^{1}_{\fB} ( \mathscr D_i ( L_i^{m-1} ( \La ) ), ( D^{(k-1)}_{\mu} \otimes \C_{\La_0} )^* ) = \ext^{1}_{\fB} ( \mathscr D_i ( L^{m-1} ( \La ) ), ( D^{(k-1)}_{\mu} \otimes \C_{\La_0} )^* ) = 0\label{eqn:ext1D}
\end{equation}
for every $\mu \in P$.

By Corollary~\ref{cor:tD} and Corollary~\ref{cor:indsurj}, we have a surjection
\[
\shom_{\fB} \bigl( \mathscr D_i ( L_i^{m-1} ( \La ) ), ( D^{(k-1)}_{\mu} \otimes \C_{\La_0} )^* \bigr)
\twoheadrightarrow
\shom_{\fB} \bigl(  L_i^{m-1} ( \La ), ( D^{(k-1)}_{\mu} \otimes \C_{\La_0} )^* \bigr)
\]
for each $i \in \tI_\af$ and $\mu \in P$.

This yields a portion of the long exact sequence:
\begin{align*}
\shom_{\fB} & \left( \mathscr D_i ( L_i^{m-1} ( \La ) ), ( D^{(k-1)}_{\mu} \otimes \C_{\La_0} )^* \right) \longrightarrow \shom_{\fB} \left( L_i^{m-1} ( \La ), ( D^{(k-1)}_{\mu} \otimes \C_{\La_0} )^* \right) \\
&\longrightarrow\ext^1_{\fB} \left( \frac{\mathscr D_i ( L_i^{m-1} ( \La ) )}{L_i^{m-1} ( \La )}, ( D^{(k-1)}_{\mu} \otimes \C_{\La_0} )^* \right) \\
& \longrightarrow \ext^1_{\fB} \left( \mathscr D_i ( L_i^{m-1} ( \La ) ), ( D^{(k-1)}_{\mu} \otimes \C_{\La_0} )^* \right) = 0,
\end{align*}
associated to the short exact sequence
\begin{equation}
0 \longrightarrow L_i^{m-1} ( \La ) \longrightarrow \mathscr D_i ( L_i^{m-1} ( \La ) ) \longrightarrow \frac{\mathscr D_i ( L_i^{m-1} ( \La ) )}{L_i^{m-1} ( \La )} \longrightarrow 0.
\label{eqn:SES-DL}
\end{equation}

From this, we deduce that
\[
\frac{\mathscr D_i ( L^{m-1} ( \La ) )}{L_i^{m-1} ( \La )} \otimes \C_{\La_0}
\]
admits a $\bD^{(k-1)}$-filtration by Theorem~\ref{thm:critfilt} at level $k{-}1$.

Moreover, we have a short exact sequence
\[
0 \longrightarrow \frac{L^{m-1} ( \La )}{L_i^{m-1} ( \La )}
\longrightarrow \frac{\mathscr D_i ( L_i^{m-1} ( \La ) )}{L_i^{m-1} ( \La )}
\longrightarrow
\bigoplus_{\substack{w \in W_\af^{\La} \\ \ell(w) = m{-}2 \\ s_i w > w}}
\bq^{\langle d,\, \La - w \La \rangle} \bD^{(k)}_{- \overline{w \La}} \longrightarrow 0.
\]

By $(\bigstar)^{m-1}$ and Lemma~\ref{lem:surjfilt}, we conclude that
\[
\frac{L^{m-1} ( \La )}{L_i^{m-1} ( \La )} \otimes \C_{\La_0}
\]
admits a $\bD^{(k-1)}$-filtration. In fact, we have the decomposition
\begin{equation}
\frac{L^{m-1} ( \La )}{L_i^{m-1} ( \La )}
\cong
\bigoplus_{\substack{w \in W_\af^{\La} \\ \ell(w) = m{-}1 \\ s_i w > w}} \bq^{\langle d,\, \La - w \La \rangle} \bD^{(k)}_{- \overline{w \La}}.
\label{eqn:idirect}
\end{equation}

Since every $w \in W_\af$ satisfies $s_i w > w$ for some $i \in \tI_\af$, it follows that every $\bD^{(k)}_{- \overline{w \La}}$ with $\ell(w) = m{-}1$ appears on the right-hand side of~\eqref{eqn:idirect} for some $i \in \tI_\af$. Thus, Corollary~\ref{cor:direct} implies $(\bigstar)_2^m$ from $(\bigstar)^{m-1}$.

This completes the induction, and we obtain $(\bigstar)$.

Therefore, each $\bD^{(k)}_\la \otimes \C_{\La_0}$ admits a $\bD^{(k-1)}$-filtration. This settles the case $k = l+1$.

For $k > l+1$, we argue by induction on $k - l$. By the induction hypothesis applied to $(k-1, l)$, every factor $\bD^{(k-1)}_\nu \otimes \C_{(k-1-l)\La_0}$ of the above filtration admits a $\bD^{(l)}$-filtration, so that
\[
\ext^1_\fB \bigl( \bD^{(k-1)}_\nu \otimes \C_{(k-1-l)\La_0}, ( D^{(l)}_\mu )^* \bigr) = 0 \qquad (\mu \in P)
\]
by Theorem~\ref{thm:critfilt}(2) at level $l$. Passing to $\bD^{(k)}_\la \otimes \C_{(k-l)\La_0}$ by d\'evissage and Lemma~\ref{lem:comm}, as in the proof of Corollary~\ref{cor:Dhwext}(2), we obtain the same vanishing there, and Theorem~\ref{thm:critfilt}(2) at level $l$ yields the assertion.
\end{proof}

\begin{cor}\label{cor:passage}
Let $k, l \in \Z_{> 0}$ be such that $k \ge l$. Assume that Theorem~\ref{thm:critfilt} holds for all levels $< k$. Then for each $\la \in P^{+}$, the module $\bW_\la^{(k)} \otimes \C_{(k-l)\La_0}$ admits a $\bW^{(l)}$-filtration.
\end{cor}

\begin{proof}
We have $\bW_\la^{(k)} = \mathscr D_{w_0} ( \bW_{\la} ^{(k)} ) = \mathscr D_{w_0} ( \bD_{\la} ^{(k)})$ by Lemma~\ref{lem:DtoWbyD}, and $W^{(k)}_\mu = D^{(k)}_{\mu_-}$ for $\mu \in P^+$ by definition. Therefore, from~\eqref{eqn:extdualone} borrowed from the proof of Proposition~\ref{prop:passage}, we deduce
\begin{align*}
0 &= \ext^{1}_{\fB} \left( \bD^{(k)}_{\la}, ( D_{\mu_-}^{(k-1)} \otimes \C_{\La_0} )^{*} \right) \\
& = \ext^{1}_{\fB} \left( \bD^{(k)}_{\la}, ( \mathscr D_{w_0} ( D_{\mu_-}^{(k-1)} \otimes \C_{\La_0} ) )^{*} \right) \\
& = \ext^{1}_{\fB} \left( \mathscr D_{w_0} ( \bD^{(k)}_{\la} ), ( D_{\mu_-}^{(k-1)} \otimes \C_{\La_0} )^{*} \right)\\
& = \ext^{1}_{\fB} \left( \bW^{(k)}_{\la}, ( W_{\mu}^{(k-1)} \otimes \C_{\La_0} )^{*} \right),
\end{align*}
where the second equality follows from the $\g$-invariance of $D_{\mu_-}^{(k-1)}$ and $\C_{\La_0}$, the third equality is Theorem~\ref{thm:adj}, and the fourth equality is the above identification. Now we apply Theorem~\ref{thm:critfilt}(4)
 to deduce the assertion for $l = k {-} 1$. Now a repeated application of this case yields the general case.
\end{proof}

\begin{thm}\label{thm:orth}
Let $k \in \Z_{> 0}$, and assume that Theorem~\ref{thm:critfilt} holds for all levels $< k$. Then for each $\La \in P^{+}_{\af, k}$ and $w \in W_\af$, we have
\begin{equation}
\ext^{\bullet}_{\fB} \left( {}^\theta L ( \La )^w, ( \C_{\mu + k\La_0} )^* \right) = 0 \qquad \text{whenever } \mu \not\succeq \overline{w \La}.\label{eqn:LCext}
\end{equation}
In particular, for every $k' \in \Z_{> 0}$ and $\mu \in P$, we have
\[
\ext^{\bullet}_{\fB} \left( {}^\theta L ( \La )^w, ( D_{\mu}^{(k')} \otimes \C_{(k - k') \La_0} )^* \right) = 0
\qquad \text{whenever } \mu \not\succeq \overline{w \La}.
\]
\end{thm}

\begin{proof}
We first prove~\eqref{eqn:LCext}.
By Corollary~\ref{cor:LtoDfilt}, the module ${}^\theta L ( \La )^w$ admits a
$\bD^{(k)}$-filtration, and by repeatedly applying Proposition~\ref{prop:passage}
this filtration can be refined into a $\bD^{(1)}$-filtration. We thus obtain a
projective system $\{ M(r) \}_r$ with surjective transition maps such that
\[
{}^\theta L ( \La )^w \cong \varprojlim_r M(r),
\]
where each $M(r)$ admits a finite filtration whose associated graded components
are grading shifts of modules of the form $\bD^{(1)}_{\nu}$.

We claim that every such factor $\bD^{(1)}_{\nu}$ satisfies
\begin{equation}
\overline{w \La} \preceq - \nu. \label{eqn:orth-piece}
\end{equation}
Indeed, since $\bD^{(1)}_{\nu}$ occurs in the $\bD^{(1)}$-filtration of
${}^\theta L ( \La )^w$, Corollary~\ref{cor:multcrit}(1) at level one, which follows from~\cite[Theorems~4.19 and~5.9(ii)]{CK18} in place of Theorems~\ref{thm:dual} and~\ref{thm:critfilt} at level one, gives
\[
0 \neq \bigl( {}^\theta L ( \La )^w : \bD^{(1)}_{\nu} \bigr)_{q^{-1}}
= \gdim \hom_{\fB} \bigl( {}^\theta L ( \La )^w, ( D^{(1)}_{-\nu} )^* \bigr).
\]
As ${}^\theta L ( \La )^w$ is a cyclic $\tb$-module generated by a vector of
$\wth$-weight $-w\La$, any nonzero homomorphism
${}^\theta L ( \La )^w \to ( D^{(1)}_{-\nu} )^*$ is nonzero on this generator,
so $-\overline{w\La}$ is a $\h$-weight of $( D^{(1)}_{-\nu} )^*$; equivalently,
$w\La \in \Psi ( D^{(1)}_{-\nu} )$. By the weight estimate in
Theorem~\ref{thm:Dorth}, we have
$\overline{\Psi ( D^{(1)}_{-\nu} )} \subset \Sigma ( -\nu )$, whence
$\overline{w\La} \preceq -\nu$. This proves~\eqref{eqn:orth-piece}.

Now assume $\mu \not\succeq \overline{w \La}$. If some factor $\bD^{(1)}_{\nu}$ satisfied
$-\nu \preceq \mu$, then combining with~\eqref{eqn:orth-piece} and the
transitivity of $\preceq$ we would obtain
$\overline{w\La} \preceq -\nu \preceq \mu$, contradicting
$\overline{w\La} \not\preceq \mu$. Hence $-\nu \not\preceq \mu$ for every factor,
and Lemma~\ref{lem:extDC} yields
\[
\ext^{\bullet}_{\fB} \bigl( \bD^{(1)}_{\nu} \otimes \C_{\La_0}, ( \C_{\mu} )^* \bigr) = 0
\]
for each factor. Applying d\'evissage to the finite $\bD^{(1)}$-filtration of
$M(r)$, we conclude
\[
\ext_{\fB}^{\bullet} \bigl( M(r), ( \C_{\mu + k \La_0} )^* \bigr) = 0
\qquad (r \ge 0).
\]
Passing to the limit via Lemma~\ref{lem:comm}, we obtain
\[
\ext_{\fB}^{\bullet} \bigl( {}^\theta L ( \La )^w, ( \C_{\mu + k \La_0} )^* \bigr)
\cong \varinjlim_r
\ext_{\fB}^{\bullet} \bigl( M(r), ( \C_{\mu + k \La_0} )^* \bigr) = 0,
\]
which proves~\eqref{eqn:LCext}.

We next deduce the second assertion. By Lemma~\ref{lem:we} and the definition of
Demazure modules, $D^{(k')}_\mu$ admits a finite filtration by one-dimensional
$\tb$-modules of the form $\bq^m \C_{\gamma + k'\La_0}$ with $\gamma \preceq \mu$
and $m \in \Z$, starting from the cyclic quotient $\C_{\mu + k'\La_0}$; the
filtration is finite since $\dim D^{(k')}_\mu < \infty$. For each such
$\gamma \preceq \mu$, the assumption $\mu \not\succeq \overline{w\La}$ and the
transitivity of $\preceq$ give $\gamma \not\succeq \overline{w\La}$. Applying the functor
\[
\ext_{\fB}^{\bullet} \bigl( {}^\theta L ( \La )^w, ( \bullet )^* \bigr)
\]
to the short exact sequences of this filtration and using~\eqref{eqn:LCext} at
each step, we deduce inductively that
\[
\ext^{\bullet}_{\fB} \bigl( {}^\theta L ( \La )^w,
( D_{\mu}^{(k')} \otimes \C_{(k - k') \La_0} )^* \bigr) = 0
\qquad \text{whenever } \mu \not\succeq \overline{w\La}.
\]
This completes the proof.
\end{proof}

\section{Proofs of the main theorems}\label{sec:proof}

We retain the setting of the previous section. This section is devoted entirely to the proofs of Theorems~\ref{thm:dual} and~\ref{thm:critfilt}. 
Since our approach proceeds by induction on the level $k$ (with the base case reducing to~\cite{CK18}), 
we assume that both Theorem~\ref{thm:dual} and Theorem~\ref{thm:critfilt} hold for all levels strictly less than $k$.

The combination of Theorems~\ref{thm:dual-ind} and~\ref{thm:crit-ind} then establishes 
Theorems~\ref{thm:dual} and~\ref{thm:critfilt} for arbitrary $k$ by induction from the case $k=1$.

\begin{thm}\label{thm:dual-ind}
Theorem~\ref{thm:dual} for level $<k$ and Theorem~\ref{thm:critfilt} for level $<k$ imply Theorem~\ref{thm:dual} for level $k$.
\end{thm}

\begin{proof}
We prove Theorem~\ref{thm:dual}(1) for level $k$. Let $\la,\mu\in P$. Choose $\La,\La' \in P_{\af,k}^+$ and $w,v\in W_\af$ such that $\la + k \La_0 = w \La$ and $\mu + k \La_0 = v \La'$. By rearranging $w$ and $v$ if necessary, we may assume that $w$ is the minimal element in $w \mathsf{Stab}_{W_\af} (\La)$, and that $v$ is the maximal element in $v \mathsf{Stab}_{W_\af} (\La')$; the latter exists since $\mathsf{Stab}_{W_\af}(\La') \subset W_\af$ is a finite subgroup by $\La' (K) > 0$.

By Theorem~\ref{thm:adj}, we have
\begin{align*}
\ext^{\bullet}_\fB ( {}^{\theta} L ( \La )^w, ( D_\mu ^{(k)} )^{*} ) 
&\cong \ext^{\bullet}_\fB ( {}^{\theta} L ( \La )^w, \mathscr D_v ( \C _{\La'} )^{*} ) \\
&\cong \ext^{\bullet}_\fB \left( \mathscr D_{v^{-1}} \left( {}^{\theta} L ( \La )^w \right), ( \C _{\La'} )^{*} \right).
\end{align*}
	
Thanks to Theorem~\ref{thm:tD}, we have
\[
\mathscr D_{v^{-1}} \left( {}^{\theta} L ( \La )^w \right) = {}^{\theta} L ( \La )
\]
if and only if a reduced expression of $v^{-1}$ contains a subexpression whose product equals $w^{-1}$. This occurs precisely when $w^{-1} \le v^{-1}$ in the Bruhat order, which is equivalent to $w \le v$. Hence, we conclude that
\begin{equation}
\ext^{\bullet}_\fB ( {}^{\theta} L ( \La )^w, ( D_\mu ^{(k)} )^{*} ) \cong 
\begin{cases} 
\ext^{\bullet}_\fB ( {}^{\theta} L ( \La ), ( \C _{\La'} )^{*} ) & \text{if } w \le v, \\
\ext^{\bullet}_\fB ( {}^{\theta} L ( \La )^u, ( \C _{\La'} )^{*} ) & \text{otherwise},
\end{cases}
\label{exttop}
\end{equation}
for some $u \in W_\af$ with $u \ne e$ and $u \La \ne \La$, by Proposition~\ref{prop:ordercomp}.

By Theorem~\ref{thm:orth} (which relies on Theorem~\ref{thm:critfilt} for levels $<k$), together with the fact that the set of classical parts of weights in $W P_{\af,k}^+$ ($= \overline{W P_{\af,k}^+}$) is closed under taking smaller elements with respect to $\prec$, we find that
\begin{equation}
\ext^{\bullet}_\fB ( {}^{\theta} L ( \La )^u, ( \C _{\La'} )^{*} ) \equiv 0 \qquad \text{if } u\La \notin W\La.\label{eqn:bigvanish}
\end{equation}

We now examine the case $u \La \in W \La$, which complements~\eqref{eqn:bigvanish}.

We apply Proposition~\ref{prop:passage} to ${}^{\theta} L ( \La )^u$ (which requires Theorem~\ref{thm:critfilt} for levels $< k$), thereby obtaining a decreasing separable filtration of ${}^{\theta} L ( \La )^u$ whose associated graded pieces are grading shifts of $\{ \bD^{(k-1)}_{\mu} \otimes \C_{-\La_0} \}_{\mu \in P}$. 
Hence we may write
\[
{}^{\theta} L ( \La )^u \;\cong\; \varprojlim_r M(r),
\]
where $\{M(r)\}_r$ is a projective system with surjective transition maps such that each $M(r) \otimes \C_{\La_0}$ admits a finite $\bD^{(k-1)}$-filtration.

In the case $\La' - \La_0 \in P_{\af,(k-1)}^+$, we have
\[
\ext^{>0}_\fB ( \bD^{(k-1)}_{\mu} \otimes \C_{- \La_0}, ( \C _{\La'} )^{*} ) \cong \ext^{>0}_\fB ( \bD^{(k-1)}_{\mu}, ( \C _{\La' - \La_0} )^{*} ) = 0
\]
by Theorem~\ref{thm:dual} for level $k{-}1$ and the identification
\[
\C_{\La'} \cong D^{(k)}_{\overline{\La'}} \cong D^{(k-1)}_{\overline{\La'}} \otimes \C_{\La_0}.
\]
It follows that
\[
\ext^{>0}_\fB ( {}^{\theta} L ( \La )^u, ( \C _{\La'} )^{*} ) = \varinjlim_r \ext^{>0}_\fB ( M(r), ( \C _{\La'} )^{*} ) \equiv 0 \qquad \text{for } u \in W
\]
by Lemma~\ref{lem:comm}.

In the case where $\La' - \La_0 \notin P_{\af,(k-1)}^+$, we necessarily have $\langle \vartheta^{\vee}, \overline{\La'} \rangle = k$. 
In this situation, we apply Proposition~\ref{prop:Demext} (with $\La$ replaced by $\La'$) and note that 
$\La' = (\La'-\La_0) + \La_0$ and $\overline{\La'} = \overline{\La' - \La_0}$. 
It then follows that
\begin{align}\nonumber
\ext^{i}_\fB \bigl( {}^{\theta} L ( \La )^u, ( \bq D ^{(k-1)}_{\overline{s_0 (\La'-\La_0)}} \otimes \C_{\La_0} )^{*} \bigr) 
& \longrightarrow 
\ext^{i+1}_\fB \bigl( {}^{\theta} L ( \La )^u, ( \C _{\La'} )^{*} \bigr)\\
\longrightarrow & \ext^{i+1}_\fB \bigl( {}^{\theta} L ( \La )^u, ( D ^{(k-1)}_{\overline{\La'}} \otimes \C_{\La_0} )^{*} \bigr)
\label{eqn:exLDs}
\end{align}
is exact for each $i \ge 0$.

As in the discussion for the case $\La' - \La_0 \in P_{\af,(k-1)}^+$, 
Proposition~\ref{prop:passage} together with Theorem~\ref{thm:dual} for level $k-1$ implies that the middle term $\ext^{i+1}({}^\theta L(\Lambda)^u,(\C_{\Lambda'})^*)$ in~\eqref{eqn:exLDs} vanishes for all $i > 0$. 
By Theorem~\ref{thm:tD}, we have
\begin{equation}
{}^{\theta} L ( \La )^u \subset 
\mathscr{D}_0 \bigl( {}^{\theta} L ( \La )^u \bigr) 
\equiv 
\mathbb{L}^{\bullet}\mathscr{D}_0 \bigl( {}^{\theta} L ( \La )^u \bigr).\label{eqn:tDrestate}
\end{equation}
Applying Corollary~\ref{cor:eDindep} with $i=0$, we may verify the
inclusion and vanishing in~\eqref{eqn:tDrestate} separately on each
direct summand in the decomposition as a $\wth$-semisimple
$\C[E_i]$-module. Proposition~\ref{prop:string} gives the necessary
and sufficient condition for these properties, and this condition
is preserved under twisting by $\C_{\La_0}$. It follows that
\[
{}^{\theta} L ( \La )^u \otimes \C_{\La_0} \subset 
\mathscr{D}_0 \bigl( {}^{\theta} L ( \La )^u \otimes \C_{\La_0} \bigr) 
\equiv 
\mathbb{L}^{\bullet}\mathscr{D}_0 \bigl( {}^{\theta} L ( \La )^u \otimes \C_{\La_0} \bigr),
\]
and
\[
\mathbb{L}^{\bullet}\mathscr{D}_0 \bigl( {}^{\theta} L ( \La )^u \otimes \C_{\La_0} \bigr) 
\;\cong\; \varprojlim_r \mathbb{L}^{\bullet}\mathscr{D}_0 \bigl( M(r) \bigr).
\]
Thus there is a surjection
\[
\hom_\fB \bigl( {}^{\theta} L ( \La )^u \otimes \C_{\La_0}, ( D ^{(k-1)}_{\overline{\La'}} )^{*} \bigr) 
\longrightarrow 
\hom_\fB \bigl( {}^{\theta} L ( \La )^u \otimes \C_{\La_0}, (  \bq D ^{(k-1)}_{\overline{s_0 (\La'-\La_0)}} )^{*} \bigr)
\]
by Theorem~\ref{thm:indsurj} and Theorem~\ref{thm:adj}. 
Hence the middle term $\ext^{i+1}({}^\theta L(\Lambda)^u,(\C_{\Lambda'})^*)$ in~\eqref{eqn:exLDs} also vanishes for $i = 0$.

Combining the above case-by-case analysis of $u \La \notin W \La$ and $u \La \in W \La$ 
(with the subcases $\La' - \La_0 \in P_{\af,(k-1)}^+$ and $\La' - \La_0 \notin P_{\af,(k-1)}^+$), 
we obtain
\[
\ext^{>0}_\fB \bigl( {}^{\theta} L ( \La )^u, ( \C_{\La'} )^{*} \bigr) = 0 
\qquad \text{for all } u \in W_\af.
\]

Since ${}^{\theta} L ( \La )^u$ is a cyclic $\tb$-module whose head is $\C_{-u \La}$, 
the corresponding $\hom$ space vanishes unless $u \La = \La'$. 
We therefore conclude that
\begin{equation}
\ext^{i}_\fB \bigl( {}^{\theta} L ( \La )^u, ( \C_{\La'} )^{*} \bigr) 
\;\cong\; \C^{\delta_{i0} \delta_{u\La, \La'}}
\label{eqn:Dworth}.
\end{equation}

Summarizing the above, we obtain
\begin{equation}
\ext^{\bullet}_\fB ( {}^{\theta} L ( \La )^w, ( D_\mu ^{(k)} )^{*} ) \cong 
\begin{cases} 
\ext^{\bullet}_\fB ( {}^{\theta} L ( \La ), ( \C _{\La'} )^{*} ) & \text{if } w \le v, \\
0 & \text{otherwise}
\end{cases}.
\label{exttop2}
\end{equation}
\medskip

We have
\begin{equation}
\La' = x \cdot \La = x ( \La + \rho ) - \rho \qquad \text{for } x \in W_\af\label{eqn:Wdot}
\end{equation}
if and only if $x = e$ and $\La = \La'$, since $\langle \al_i^{\vee}, \rho \rangle > 0$ for each $i \in \tI_\af$. We have $\rho(K) = h^{\vee} > 0$ (the dual Coxeter number),
\[
P_{\af, k}^+ + h^{\vee} \La_0 \subset P^+_{\af, (k + h^{\vee})},
\]
and that $P^+_{\af,(k + h^{\vee})}$ is the fundamental domain for the action of $W_\af$ on $P_{\af,(k + h^{\vee})}$. Therefore, by the BGG resolution~\eqref{eqn:BGG-res} applied to the first factor on the right-hand side of~\eqref{exttop2}, we obtain
\begin{equation}
\ext^{i}_\fB ( {}^{\theta} L ( \La )^w, ( D_\mu ^{(k)} )^{*} ) \cong 
\begin{cases} 
\C & \text{if } i = 0,\ w \le v,\ \La = \La', \\
0 & \text{otherwise}
\end{cases}.
\label{eqn:LDdual}
\end{equation}
	
Since the maps in~\eqref{eqn:LDdual} are induced from the right-hand side of~\eqref{exttop2} through successive applications of the Demazure functors, the inclusion
\[
{}^{\theta} L ( \La )^w \subset {}^{\theta} L ( \La )^{w'} \qquad (w' \in W_\af)
\]
induces an isomorphism
\begin{equation}
\hom_\fB \bigl( {}^{\theta} L ( \La )^w, ( D_\mu ^{(k)} )^{*} \bigr) \cong 
\hom_\fB \bigl( {}^{\theta} L ( \La )^{w'}, ( D_\mu ^{(k)} )^{*} \bigr)\label{eqn:homthick}
\end{equation}
whenever $w' \le w \le v$. Since $\{ {}^{\theta} L ( \La )^u \}_{u \in W_\af}$ gives rise to a $\bD^{(k)}$-filtration of ${}^{\theta} L ( \La )$ described in Corollary~\ref{cor:LtoDfilt}, taking the associated graded of this filtration yields
\begin{equation}
\ext^{i}_\fB ( \bD _\la^{(k)}, ( D_\mu ^{(k)} )^{*}) \cong 
\begin{cases} 
\C & \text{if } i = 0,\ \la + \mu = 0, \\
0 & \text{otherwise}
\end{cases},
\end{equation}
where $\shom$-part follows from~\eqref{eqn:homthick}, and $\ext^{>0}$-part follows from~\eqref{eqn:LDdual}. This completes the proof of Theorem~\ref{thm:dual}(1) for level $k$.

\medskip

We now prove Theorem~\ref{thm:dual}(2) for level $k$. By Proposition~\ref{prop:DCext}, we have
\[
\ext^{\bullet}_{\fC} \bigl( \bW_\la^{(k)}, ( W_{\mu}^{(k)} )^* \bigr) 
= 
\ext^{\bullet}_{\fB} \bigl( \bW_\la^{(k)}, ( W_{\mu}^{(k)} )^* \bigr)
\qquad \text{for } \la, \mu \in P^+.
\]

Note that $W_{\mu}^{(k)} = D_{w_0 \mu}^{(k)}$, and that $\bW_\la^{(k)}$ admits a filtration by $\bD^{(k)}_{w \la}$ with multiplicity one, by Proposition~\ref{prop:WDfilt}. Hence, by the $\ext^{>0}$-vanishing part of the first assertion, we obtain
\[
\ext^{\bullet}_{\fB} \bigl( \bW_\la^{(k)}, ( W_{\mu}^{(k)} )^* \bigr) 
\cong 
\bigoplus_{\nu \in W\la} 
\hom_{\fB} \bigl( \bD_{\nu}^{(k)}, ( D_{w_0 \mu}^{(k)} )^* \bigr).
\]

This further implies
\[
\ext^{i}_{\fB} \bigl( \bW_\la^{(k)}, ( W_{\mu}^{(k)} )^* \bigr) 
\cong 
\begin{cases} 
\C & \text{if } \la + \mu_-=0,\ i = 0, \\
0 & \text{otherwise}
\end{cases}
\]
by the $\hom$-part of the first assertion. This completes the proof of Theorem~\ref{thm:dual}(2) for level $k$.
\end{proof}

\begin{thm}\label{thm:crit-ind}
Theorems~\ref{thm:dual} for level $\le k$ and Theorem~\ref{thm:critfilt} for level $<k$ imply Theorem~\ref{thm:critfilt} for level $k$.
\end{thm}

\begin{proof}
If the module $M$ admits one of the four types of filtrations, then we have a presentation
\[
M \cong \varprojlim_r M(r) \in \fB_{\bd}
\]
by a projective system $\{M(r)\}_r$ with surjective transition maps such that each term $M(r)$ admits the same type of filtration.

Thus, by Theorem~\ref{thm:dual} and Lemma~\ref{lem:comm}, the desired vanishing of higher $\ext$-groups follows by passing to the limit, since each finite stage $M(r)$ has the same type of filtration and hence has vanishing higher $\ext$ against the corresponding dual family.

In general, if $\mathrm{Ext}^{>0} ( M, N ) = 0$ for some modules $M$ and $N$, then in particular $\mathrm{Ext}^1 ( M, N ) = 0$. Therefore, it suffices to show that the vanishing of $\mathrm{Ext}^1$ with respect to one of the four series of modules implies the existence of the corresponding filtration.
	
Since the cases of $\bW^{(k)}$ and $W^{(k)}$ are analogous to those of $\bD^{(k)}$ and $D^{(k)}$ if we replace Theorem~\ref{thm:dual}(1) with Theorem~\ref{thm:dual}(2) and Proposition~\ref{prop:Dhw} with Corollary~\ref{cor:Whw}, we focus on the cases of $\bD^{(k)}$ and $D^{(k)}$, and leave the details of the remaining cases to the reader.

\medskip

We prove Theorem~\ref{thm:critfilt}(1) for level $k$. Suppose $M \in \fB_\bd$ satisfies
\[
\ext^1_{\fB} ( \bD^{(k)}_\la, M^{\vee} ) = 0 \qquad \text{for all } \la \in P.
\]
Consider the set
\[
H := \{ ( \la, s ) \in P \times \Z \mid \Hom_{\fB} ( \bq^s \bD^{(k)}_\la, M^{\vee} ) \neq 0 \},
\]
and choose an element $( \mu, m ) \in H$ such that
\[
- \la + k \La_0 - s \delta \not\blacktriangleleft - \mu + k \La_0 - m \delta \qquad ( \la, s ) \in H \setminus \{ (\mu, m ) \}.
\]
Since the degrees of $M$ are bounded from below and $(\bullet,\bullet)_0$ is positive definite, the $\mathbf q$-values of the weights of $M$, which is of level $k>0$, are bounded from below. Moreover, for every lower bound on the degree and every upper bound on the $\mathbf q$-value, there are only finitely many elements of $P\times\mathbb Z$ satisfying both bounds. Hence, such a pair $(\mu,m)$ necessarily exists.

\begin{claim}\label{clm:Mwt}
The $\wth$-weight $-\mu + k \La_0 - m \delta$ appears in $\hd M$.
\end{claim}

\begin{proof}
Note that every $\wth$-weight of $\soc M^{\vee}$ gives rise to an element of $H$ by $\hd \bD^{(k)}_\la = \C_{\la - k \La_0}$. We have a map
\[
f : M \longrightarrow \bq^{-m} ( \bD^{(k)}_\mu )^{\vee}.
\]
Every $\wth$-weight of $\bq^{-m} ( \bD^{(k)}_\mu )^{\vee}$ satisfies $\blacktriangleleft - \mu+ k \La_0 - m \delta$ by Proposition~\ref{prop:Dhw}(2) (proved for level $k$ in Proposition~\ref{prop:Dhw-ind}). As the $\wth$-weight $- \mu + k \La_0 - m \delta$-part of $\bq^{-m} ( \bD^{(k)}_\mu )^{\vee}$ is one-dimensional, we conclude
\[
\Hom_{\wth} ( \bq^{-s} \C_{-\la + k \La_0}, \bq^{-m} ( \ker ( \bD^{(k)}_\mu \to \C_{\mu - k \La_0} ) )^{\vee} ) = 0 \qquad (\la, s) \in H.
\]
It follows that
\[
f ( \hd M ) \subset \bq^{-m} \C_{- \mu + k \La_0} \subset \bq^{-m} ( \bD^{(k)}_\mu )^{\vee}
\]
as a $\wth$-module, which implies $\mathrm{Im} \, f = \bq^{-m} \C_{- \mu + k \La_0}$. This proves the claim.
\end{proof}

We return to the proof of Theorem~\ref{thm:critfilt}(1) for level $k$. We have a projective cover $f_1 \colon Q \twoheadrightarrow M$ in $\fB_{\bd}$ whose $d$-grading is bounded from below. Every projective module in $\fB_\bd$ is a direct sum of projective covers
of one-dimensional modules, as described in~\S\ref{subsubsec:Bbasic}. By Claim~\ref{clm:Mwt}, the projective module $Q$ has a direct summand $\bq^{-m} Q_{- \mu + k \La_0}$. Since $D_{- \mu}^{(k)}$ is generated by its $\wth$-weight $- \mu + k \La_0$-part, we also have a surjective map
\[
f_2 \colon Q \longrightarrow \bq^{-m} D_{-\mu}^{(k)},
\]
by extending the quotient map $\bq^{-m} Q_{- \mu + k \La_0} \twoheadrightarrow \bq^{-m} D_{-\mu}^{(k)}$ by zero on a complementary direct summand of $Q$. We then define the maximal quotient
\[
M' := Q \big/ \left( \ker f_1 + \ker f_2 \right),
\]
which is also the maximal quotient of $\bq^{-m} D_{-\mu}^{(k)}$ that admits a (compatible) surjection from $M$. Set $M'' := \ker ( M \to M' )$.

\smallskip

Assume, for contradiction, that $M'$ is a proper quotient of $\bq^{-m} D_{-\mu}^{(k)}$. Then, we have some $- \gamma + k \La_0 - c \delta \in P_\af$ such that
\begin{equation}
- \gamma + k \La_0 - c \delta \neq - \mu + k \La_0 - m \delta, \quad
- \gamma + k \La_0 - c \delta \blacktriangleleft - \mu + k \La_0 - m \delta, \label{eqn:DMgmcomp}
\end{equation}
and
\[
0 \rightarrow \C_{- \gamma + k \La_0 - c \delta} \rightarrow \widetilde{M}' \rightarrow M' \rightarrow 0
\]
is a non-split short exact sequence in $\fB_{\bd}$ arising from a larger quotient $\widetilde{M}'$ of $\bq^{-m} D_{-\mu}^{(k)}$. Pulling back this short exact sequence along the quotient map
\[
M\twoheadrightarrow M',
\]
we obtain a short exact sequence
\begin{equation}
0 \rightarrow \C_{- \gamma + k \La_0 - c \delta} \rightarrow \widetilde{M} \rightarrow M \rightarrow 0.\label{eqn:DMext}
\end{equation}
Here $\widetilde{M}'$ cannot be a quotient of $M$ by the maximality of $M'$ (as a common quotient of $M$ and $\bq^{-m} Q_{-\mu + k \La_0}$ through $f_1,f_2$). It follows that~\eqref{eqn:DMext} defines a non-split short exact sequence. Consequently, we have
\begin{equation}
0 \neq \Ext^1_{\fB} ( M, \C_{- \gamma + k \La_0 - c \delta} ) \cong \Ext^1_{\fB} ( \C_ {\gamma - k \La_0 + c \delta}, M^{\vee} ).\label{eqn:proh-ex}
\end{equation}

We set $K := \ker ( \bD_{\gamma}^{(k)} \to \C_{\gamma - k \La_0})$. Assume that
\[
\Hom_{\fB} ( \bq^c K, ( M )^{\vee} ) \neq 0.
\]
Then, $\soc ( M )^{\vee}$ contains $\C_{\kappa - k \La_0 + e \delta}$ as its direct summand, and we have
\[
- \kappa + k \La_0 - e \delta \blacktriangleleft - \gamma + k \La_0 - c \delta \qquad \bigl(  \blacktriangleleft - \mu + k \La_0 - m \delta \bigr)
\]
by Proposition~\ref{prop:Dhw}(2) (and~\eqref{eqn:DMgmcomp}). Since $\hd \bq^e \bD_{\kappa}^{(k)} = \C_{\kappa - k \La_0 + e \delta}$, we have
\[
\Hom_{\fB} ( \bq^e \bD_{\kappa}^{(k)}, ( M )^{\vee} ) \neq 0.
\]
This contradicts the minimality of $( \mu, m ) \in H$. Therefore, we necessarily have
\begin{equation}
\Hom_{\fB} ( \bq^c K, ( M )^{\vee} ) = 0. \label{eqn:Kvanish}
\end{equation}

Thus, the short exact sequence associated to
\[
0 \rightarrow \bq^c K \rightarrow \bq^c \bD_{\gamma}^{(k)} \rightarrow \C_{\gamma - k \La_0 + c \delta} \rightarrow 0
\]
gives an exact sequence
\[
0 = \Hom_{\fB} ( \bq^c K, ( M )^{\vee} ) \rightarrow \Ext^1_{\fB} ( \C_{\gamma - k \La_0 + c \delta}, ( M )^{\vee} ) \rightarrow \Ext^1_{\fB} ( \bq^c \bD_{\gamma}^{(k)}, ( M )^{\vee} ).
\]
Combined with~\eqref{eqn:proh-ex}, we derive
\[
 \Ext^1_{\fB} ( \bq^c \bD_{\gamma}^{(k)}, ( M )^{\vee} ) \neq 0,
\]
which contradicts the assumption.

Therefore, $M'$ cannot be a proper quotient of $\bq^{-m} D_{-\mu}^{(k)}$.

\smallskip

As a consequence, we have $M' = \bq^{-m} D_{-\mu}^{(k)}$. Consider (a portion of) the long exact sequence
\[
\ext^1_{\fB} ( \bD^{(k)}_{\la}, M^{\vee}) 
\longrightarrow 
\ext^1_{\fB} ( \bD^{(k)}_{\la}, (M'')^{\vee}) 
\longrightarrow 
\ext^2_{\fB} ( \bD^{(k)}_{\la}, (D_{-\mu}^{(k)})^{\vee})
\]
for every $\la \in P$. Since we have $\ext^1_{\fB} ( \bD^{(k)}_{\la}, M^{\vee}) = 0$ and $\ext^2_{\fB} ( \bD^{(k)}_{\la}, (D_{-\mu}^{(k)})^{\vee}) = 0$, it follows that
\[
\ext^1_{\fB} ( \bD^{(k)}_{\la}, (M'')^{\vee}) = 0 \qquad \text{for all } \la \in P.
\]

Thus the problem reduces to the case of $M''$.
By repeating this procedure, we obtain a decreasing filtration
\[
M = M^0 \supset M^1 \supset M^2 \supset \cdots
\]
such that each quotient $M^i / M^{i+1}$ is isomorphic to a grading shift of some $D^{(k)}_{\gamma}$ with $\gamma \in P$.

Let $d_0 \in \Z$ be the minimal degree such that $M_{d_0} \neq 0$. We set
\[
S_0 := \{ \gamma \in P \mid \Hom_{\h} (\C_{\gamma}, M_{d_0} ) \neq 0\}. 
\]
Then, $|S_0| < \infty$ as we have $\dim M_{d_0} < \infty$ by $M \in \fB_\bd$. We set
\[
S_1 := \{ ( \kappa, l ) \in P \times \Z_{\ge d_0}  \mid \exists \gamma \in S_0 \text{ such that } \mathbf{q} ( \kappa + k \La_0 + l \delta ) \le \mathbf{q} ( \gamma + k \La_0 + d_0 \delta ) \}.
\]
Since $S_0$ is finite, the finiteness property of $\mathbf q$ discussed above implies $|S_1|<\infty$.

While $( M^i )_{d_0} \neq 0$, the head of $M^i$ contains a $\wth$-weight of the
form $\gamma + k \La_0 + d_0 \delta$ with $\gamma \in S_0$, and the
$\blacktriangleleft$-minimality in the choice of $( \mu_i, m_i )$ with $M^i/M^{i+1} \cong \bq^{-m_i} D_{-\mu_i}^{(k)}$ forces
$\mathbf{q} ( - \mu_i + k \La_0 - m_i \delta ) \le
\mathbf{q} ( \gamma + k \La_0 + d_0 \delta )$, i.e.
$( - \mu_i, - m_i ) \in S_1$. Since the
$( - \mu_i + k \La_0 - m_i \delta )$-weight space of $M^i$ strictly decreases
when we pass to $M^{i+1}$, this can happen only finitely many times.

Consequently, we have $(M^i)_{d_0} = 0$ for $i > n_0$ for some $n_0 \in \Z_{> 0}$.
It follows that for each $m \in \Z$ there exists $i_m \gg 0$ such that 
$(M^i)_{j} = 0$ for all $j < m$ whenever $i > i_m$. 
This implies that the filtration $\{M^i\}_i$ is separable. 
Therefore,
\[
M \cong \varprojlim_i M/M^i \;\in\; \fB_{\bd}.
\]
Hence Theorem~\ref{thm:critfilt}(1) holds for level $k$.
\medskip
	
We consider Theorem~\ref{thm:critfilt}(2) for level $k$. 
By~\eqref{eqn:transpose} and the Yoneda interpretation (\S\ref{subsubsec:Yoneda}), 
our condition on $M$ is equivalent to
\[
\ext^1_{\fB} \bigl( M, ( D_{\la}^{(k)} )^{*} \bigr) 
\cong \ext^1_{\fB} \bigl( D_{\la}^{(k)}, M^{\vee} \bigr) = 0 
\qquad (\la \in P),
\]
while the analogous weight estimate to Claim~\ref{clm:Mwt} requires Proposition~\ref{prop:Dhw}(1), which is unconditionally proved in Proposition~\ref{prop:Dhw-ind} (and hence here we do not need to assume Theorem~\ref{thm:dual}(1) for level $k$, though we use it in other places).

Thus, arguments parallel to those in the proof of Theorem~\ref{thm:critfilt}(1) for level $k$ 
imply Theorem~\ref{thm:critfilt}(2) for level $k$.

This completes the proof.
\end{proof}

\section{Branching rules}\label{sec:branch}
We continue to work in the setting of the previous section.

\begin{thm}\label{thm:branch-b}
Let $k \in \Z_{>0}$ and $\la \in P$. Then:
	\begin{enumerate}
	\item The graded $\tb$-module $\bD^{(k+1)}_\la \otimes \C_{\La_0}$ admits a $\bD^{(k)}$-filtration;
	\item The graded $\tb$-module $D^{(k)}_\la \otimes \C_{\La_0}$ admits a $D^{(k+1)}$-filtration.
	\end{enumerate}
Moreover, for every $\mu \in P$, we have
\[
( \bD^{(k+1)}_\la : \bD^{(k)}_\mu )_q = ( D^{(k)}_\mu : D^{(k+1)}_\la )_q.
\]
\end{thm}

\begin{rem}
Theorem~\ref{thm:branch-b}(2) was established by Joseph~\cite{Jos06} 
(see also Naoi~\cite[Remark~4.15]{Nao12}) in the case where $\g$ is of type $\mathsf{ADE}$.
\end{rem}

\begin{proof}[Proof of Theorem~\ref{thm:branch-b}]
In view of Theorem~\ref{thm:critfilt}, the first two assertions are equivalent to the vanishing
\begin{equation}
\ext^{1}_{\fB} ( \bD^{(k+1)}_\la, ( D_{\mu}^{(k)} \otimes \C_{\La_0} )^{*} ) = 0 \qquad \text{for all } \la, \mu \in P, \label{eqn:extdual}
\end{equation}
which is exactly the vanishing established in~\eqref{eqn:extdualone}.

The final assertion follows from the fact that
\[
\left\langle \bD^{(k+1)}_\la,\ D^{(k+1)}_\mu \right\rangle_{\mathtt{EP}} = \delta_{\la,\mu} = \left\langle \bD^{(k)}_\la,\ D^{(k)}_\mu \right\rangle_{\mathtt{EP}} \qquad \text{for all } \la, \mu \in P,
\]
by Theorem~\ref{thm:dual}. Hence, the transition matrices between their graded characters are transposes of each other by Lemma~\ref{lem:top}.
\end{proof}

\begin{cor}\label{cor:kumar}
Let $k \in \Z_{>0}$ and $\la \in P$. Then the module $D^{(k)}_\la$, regarded as a $\gb$-module, admits a filtration by Demazure modules of $\g$.
\end{cor}

\begin{proof}
By repeatedly applying Theorem~\ref{thm:branch-b}, we obtain a $D^{(l)}$-filtration of $D^{(k)}_{\la}$ for each $l \ge k$ such that $\langle \vartheta^{\vee}, \mu_+ \rangle < l$ for every $\mu \in \overline{\Psi ( D^{(k)}_{\la} )}$. 

Each filtration piece of the form $D^{(l)}_{\nu}$ (for some $\nu \in P$, up to grading shifts) is contained in $V_{\nu_+}$, since $\nu_+ + l \La_0 \in P^+_{\af,l}$. In particular, each such $D^{(l)}_{\nu}$ is a Demazure module of $\g$ when viewed as a $\gb$-module. 

This proves the claim.
\end{proof}

\begin{lem}\label{lem:Pfilt}
Let $k \in \Z_{>0}$. For each $\la \in P^+$, the module $P_{\la - k \La_0}$ admits a $\bW^{(k)}$-filtration.	
\end{lem}

\begin{proof}
Since $P_{\la-k\Lambda_0}$ is projective in $\fC$, we have
\[
\ext^1_\fC \bigl( P_{\la-k\Lambda_0}, (W^{(k)}_\nu)^* \bigr)=0
\]
for every $\nu\in P^+$. Hence, by Theorem~\ref{thm:critfilt}, the module $P_{\la - k \La_0}$ admits a $\bW^{(k)}$-filtration.
\end{proof}

\begin{thm}\label{thm:branch}
Let $k \in \Z_{>0}$ and $\la \in P^+$. Then:
	\begin{enumerate}
	\item The graded $\tg_{\ge 0}$-module $\bW^{(k+1)}_\la \otimes \C_{\La_0}$ admits a $\bW^{(k)}$-filtration;
	\item The graded $\tg_{\ge 0}$-module $W^{(k)}_\la \otimes \C_{\La_0}$ admits a $W^{(k+1)}$-filtration.
	\end{enumerate}
Moreover, for each $\mu \in P^+$, we have
\begin{equation}
( P_\la : \bW^{(k)}_\mu )_q = [ W^{(k)}_\mu : V_\la ]_q
\quad \text{and} \quad
( \bW^{(k+1)}_\la : \bW^{(k)}_\mu )_q = ( W^{(k)}_\mu : W^{(k+1)}_\la )_q. \label{eqn:hBGG-rec}
\end{equation}
\end{thm}

\begin{proof}
In view of Theorem~\ref{thm:critfilt}, the first two assertions are equivalent to
\begin{equation}
\ext^{>0}_{\fC} \big( \bW^{(k+1)}_\la, ( W_{\mu}^{(k)} \otimes \C_{\La_0} )^{*} \big) = 0 
\qquad \text{for all } \la, \mu \in P^+. \label{extWWone}
\end{equation}

Since we have 
\[
W_{\mu}^{(k)} = D_{w_0 \mu}^{(k)} = \mathbb L^{\bullet} \mathscr D_{w_0} \big( D_{\mu}^{(k)} \big)
\quad \text{and} \quad 
\bW_{\la}^{(k)} = \mathbb L^{\bullet} \mathscr D_{w_0} \big( \bD_{\la}^{(k)} \big),
\]
Theorem~\ref{thm:adj} reduces \eqref{extWWone} to a special case of \eqref{eqn:extdual}.

The final assertion follows from the equalities
\begin{align*}
\left\langle \bW^{(k+1)}_\la,\ W^{(k+1)}_\mu \right\rangle_{\mathtt{EP}} &= \delta_{\la,\mu} 
= \left\langle \bW^{(k)}_\la,\ W^{(k)}_\mu \right\rangle_{\mathtt{EP}}, \\
\left\langle P_\la,\ V_\mu \right\rangle_{\mathtt{EP}} &= \delta_{\la,\mu}
\qquad \text{for all } \la, \mu \in P^+,
\end{align*}
which follow from Theorem~\ref{thm:dual} and Proposition~\ref{prop:DCext}. Hence, the transition matrices between their graded characters are transposes of each other by Lemma~\ref{lem:top}.
\end{proof}

\section{Applications}\label{sec:app}

We retain the setting of the previous section.

\begin{lem}\label{lem:idDL}
Let $k \in \Z_{>0}$ and $\la \in P^+$. If $\langle \vartheta^{\vee}, \la \rangle < k$, then we have an isomorphism
\[
\bW^{(k)}_{\la} \otimes \C_{\La_0} \cong {}^{\theta} L ( - \la_- + (k-1) \La_0 )
\]
as $\tb$-modules.
\end{lem}

\begin{proof}
Consider the surjection
\begin{equation}\label{quotW1st}
\pi : {}^{\theta} L ( - \la_- + k \La_0 ) \longrightarrow \bW^{(k)}_\la,
\end{equation}
which exists by the definition of $\bW^{(k)}_\la$~\eqref{eqn:defW}. A further examination of~\eqref{eqn:defW} shows that its kernel is the sum of thick Demazure modules ${}^{\theta} L ( - \la_- + k \La_0 )^u \subset {}^{\theta} L ( - \la_- + k \La_0 )$, where $u$ runs over elements in $\bigcup_{\beta \neq 0} W t_\beta W$. In particular, this kernel contains the submodule generated by the $\theta$-twist of the vector $\bv_{s_0(-\la_- + k \La_0)}$ as a $\tb$-module.

The defining relations of ${}^{\theta} L ( - \la_- + (k-1) \La_0 )$ are given by the kernel of its projective cover in $\fB$. More precisely,
\[
{}^{\theta} M ( - \la_- + (k-1) \La_0 ) \twoheadrightarrow {}^{\theta} L ( - \la_- + (k-1) \La_0 )
\]
is the projective cover in $\fB$, and hence its kernel describes all defining relations of
${}^{\theta} L ( - \la_- + (k-1) \La_0 )$. By the first step of the BGG resolution~\eqref{eqn:BGG-res}, this kernel is generated by vectors of $\wth$-weights
\[
- s_i ( - \la_- + (k-1) \La_0 ) + \al_i \qquad (i \in \tI_\af).
\]

For $i \neq 0$, they are $d$-degree zero and are outside of
\[
\Psi ( V_{\la - (k-1) \La_0} ) = \Psi ( {}^{\theta} L ( - \la_- + (k-1) \La_0 )_0 ),
\]
and hence they give rise to a relation in the degree zero part of $\bW^{(k)}_{\la}$. 

For $i=0$, the corresponding weight before twisting by $\theta$ is the extremal weight
$s_0(-\la_-+k\Lambda_0)$ of $L(-\la_-+k\Lambda_0)$, and hence generates the thick Demazure submodule
${}^{\theta}L(-\la_-+k\Lambda_0)^{s_0}\otimes\C_{\Lambda_0}$. We have $s_0 ( - \la_- + k \La_0 ) \neq - \la_- + k \La_0$, since
\[
- \langle \vartheta^{\vee}, \la_- \rangle = \langle \vartheta^{\vee}, \la \rangle < k.
\]
Hence, we have ${}^{\theta} L ( - \la_- + k \La_0 )^{s_0} \subsetneq {}^{\theta} L ( - \la_- + k \La_0 )$, and the kernel of the quotient map ${}^{\theta} L ( - \la_- + k \La_0 ) \twoheadrightarrow {}^{\theta} L ( - \la_- + (k-1) \La_0 )$ is ${}^\theta L ( - \la_- + k \La_0 )^{s_0}$.

Therefore, the cyclic vector of ${}^\theta L ( - \la_- + k \La_0 )^{s_0}$ as $\tb$-module, which necessarily has $\wth$-weight $s_0 ( - \la_- + k \La_0 )$, gives rise to zero in $\bW^{(k)}_{\la}$ when it is realized as a quotient of ${}^{\theta} L ( - \la_- + k \La_0 ) $. This yields a surjection
\begin{equation}\label{eqn:surjcomp}
{}^{\theta} L ( - \la_- + (k-1) \La_0 ) \longrightarrow \bW^{(k)}_\la \otimes \C_{\La_0}.
\end{equation}

Note that $s_0 \not\le w \in W_\af$ if and only if $w \in W$. By Corollary~\ref{cor:dominant}, any $\g$-stable proper thick Demazure submodule $L ( \La )^w$ of $L ( \La )$ must satisfy $w \notin W$ and hence $w \ge s_0$. Thus, every proper $\g$-stable thick Demazure submodule of
${}^\theta L ( - \la_- + k \La_0 )$
is contained in
${}^\theta L ( - \la_- + k \La_0 )^{s_0}$.
Since the kernel of~\eqref{quotW1st} is the sum of such submodules, we obtain
\[
\ker\pi \subset {}^\theta L ( - \la_- + k \La_0 )^{s_0}.
\]
Together with the opposite inclusion proved above, this yields
\[
\ker\pi = {}^\theta L ( - \la_- + k \La_0 )^{s_0}.
\]
Hence the surjection~\eqref{eqn:surjcomp} is an isomorphism.
\end{proof}

\begin{rem}\label{rem:theta}
The inner product $\langle\bullet,\bullet \rangle_{\mathtt{EP}}$,
defined by the Euler--Poincar\'e characteristic in terms of
$\ext^{\bullet}_{\fB}$, coincides, up to a constant factor, with the
$t=0$ specialization of the Cherednik inner product; see~\cite[Appendix]{FKM}
for the non-symmetric setting. In the symmetric case, the corresponding
identification was established in~\cite{CI15}. Under this inner product, the
$t=0$ specializations of Macdonald polynomials are mutually orthogonal up to
normalization factors.

By the Frenkel--Kac construction~\cite{FK80,Seg81} and
Lemma~\ref{lem:idDL}, the graded character $\gch \bW_{\la}^{(2)}$ is
essentially a theta function when $\langle \vartheta^{\vee}, \la \rangle \le 1$.
Hence, the case $k = 2$ of Theorem~\ref{thm:dual}(2) asserts that theta
functions arise as duals of Demazure characters with respect to the Hall inner
product (\cite[\S5.1]{Mac03}), which is the Euler--Poincar\'e characteristic defined in terms of
$\ext^{\bullet}_{\fC}$.

In the same vein, the family $\{\gch \bW^{(k)}_\la \}_{\la \in P^+}$ for
$k \ge 2$ may be viewed as a natural enhancement of certain spaces of modular
forms, as discussed in~\cite[\S13]{Kac}.
\end{rem}

\begin{rem}\label{rem:limit}
The proofs of Lemma~\ref{lem:idDL} and Corollary~\ref{cor:kumar} imply the identities
\[
V_\la = W^{(k)}_\la \quad \text{for} \quad k \gg 0, \qquad
P_\la = \varprojlim_{k} \bW_\la^{(k)} \quad \text{for each} \quad \la \in P^+.
\]
Accordingly, we can interpret asymptotic behaviors of~\eqref{eqn:()q} as
\[
\lim_{k \to \infty} (\bullet : W^{(k)}_{\mu})_q = [\bullet : V_\mu]_q,
\qquad
\lim_{k \to \infty} (\bW^{(k)}_\la : \bullet )_q = (P_\la : \bullet )_q,
\]
whenever these expressions are well-defined.
\end{rem}

We define the level-$k$-restricted Kostka polynomials as
\begin{equation}
X^{(k)}_{\mu, \La}(q) := ( {}^{\theta} L ( \La ) : \bW^{(1)}_{-\mu_-} )_q \in \Z_{\ge 0} [q] \qquad \mu \in P^+, \La \in P_{\af, k}^+.\label{eqn:Xkdef}
\end{equation}
In case $\langle d, \La \rangle = 0$, this polynomial coincides with the level-$k$-restricted Kostka polynomial defined by~\cite{Oka,HKOTY,HKKOTY,HKOTT,SS01} in terms of affine crystals~\cite{KKMMNN,KMS95} by~\cite[Corollaries~5.12 and~3.6]{KL17}. Here we further reinterpret these polynomials using our results:

\begin{cor}\label{cor:kostka}
Assume that $\g$ is of type $\mathsf{ADE}$. Let $k \in \Z_{>0}$ and let $\La \in P_{\af, k}^+$ satisfy $\langle d, \La \rangle = 0$. For each $\mu \in P^+$, we have
\begin{equation}\label{eqn:Xk}
X^{(k)}_{\mu, \La} ( q ) = ( W^{(1)}_{\mu} : W^{(k+1)}_{\overline{\La}} )_q \in \Z_{\ge 0} [q].
\end{equation}
Moreover, the right-hand side of~\eqref{eqn:Xk} makes sense for arbitrary types of $\g$, and always has nonnegative coefficients.
\end{cor}

\begin{cor}\label{cor:KR}
Let $k \in \Z_{>0}$. Let $M \in \fC_\bd$ be a finite-dimensional $W^{(k)}$-filtered module on which $K$ acts by $k$, and assume that
\[
\bigcup_{w \in W_\af} \mathscr D_{w} ( M ) \cong \bigoplus_{\La \in P_\af^+} L ( \La )^{\oplus m_\La ( M )} \qquad \text{with}\quad m_{\La} ( M ) \in \Z_{\ge 0}.
\]
Then, for each $\la \in P^+$ with $\langle \vartheta^{\vee}, \la \rangle \le k$, we have
\begin{equation}\label{eqn:KRchar}
\sum_{\overline{\La} = \la} m_{\La} ( M )\, q^{\langle d, \La \rangle} = ( M : W_\la^{(k+1)} )_q.
\end{equation}
\end{cor}

\begin{rem}\label{rem:kostka}
\textbf{(1)} Since the right-hand side of~\eqref{eqn:Xk} is well-defined for every $\La \in P_\af$ with $\overline{\La} \in P^+$, Corollary~\ref{cor:kostka} embeds the set of level-$k$-restricted Kostka polynomials into the following family of polynomials indexed by $\la, \mu \in P^+$:
\[
( \bW^{(k+1)}_{\la} : \bW^{(1)}_{\mu} )_q 
\stackrel{(\ref{eqn:hBGG-rec})}{=} 
( W^{(1)}_{\mu} : W^{(k+1)}_{\la} )_q 
\in \Z_{\ge 0}[q].
\]

\textbf{(2)} Suppose $\g$ is not of type $\mathsf{E}_7$, $\mathsf{E}_8$, or $\mathsf{F}_4$. Then, the graded character $\gch W^{(l)}_{l \la}$ ($l \le k$, $\la \in P^+$) coincides with the character of the tensor product $B$ of Kirillov–Reshetikhin crystals of level $l$ (see \cite[Theorem 5.1]{ST12} and \cite{FSS07,Nao18,BS20,NS21}). It follows that the left-hand side of~\eqref{eqn:KRchar} equals the character of the set of level-$k$-restricted highest-weight elements in $B$ (cf.~\cite[(3.9)]{HKOTT} or \cite[Definition 5.5]{KL17}).

More generally, for certain tensor products of Kirillov–Reshetikhin crystals (together with a construction of the corresponding module $M$), the left-hand side of~\eqref{eqn:KRchar} again represents the character of the set of level-$k$-restricted highest-weight elements~\cite{Nao13}. At least in these cases, identity~\eqref{eqn:KRchar} provides a module-theoretic interpretation of level-$k$-restricted generalized Kostka polynomials, thereby extending Corollary~\ref{cor:kostka} (see also~\cite{Oka,HKOTY,HKKOTY,HKOTT}).
\end{rem}

\begin{proof}[Proof of Corollary~\ref{cor:kostka}]
In view of Lemma~\ref{lem:idDL}, the right-hand side of~\eqref{eqn:Xkdef} satisfies
\[
( {}^{\theta} L ( \La ) : \bW^{(1)}_{-\mu_-} )_q = ( \bW^{(k+1)}_{- (\overline{\La})_-} : \bW^{(1)}_{-\mu_-} )_q.
\]
We now apply the affine Dynkin diagram automorphism inducing the finite Dynkin diagram automorphism $i \mapsto i'$ characterized by $\alpha_{i'}=-w_0\alpha_i$ to deduce 
\[
( \bW^{(k+1)}_{- (\overline{\La})_-} : \bW^{(1)}_{-\mu_-} )_q 
= ( \bW^{(k+1)}_{\overline{\La}} : \bW^{(1)}_{\mu} )_q.
\]
Applying Theorem~\ref{thm:branch}, this yields
\[
( \bW^{(k+1)}_{- (\overline{\La})_-} : \bW^{(1)}_{-\mu_-} )_q 
= ( \bW^{(k+1)}_{\overline{\La}} : \bW^{(1)}_{\mu} )_q 
= ( W_\mu^{(1)} : W_{\overline{\La}}^{(k+1)} )_q \in \Z_{\ge 0}[q].
\]
The last equality holds for arbitrary types of $\g$, completing the proof.
\end{proof}

\begin{proof}[Proof of Corollary~\ref{cor:KR}]
By Theorem~\ref{thm:thinD}, we have $\mathbb L^{\bullet} \mathscr D_{w}(W) \cong \mathscr D_{w}(W)$ for each $w \in W_\af$.

Now we calculate as follows:
\begin{align*}
\ext_{\fC}^{i} \big( {}^{\theta} L ( \la + k \La_0 ), M^* \big) 
&\cong \ext_{\fC}^{i} \big( {}^{\theta} L ( \la + k \La_0 ), \mathscr D_w(M)^* \big) \\
&\cong \varprojlim_{w} \ext_{\fC}^{i} \big( {}^{\theta} L ( \la + k \La_0 ), \mathscr D_w(M)^* \big) \\
&\cong \bigoplus_{\La \in P_\af^+} \ext_{\fC}^{i} \big( {}^{\theta} L ( \la + k \La_0 ), {}^{\theta} L ( \La ) \big)^{\oplus m_\La(M)} \\
&\cong \bigoplus_{\La \in P_\af^+} \varprojlim_{w} \ext_{\fC}^{i} \big( {}^{\theta} L ( \la + k \La_0 ), \mathscr D_w(\C_{\La})^* \big)^{\oplus m_\La(M)} \\
&\cong \bigoplus_{\La \in P_\af^+} \ext_{\fC}^{i} \big( {}^{\theta} L ( \la + k \La_0 ), \C_{\La}^* \big)^{\oplus m_\La(M)} \\
&\cong
\begin{cases}
\displaystyle \bigoplus_{\overline{\La} = \la} \C_{-\langle d, \La \rangle}^{\oplus m_\La(M)} & (i = 0), \\
0 & (i > 0).
\end{cases}
\end{align*}
Here the first and fifth isomorphisms follow from Theorem~\ref{thm:adj} and $\mathscr D_i ( {}^{\theta} L ( \La ) ) \cong {}^{\theta} L ( \La)$ for each $\La \in P^+_\af$ and $i \in \tI_\af$ (Theorem~\ref{thm:thinD}). The second and fourth are obtained by passing to the limit over
$w\in W_\af$, while the third follows from the assumption. Finally, the sixth isomorphism follows from the BGG resolution~\eqref{eqn:BGG-res}.

On the other hand, by Lemma~\ref{lem:idDL} and Corollary~\ref{cor:multcrit}, we obtain:
\[
\gdim \shom_{\fC} \big( {}^{\theta} L ( \la + k \La_0 ), M^* \big)
= \gdim \shom_{\fC} \big( \bW_{-\la_-}^{(k+1)}, M^* \big)
= ( M : W^{(k+1)}_\la )_{q^{-1}}.
\]

Comparing this with the graded dimension computed above and then replacing $q^{-1}$ with $q$, we obtain the asserted identity (cf.~\eqref{eqn:qtoqinv}).
\end{proof}

\begin{cor}\label{cor:num}
Let $k \in \Z_{>0}$ and $\la, \mu \in P^+$. If $\langle \vartheta^{\vee}, \la \rangle \le k$, then
\[
(W^{(k)}_{\mu} : W^{(k+1)}_\la )_q =
\begin{cases}
q^{m} & \text{if } \exists\, w \in W_\af \text{ such that } \mu + k \La_0 - m \delta = w ( \la + k \La_0 ), \\
0 & \text{otherwise}.
\end{cases}
\]
In addition, the module $W^{(k+1)}_\la$ appears only in the socle of $W^{(k)}_{\mu}$ (with respect to the $W^{(k+1)}$-filtration).
\end{cor}

\begin{proof}
By Lemma~\ref{lem:idDL} and the definition of $\bW^{(k)}$-modules, we have
\begin{align*}
(\bW^{(k+1)}_\la : \bW^{(k)}_{\mu})_q
& = ( {}^{\theta} L (-\la_- + k \La_0 ) : \bW^{(k)}_{\mu} )_q\\
& = 
\begin{cases}
q^{m} & \text{if } \exists\, w \in W_\af \text{ such that } \mu + k \La_0 - m \delta = w ( \la + k \La_0 ), \\
0 & \text{otherwise}.
\end{cases}
\end{align*}

Applying Theorem~\ref{thm:branch} yields the claimed multiplicity.

Now observe that $\la + k\La_0 \in P^+_{\af,k}$ implies $W^{(k)}_\la = V_\la$. Since we have
\[
V_\la \subset W^{(k+1)}_\la \subset W^{(k)}_\la,
\]
it follows that
\[
V_\la = W^{(k+1)}_\la = W^{(k)}_\la \subset \bq^{-m} W^{(k)}_\mu \subset L(\la + k \La_0)
\]
as $\tb$-modules. Since the multiplicity is at most one by the first part, this occurrence must account for the entire contribution of $W^{(k+1)}_\la$ in the $W^{(k+1)}$-filtration of $W^{(k)}_\mu$ and therefore it occurs in the socle. This proves the assertion.
\end{proof}

\begin{rem}\label{rem:num}
Let $k \in \Z_{>0}$ and $\la \in P^+$ with $\langle \vartheta^{\vee}, \la \rangle < k$. Then we have a short exact sequence of $\tb$-modules:
\[
0 \rightarrow \C_\la \longrightarrow \mathscr D_0 ( \C_{\la + k \La_0} ) \otimes \C_{-k \La_0}
\longrightarrow \bq^{-1} \mathscr D_0 ( \C_{\la + \vartheta + (k+1) \La_0} ) \otimes \C_{-(k+1)\La_0} \rightarrow 0.
\]

Applying $\mathscr D_{w_0}$ to this sequence, we obtain another short exact sequence:
\[
0 \rightarrow \bq^{k - \langle \vartheta^{\vee}, \la \rangle } W^{(k+1)}_{\la} \longrightarrow W^{(k)}_\mu \longrightarrow W^{(k+1)}_\mu \otimes \C_{- \La_0} \rightarrow 0,
\]
where $\mu = \la - \langle \vartheta^{\vee}, \la \rangle \vartheta + k \vartheta$.

The multiplicity computation in Corollary~\ref{cor:num} extends this to arbitrary $\mu \in P^+$.
\end{rem}

The following result is a higher-level analogue of the main result of~\cite{KL17}:

\begin{cor}\label{cor:int-filt}
Let $k, l \in \Z_{>0}$ be such that $k \ge l - 1$. For each $\La \in P^+_{\af,k}$, the module ${}^\theta L ( \La )$ admits a $\bW^{(l)}$-filtration.
\end{cor}

\begin{proof}
This follows by combining Lemma~\ref{lem:idDL} with Theorem~\ref{thm:branch}(1).
\end{proof}

\begin{cor}\label{cor:Ww}
Let $k, l \in \Z_{>0}$. For each $\la \in P^+$, the module $\bW_{\la}^{(k)}$ admits a $W^{(l)}$-filtration.
\end{cor}

\begin{proof}
By Corollary~\ref{cor:passage}, we have
\[
\ext^1_{\fC} \left( \bW_{\la}^{(k)}, ( W_{\mu}^{(1)} \otimes \C_{(k-1)\La_0} )^* \right) = 0 \qquad \text{for all } \la, \mu \in P^+.
\]
In particular, this implies that each $W^{(1)}_\mu$ admits a $W^{(k)}$-filtration by Theorem~\ref{thm:critfilt}(3) applied to $M = W_{\mu}^{(1)} \otimes \C_{(k-1)\La_0}$.

Moreover, Theorem~\ref{thm:critfilt}(4) applied to $M = \bW_{\la}^{(k)}$ shows that $\bW^{(k)}_\la$ admits a $\bW^{(1)}$-filtration. Since $\bW^{(1)}_\la$ is a self-extension of $W^{(1)}_\la$ by Theorem~\ref{thm:CKmain1}, it in particular admits a decreasing separable $W^{(1)}$-filtration. Therefore, we conclude that $\bW^{(k)}_\la$ admits a $W^{(l)}$-filtration.
\end{proof}

\appendix
\renewcommand{\thesection}{\Alph{section}}
\setcounter{section}{1}
\renewcommand{\theequation}{\Alph{section}.\arabic{equation}}
\setcounter{equation}{0}
\setcounter{thm}{0}

\addcontentsline{toc}{section}{Appendix: Higher level $q$-Cauchy kernel identities}
\begin{flushleft}
{\normalsize\textbf{Appendix: Higher level $q$-Cauchy kernel identities}}
\end{flushleft}

In this appendix, we derive several numerical identities from the main results of the paper. We adopt the same conventions as in the main body (see \S\ref{sec:module}).

In addition, we define $\gch_x$ and $\gch_y$ as variants of $\gch$ in which the $\h$-character $e^\lambda$ is recorded independently as $x^\lambda$ (for $\gch_x$) and $y^\lambda$ (for $\gch_y$), while the $q$-grading is recorded in the same way in both cases.

Throughout this appendix, all equalities are understood in the sense of formal series.

\begin{prop}\label{prop:appprop}
For each $k \in \Z_{>0}$, we have the identity
\[
\sum_{\la \in P^+} ( \gch\nolimits_x V_\la^* ) \cdot \gch\nolimits_y P_\la
= \sum_{\la \in P^+} ( \gch\nolimits_x W_{-w_0\la}^{(k)} ) \cdot \gch\nolimits_y \bW_\la^{(k)}.
\]
\end{prop}

\begin{proof}
By Remark~\ref{rem:limit}, for each fixed $\la \in P^+$, the characters
\[
\gch\nolimits_x W_{-w_0\la}^{(k)}
\quad\text{and}\quad
\gch\nolimits_y \bW_\la^{(k)}
\]
stabilize coefficientwise to $\gch\nolimits_x V_\la^*$ and $\gch\nolimits_y P_\la$, respectively, as $k \to \infty$. Hence it suffices to show that the right-hand side is independent of $k$.

Using~\eqref{eqn:hBGG-rec}, we compute
\begin{align*}
\sum_{\la \in P^+} \gch\nolimits_x W_{-w_0\la}^{(k)} \cdot \gch\nolimits_y \bW_\la^{(k)}
&= \sum_{\la, \mu \in P^+}
(W_{-w_0\la}^{(k)} : W_{-w_0\mu}^{(k+1)})_q \cdot
\gch\nolimits_x W_{-w_0\mu}^{(k+1)} \cdot \gch\nolimits_y \bW_\la^{(k)} \\
&= \sum_{\la, \mu \in P^+}
(\bW_{-w_0\mu}^{(k+1)} : \bW_{-w_0\la}^{(k)})_q \cdot
\gch\nolimits_x W_{-w_0\mu}^{(k+1)} \cdot \gch\nolimits_y \bW_\la^{(k)} \\
&= \sum_{\la, \mu \in P^+}
(\bW_{\mu}^{(k+1)} : \bW_{\la}^{(k)})_q \cdot
\gch\nolimits_x W_{-w_0\mu}^{(k+1)} \cdot \gch\nolimits_y \bW_\la^{(k)} \\
&= \sum_{\mu \in P^+}
\gch\nolimits_x W_{-w_0\mu}^{(k+1)} \cdot \gch\nolimits_y \bW_\mu^{(k+1)},
\end{align*}
where the last equality follows from the expansion of $\gch \bW_\mu^{(k+1)}$ in the basis $\{\gch \bW_\la^{(k)}\}_{\la \in P^+}$.
Thus the right-hand side is independent of $k$, as desired.
\end{proof}

\begin{prop}\label{prop:apppropb}
For each $k \in \Z_{>0}$, we have the identity
\[
\sum_{\la \in P} \ch_{x} \C_{-\la} \cdot \gch\nolimits_y Q_\la
= \sum_{\la \in P} ( \gch\nolimits_x D_{-\la}^{(k)} ) \cdot \gch\nolimits_y \bD_\la^{(k)}.
\]
\end{prop}

\begin{proof}
By Theorem~\ref{thm:branch-b}, an argument analogous to the proof of Proposition~\ref{prop:appprop} shows that
\[
\sum_{\la \in P} ( \gch\nolimits_x D_{-\la}^{(k+1)} ) \cdot \gch\nolimits_y \bD_\la^{(k+1)}
= \sum_{\la \in P} ( \gch\nolimits_x D_{-\la}^{(k)} ) \cdot \gch\nolimits_y \bD_\la^{(k)}.
\]
Hence it is enough to identify this common value in the limit $k \to \infty$.

The degree-zero part with respect to the $d$-grading is shared by all these expressions and corresponds to the characters of
\[
D_\la := U(\gb) \cdot \bv_\la \subset V_\la, \qquad
\mathring{D}_\la := \frac{D_\la}{\sum_{\mu \in W\la,\; D_\la \supsetneq D_\mu} D_\mu}.
\]
We have
\[
\lim_{k \to \infty} \gch\nolimits_x D_\la^{(k)} = \gch\nolimits_x D_\la,
\]
since $(\la + k \La_0) \in W P_\af^+$ for $k \gg 0$ and fixed $\la \in P$.

Similarly, for each fixed $d$-degree, the quotient defining $\bD_\la^{(k)}$ eventually imposes no further relation in that degree as $k \to \infty$. Hence
\[
\lim_{k \to \infty} \gch\nolimits_y \bD_\la^{(k)}
= \gch\nolimits_y \left( U(\tb) \otimes_{U(\gb + \wth)} \mathring{D}_\la \right).
\]
It follows that
\[
\gch\nolimits_y U ( \tb / (\gn + \wth) ) \cdot \gch\nolimits_y \mathring{D}_\la
= \lim_{k \to \infty} \gch\nolimits_y \bD_\la^{(k)} \qquad (\la \in P),
\]
and hence
\[
\gch\nolimits_y U ( \tb / (\gn + \wth) )
\cdot \sum_{\la \in P} ( \gch\nolimits_x D_{-\la} ) \cdot \gch\nolimits_y \mathring{D}_\la
= \sum_{\la \in P} ( \gch\nolimits_x D_{-\la}^{(k)} ) \cdot \gch\nolimits_y \bD_\la^{(k)}.
\]

By~\cite[Theorem~1.13]{vdK89}, we have
\[
\sum_{\la \in P} ( \gch\nolimits_x D_{-\la} ) \cdot \gch\nolimits_y \mathring{D}_\la
= \ch_{x,y} \C[B]^{\vee},
\]
where $\ch_{x,y}$ denotes the character as an $(\h \oplus \h)$-module. On the other hand,
\[
\ch_{x,y} \C[B]^{\vee}
= \sum_{\la \in P} \ch_x \C_{-\la} \cdot \ch_y \left( U(\gb) \otimes_{U(\h)} \C_\la \right)
\]
by direct computation.

Combining these identities, we obtain
\begin{align*}
\sum_{\la \in P} ( \gch\nolimits_x D_{-\la}^{(k)} ) \cdot \gch\nolimits_y \bD_\la^{(k)}
&= \gch\nolimits_y U(\tb / (\gn + \wth)) \cdot
\sum_{\la \in P} ( \gch\nolimits_x D_{-\la} ) \cdot \gch\nolimits_y \mathring{D}_\la \\
&= \gch\nolimits_y U(\tb / (\gn + \wth)) \cdot
\sum_{\la \in P} \ch_x \C_{-\la} \cdot \ch_y ( U(\gn) \otimes \C_\la ) \\
&= \sum_{\la \in P} \ch_x \C_{-\la} \cdot \gch\nolimits_y Q_\la,
\end{align*}
as required.
\end{proof}

Let $G$ be the connected, simply connected, simple algebraic group over $\C$ whose Lie algebra is $\g$, and let $N$ and $H$ denote the connected closed subgroups of $G$ whose Lie algebras are $\gn$ and $\h$, respectively.

Let $\widetilde{G}$ denote the (pro-)algebraic group whose Lie algebra is $[\tg_{\ge 0}, \tg_{\ge 0}]$, and let $\C[\widetilde{G}]$ be its coordinate ring. Let $\widetilde{B}$ denote its closed subgroup scheme with Lie algebra $\tb \cap [\tg_{\ge 0}, \tg_{\ge 0}]$, and let $\C[\widetilde{B}]$ be its coordinate ring.

Both $\widetilde{G}$ and $\widetilde{B}$ admit actions of $H \times H$ as well as of $\C^\times$, the latter recording the grading. Accordingly, we define the graded character $\gch_{x,y}$, which records the left $H$-action via $x^\bullet$, the right $H$-action via $y^\bullet$, and the common grading via the variable $q$.

The restricted dual $(\bullet)^\vee$ in Theorem~\ref{thm:appmain} is understood as the restricted dual taken with respect to the structure of an $(H \times H \times \Gm)$-module. Note that this differs from the restricted dual considered in~\S\ref{rep-current}.

\begin{thm}\label{thm:appmain}
For each $k \in \Z_{>0}$, we have the identities
\begin{align*}
\gch\nolimits_{x,y} \, \C[\widetilde{G}]^\vee
&= \sum_{\la \in P^+} ( \gch\nolimits_x W_{-w_0\la}^{(k)} ) \cdot \gch\nolimits_y \bW_\la^{(k)}, \\[2pt]
\gch\nolimits_{x,y} \, \C[\widetilde{B}]^\vee
&= \sum_{\la \in P} ( \gch\nolimits_x D_{-\la}^{(k)} ) \cdot \gch\nolimits_y \bD_\la^{(k)}.
\end{align*}
\end{thm}

\begin{rem}
When $k = 1$, Theorem~\ref{thm:appmain} recovers the identities established in~\cite{FKM19} and~\cite{FMO23,FKMO}, respectively. As explained in Remark~\ref{rem:theta}, an identity involving theta functions arises when $\g = \mathfrak{sl}(2)$ and $k = 2$.
\end{rem}

\begin{proof}[Proof of Theorem~\ref{thm:appmain}]
We set
\[
\widetilde{G}_+ := \ker\bigl( \widetilde{G} \stackrel{\mathtt{ev}_0}{\longrightarrow} G \bigr),
\]
where $\mathtt{ev}_0$ denotes the evaluation map at $z = 0$. Then
\[
\C[\widetilde{G}] \cong \C[G] \otimes \C[\widetilde{G}_+],
\qquad
\C[\widetilde{B}] \cong \C[H \cdot N] \otimes \C[\widetilde{G}_+].
\]

By the algebraic Peter--Weyl theorem, we have
\[
\C[\widetilde{G}] \cong \bigoplus_{\la \in P^+} V_\la \boxtimes \bigl( V_\la^* \otimes \C[\widetilde{G}_+] \bigr).
\]
Moreover, for each $\la \in P^+$, we have an isomorphism of right $\tg_{\ge 0}$-modules
\[
P_\la \xrightarrow{\sim} V_\la \otimes \C[\widetilde{G}_+]^\vee,
\]
compatible with the left $G$-action. Taking characters as $(H \times H \times \Gm)$-modules, we obtain
\[
\gch\nolimits_{x,y}\, \C[\widetilde{G}]^\vee
= \sum_{\la \in P^+} \bigl( \gch\nolimits_x V_\la^* \bigr) \cdot \gch\nolimits_y P_\la.
\]
Thus, Proposition~\ref{prop:appprop} yields the first identity.

Similarly, we have
\[
\C[\widetilde{B}]^\vee \cong \bigoplus_{\la \in P}
\C_{-\la} \boxtimes \bigl( \C_\la \otimes \C[N \cdot \widetilde{G}_+]^\vee \bigr),
\]
and, for each $\la \in P$, an isomorphism of right $\tb$-modules
\[
Q_\la \xrightarrow{\sim} \C_\la \otimes \C[N \cdot \widetilde{G}_+]^\vee,
\]
compatible with the left $H$-action. Taking characters again as $(H \times H \times \Gm)$-modules, we obtain
\[
\gch\nolimits_{x,y}\, \C[\widetilde{B}]^\vee
= \sum_{\la \in P} \bigl( \gch\nolimits_x \C_{-\la} \bigr) \cdot \gch\nolimits_y Q_\la.
\]
Hence, Proposition~\ref{prop:apppropb} yields the second identity.
\end{proof}

\medskip

{\small
{\bf Acknowledgments:}
This work was supported in part by JSPS KAKENHI Grant Number JP19H01782 and JP24K21192. The author is indebted to Rekha Biswal for suggesting Corollary~\ref{cor:KR} and Remark~\ref{rem:kostka}(2), and to Masato Okado for informing
him about one-dimensional sums and drawing his attention to~\cite{Oka}. The author would also like to thank Shrawan Kumar for helpful email correspondence. Finally, the author is indebted to the referee who has made many constructive suggestions throughout this paper.}

\medskip
{\small
{\bf Conflict of interest:} None.}

\medskip

{\small
{\bf Data availability statement:} My manuscript has no associated data.}

{\footnotesize
\bibliography{ref}

\begin{thebibliography}{10}

\bibitem{BBCKL}
Matthew Bennett, Arkady Berenstein, Vyjayanthi Chari, Anton Khoroshkin, and
  Sergey Loktev.
\newblock Macdonald polynomials and {BGG} reciprocity for current algebras.
\newblock {\em Selecta Math. (N.S.)}, 20(2):585--607, 2014.

\bibitem{BCM}
Matthew Bennett, Vyjayanthi Chari, and Nathan Manning.
\newblock B{GG} reciprocity for current algebras.
\newblock {\em Adv. Math.}, 231(1):276--305, 2012.

\bibitem{BK21}
Rekha Biswal and Deniz Kus.
\newblock A combinatorial formula for graded multiplicities in excellent
  filtrations.
\newblock {\em Transform. Groups}, 26(1):81--114, 2021.

\bibitem{BS20}
Rekha Biswal and Travis Scrimshaw.
\newblock Existence of {K}irillov-{R}eshetikhin crystals for multiplicity-free
  nodes.
\newblock {\em Publ. Res. Inst. Math. Sci.}, 56(4):761--778, 2020.

\bibitem{BM95}
Maegan~K. Bos and Kailash~C. Misra.
\newblock An application of crystal bases to representations of affine {L}ie
  algebras.
\newblock {\em Journal of Algebra}, 173:436--458, 1995.

\bibitem{Bou02}
Nicolas Bourbaki.
\newblock {\em Elements of Mathematics. Lie Groups and Lie algebras. Chapters
  4--6}.
\newblock Springer-Verlag, Berlin, 2002.

\bibitem{Bou08}
Nicolas Bourbaki.
\newblock {\em Lie Groups and Lie Algebras: Chapters 7--9}.
\newblock Elements of Mathematics. Springer-Verlag, Berlin, 2008.
\newblock Softcover reprint of the 2005 English translation.

\bibitem{BS24}
Jonathan Brundan and Catharina Stroppel.
\newblock Semi-infinite highest weight categories.
\newblock {\em Mem. Amer. Math. Soc.}, 293(1459):vii+152, 2024.

\bibitem{CE56}
Henri Cartan and Samuel Eilenberg.
\newblock {\em Homological Algebra}.
\newblock Princeton University Press, Princeton, NJ, 13th printing, first
  paperback printing edition, 1956.
\newblock Reprinted 1999.

\bibitem{CFK}
Vyjayanthi Chari, Ghislain Fourier, and Tanusree Khandai.
\newblock A categorical approach to {W}eyl modules.
\newblock {\em Transform. Groups}, 15:517--549, 2010.

\bibitem{CG07}
Vyjayanthi Chari and Jacob Greenstein.
\newblock Current algebras, highest weight categories and quivers.
\newblock {\em Adv. Math.}, 216(2):811--840, 2007.

\bibitem{CI15}
Vyjayanthi Chari and Bogdan Ion.
\newblock B{GG} reciprocity for current algebras.
\newblock {\em Compos. Math.}, 151(7):1265--1287, 2015.

\bibitem{CL06}
Vyjayanthi Chari and Sergei Loktev.
\newblock Weyl, {D}emazure and fusion modules for the current algebra of
  {$\mathfrak{sl}_{r+1}$}.
\newblock {\em Adv. Math.}, 207(2):928--960, 2006.

\bibitem{CP01}
Vyjayanthi Chari and Andrew Pressley.
\newblock Weyl modules for classical and quantum affine algebras.
\newblock {\em Represent. Theory}, 5:191--223 (electronic), 2001.

\bibitem{CSSW}
Vyjayanthi Chari, Lisa Schneider, Peri Shereen, and Jeffrey Wand.
\newblock Modules with {D}emazure flags and character formulae.
\newblock {\em SIGMA Symmetry Integrability Geom. Methods Appl.}, 10(032):16pp,
  2014.

\bibitem{Che95}
Ivan Cherednik.
\newblock {N}onsymmetric {M}acdonald {P}olynomials.
\newblock {\em International Mathematics Research Notices}, 2(10), 1995.

\bibitem{CF13}
Ivan Cherednik and Boris Feigin.
\newblock Rogers-{R}amanujan type identities and {N}il-{DAHA}.
\newblock {\em Adv. Math.}, 248:1050--1088, 2013.

\bibitem{CK18}
Ivan Cherednik and Syu Kato.
\newblock Nonsymmetric {R}ogers-{R}amanujan sums and thick {D}emazure modules.
\newblock {\em Adv. in Math.}, 374:Article number 107335, 2020.

\bibitem{Chi22}
Masahiro Chihara.
\newblock Demazure slices of type {$\mathsf A^{(2)}_{2\ell}$}.
\newblock {\em Algebras and Representation Theory}, 25:491--519, 2022.

\bibitem{CPS88}
E.~Cline, B.~Parshall, and L.~Scott.
\newblock Finite dimensional algebras and highest weight categories.
\newblock {\em J. reine angew. Math.}, 391:85--99, 1988.

\bibitem{FKM}
Evgeny Feigin, Syu Kato, and Ievgen Makedonskyi.
\newblock Representation theoretic realization of non-symmetric {M}acdonald
  polynomials at infinity.
\newblock {\em J. reine angew. Math.}, 767:181--216, 2020.
\newblock arXiv: 1703.04108.

\bibitem{FKM19}
Evgeny Feigin, Anton Khoroshkin, and Ievgen Makedonskyi.
\newblock Duality theorems for current groups.
\newblock {\em Israel Journal of Mathematics}, 248:441--479, 2022.

\bibitem{FKMO}
Evgeny Feigin, Anton Khoroshkin, Ievgen Makedonskyi, and Daniel Orr.
\newblock Peter-{W}eyl theorem for {I}wahori groups and highest weight
  categories.
\newblock {\em Transformation Groups}, to appear.

\bibitem{FMO18}
Evgeny Feigin, Ievgen Makedonskyi, and Daniel Orr.
\newblock Generalized {W}eyl modules and nonsymmetric $q$-{W}hittaker
  functions.
\newblock {\em Adv. Math.}, 339:997--1033, 2018.

\bibitem{FMO23}
Evgeny Feigin, Ievgen Makedonskyi, and Daniel Orr.
\newblock Nonsymmetric $q$-{C}auchy identity and representations of the
  {I}wahori algebra.
\newblock {\em Kyoto J. Math.}, to appear.

\bibitem{FL07}
G.~Fourier and P.~Littelmann.
\newblock Weyl modules, {D}emazure modules, {KR}-modules, crystals, fusion
  products and limit constructions.
\newblock {\em Adv. Math.}, 211(2):566--593, 2007.

\bibitem{FMS13}
G.~Fourier, N.~Manning, and P.~Senesi.
\newblock Global {W}eyl modules for the twisted loop algebra.
\newblock {\em Abh. Math. Semin. Univ. Hamb.}, 83:53--82, 2013.

\bibitem{FSS07}
Ghislain Fourier, Anne Schilling, and Mark Shimozono.
\newblock Demazure structure inside {K}irillov-{R}eshetikhin crystals.
\newblock {\em J. Algebra}, 309(1):386--404, 2007.

\bibitem{FK80}
Igor Frenkel and Victor Kac.
\newblock Basic representations of affine {L}ie algebras and dual resonance
  models.
\newblock {\em Invent. Math.}, 62:23--66, 1980/1981.

\bibitem{Gro57}
Alexander Grothendieck.
\newblock Sur quelques points d'alg\`ebre homologique.
\newblock {\em Tohoku Math. J. (2)}, 9:119--221, 1957.

\bibitem{HKOTY}
G.~Hatayama, A.~Kuniba, M.~Okado, T.~Takagi, and Y.~Yamada.
\newblock Remarks on fermionic formula.
\newblock In {\em Recent developments in quantum affine algebras and related
  topics ({R}aleigh, {NC}, 1998)}, volume 248 of {\em Contemp. Math.}, pages
  243--291. Amer. Math. Soc., Providence, RI, 1999.

\bibitem{HKKOTY}
Goro Hatayama, Anatol~N. Kirillov, Atsuo Kuniba, Masato Okado, Taichiro Takagi,
  and Yasuhiko Yamada.
\newblock Character formulae of {$\widehat{\rm sl}_n$}-modules and
  inhomogeneous paths.
\newblock {\em Nuclear Phys. B}, 536(3):575--616, 1999.

\bibitem{HKOTT}
Goro Hatayama, Atsuo Kuniba, Masato Okado, Taichiro Takagi, and Zengo Tsuboi.
\newblock Paths, crystals and fermionic formulae.
\newblock In {\em Math{P}hys odyssey, 2001}, volume~23 of {\em Prog. Math.
  Phys.}, pages 205--272. Birkh\"auser Boston, Boston, MA, 2002.

\bibitem{Hec87}
G.~J. Heckman.
\newblock Root systems and hypergeometric functions. {II}.
\newblock {\em Compositio Mathematica}, 64(3):353--373, 1987.

\bibitem{Ion03}
Bogdan Ion.
\newblock Nonsymmetric {M}acdonald polynomials and {D}emazure characters.
\newblock {\em Duke Math. J.}, 116(2):299--318, 2003.

\bibitem{Jos85}
Anthony Joseph.
\newblock On the {D}emazure character formula.
\newblock {\em Ann. Sci. \'Ec. Norm. Sup\'er. (4)}, 18(3):389--419, 1985.

\bibitem{Jos03}
Anthony Joseph.
\newblock A decomposition theorem for {D}emazure crystals.
\newblock {\em J. Algebra}, 265(2):562--578, 2003.

\bibitem{Jos06}
Anthony Joseph.
\newblock Modules with a {D}emazure flag.
\newblock In {\em Studies in Lie theory}, volume 243 of {\em Progress in
  Mathematics}. Birkh\"auser Boston, Boston, MA, 2006.

\bibitem{Kac}
Victor~G. Kac.
\newblock {\em Infinite-dimensional {L}ie algebras}.
\newblock Cambridge University Press, Cambridge, third edition, 1990.

\bibitem{KKMMNN}
Seok-Jin Kang, Masaki Kashiwara, Kailash~C. Misra, Tetsuji Miwa, Toshiki
  Nakashima, and Atsushi Nakayashiki.
\newblock Perfect crystals of quantum affine {L}ie algebras.
\newblock {\em Duke Mathematical Journal}, 68(3):499--607, 1992.

\bibitem{KMS95}
M.~Kashiwara, T.~Miwa, and E.~Stern.
\newblock Decomposition of $q$-deformed {F}ock spaces.
\newblock {\em Selecta Mathematica}, 1(4):787--805, 1995.
\newblock Published: 01 December 1995.

\bibitem{Kas93}
Masaki Kashiwara.
\newblock The crystal base and {L}ittelmann's refined {D}emazure character
  formula.
\newblock {\em Duke Math. J.}, 71(3):839--858, 1993.

\bibitem{Kas05}
Masaki Kashiwara.
\newblock Level zero fundamental representations over quantized affine algebras
  and {D}emazure modules.
\newblock {\em Publ. Res. Inst. Math. Sci.}, 41(1):223--250, 2005.

\bibitem{KS09}
Masaki Kashiwara and Mark Shimozono.
\newblock Equivariant {$K$}-theory of affine flag manifolds and affine
  {G}rothendieck polynomials.
\newblock {\em Duke Math. J.}, 148(3):501--538, 2009.

\bibitem{Kat18}
Syu Kato.
\newblock Demazure character formula for semi-infinite flag varieties.
\newblock {\em Math. Ann.}, 371(3):1769--1801, 2018.

\bibitem{Kat18b}
Syu Kato.
\newblock Frobenius splitting of thick flag manifolds of {K}ac-{M}oody
  algebras.
\newblock {\em Int. Math. Res. Not. IMRN}, 2020(17):5401--5427, 2020.
\newblock arXiv:1707.03773.

\bibitem{KL17}
Syu Kato and Sergey Loktev.
\newblock A {W}eyl module stratification of integrable representations.
\newblock {\em Comm. Math. Phys.}, 368:113--141, 2019.
\newblock arXiv:1712.03508.

\bibitem{Kum02}
Shrawan Kumar.
\newblock {\em Kac-{M}oody groups, their flag varieties and representation
  theory}, volume 204 of {\em Progress in Mathematics}.
\newblock Birkh\"auser Boston, Inc., Boston, MA, 2002.

\bibitem{Mac95}
I.~G. Macdonald.
\newblock {\em Symmetric functions and {H}all polynomials}.
\newblock Oxford Mathematical Monographs. The Clarendon Press, Oxford
  University Press, New York, second edition, 1995.
\newblock With contributions by A. Zelevinsky, Oxford Science Publications.

\bibitem{Mac95B}
Ian~G. Macdonald.
\newblock Affine {H}ecke algebras and orthogonal polynomials.
\newblock In {\em S{\'e}minaire {B}ourbaki, {V}ol. 1994/95}, number 797, pages
  1--18. Soci{\'e}t{\'e} Math{\'e}matique de France, Paris, 1995.
\newblock Ast{\'e}risque No. 237 (1996), Exp. No. 797, 3.

\bibitem{Mac03}
Ian~G. Macdonald.
\newblock {\em Affine Hecke Algebras and Orthogonal Polynomials}, volume 128 of
  {\em Cambridge Tracts in Mathematics}.
\newblock Cambridge University Press, 2003.

\bibitem{Mat88}
Olivier Mathieu.
\newblock Formules de caract\`eres pour les alg\`ebres de {K}ac-{M}oody
  g\'en\'erales.
\newblock {\em Ast\'erisque}, 159--160:1--267, 1988.

\bibitem{Mat89c}
Olivier Mathieu.
\newblock Frobenius action on the {$B$}-cohomology.
\newblock In {\em Infinite-dimensional {L}ie algebras and groups
  ({L}uminy-{M}arseille, 1988)}, volume~7 of {\em Adv. Ser. Math. Phys.} World
  Sci. Publ., Teaneck, NJ, 1989.

\bibitem{Nao12}
Katsuyuki Naoi.
\newblock Weyl modules, {D}emazure modules and finite crystals for non-simply
  laced type.
\newblock {\em Adv. in Math.}, 229(2):875--934, 2012.

\bibitem{Nao13}
Katsuyuki Naoi.
\newblock Demazure crystals and tensor products of perfect
  {K}irillov-{R}eshetikhin crystals with various levels.
\newblock {\em J. Algebra}, 374:1--26, 2013.

\bibitem{Nao18}
Katsuyuki Naoi.
\newblock Existence of {K}irillov-{R}eshetikhin crystals of type {$G_2^{(1)}$}
  and {$D_4^{(3)}$}.
\newblock {\em J. Algebra}, 512:47--65, 2018.

\bibitem{NS21}
Katsuyuki Naoi and Travis Scrimshaw.
\newblock Existence of {K}irillov-{R}eshetikhin crystals for near adjoint nodes
  in exceptional types.
\newblock {\em J. Pure Appl. Algebra}, 225(5):Paper No. 106593, 38, 2021.

\bibitem{Oka}
Masato Okado.
\newblock mimeo.

\bibitem{Pol89}
Patrick Polo.
\newblock Vari\'et\'es de {S}chubert et excellentes filtrations.
\newblock In {\em Orbites unipotentes et repr\'esentations, III}, number
  173-174 in Ast\'erisque. 1989.

\bibitem{San00}
Yasmine~B. Sanderson.
\newblock On the connection between {M}acdonald polynomials and {D}emazure
  characters.
\newblock {\em J. Algebraic Combin.}, 11(3):269--275, 2000.

\bibitem{SW99}
A.~Schilling and S.~O. Warnaar.
\newblock Inhomogeneous lattice paths, generalized {K}ostka polynomials and
  $\mathsf{A}_{n-1}$-supernomials.
\newblock {\em Commun. Math. Phys}, 202:359--401, 1999.

\bibitem{SS01}
Anne Schilling and Mark Shimozono.
\newblock Fermionic formulas for level-restricted generalized {K}ostka
  polynomials and coset branching functions.
\newblock {\em Comm. Math. Phys.}, 220(1):105--164, 2001.

\bibitem{ST12}
Anne Schilling and Peter Tingley.
\newblock Demazure crystals, {K}irillov-{R}eshetikhin crystals, and the energy
  function.
\newblock {\em Electron. J. Combin.}, 19(2):Paper 4, 42, 2012.

\bibitem{Seg81}
Graeme Segal.
\newblock Unitary representations of some infinite-dimensional groups.
\newblock {\em Comm. Math. Phys.}, 80(3):301--342, 1981.

\bibitem{vdK89}
Wilberd van~der Kallen.
\newblock Longest weight vectors and excellent filtrations.
\newblock {\em Math. Zeit.}, 201(no.1):19--31, 1989.

\bibitem{Ver96}
Jean-Louis Verdier.
\newblock Des cat{\'e}gories d{\'e}riv{\'e}es des cat{\'e}gories
  ab{\'e}liennes.
\newblock {\em Ast{\'e}risque}, 239, 1996.

\bibitem{Wei94}
Charles~A. Weibel.
\newblock {\em An Introduction to Homological Algebra}, volume~38 of {\em
  Cambridge Studies in Advanced Mathematics}.
\newblock Cambridge University Press, Cambridge, 1994.

\end{thebibliography}
\bibliographystyle{hplain}}
\end{document}